\documentclass[aps,pra,superscriptaddress,nofootinbib]{revtex4}

\pdfoutput=1
\usepackage[intlimits]{amsmath}
\usepackage{amsfonts}

\usepackage{fancyvrb}
\usepackage{enumitem}
\usepackage{psfrag}

\usepackage[utf8]{inputenc}
\usepackage{longtable}

\usepackage{amsthm}
\usepackage{float}
\usepackage{tikz}

\makeatletter
\newcommand\Row[1]{%
  \par\nobreak\nointerlineskip\vskip-\fboxrule%
  \@tfor\@tempa:=#1 \do {\csname ChessBox\@tempa\endcsname\kern-\fboxrule}}
\define@key{chessB}{side}{\def\ChessSide{#1}}
\define@key{chessB}{colori}{\def\ChessColori{#1}}
\define@key{chessB}{colorii}{\def\ChessColorii{#1}}
\setkeys{chessB}{
  side=1.5em,
  colori=black!70,
  colorii=white}
\makeatother

\newenvironment{Chessboard}[1][]
  {\setkeys{chessB}{#1}%
  \par\medskip\setlength\parindent{0pt}}
  {\par\medskip}

\advance\voffset by -1cm

\usepackage[utf8]{inputenc}

\newtheorem{theorem}{Theorem}

\newcommand\beq            {\begin{equation}}
\newcommand\be            {\begin{equation}}
\newcommand\bea           {\begin{equation}\begin{array}l\displaystyle}
\newcommand\ee            {\end{equation}}
\newcommand\eeq            {\end{equation}}
\newcommand\bes           {\begin{subequations}}
\newcommand\esu           {\end{subequations}}

\renewcommand{\(}{\left(}
\renewcommand{\)}{\right)}
\renewcommand{\[}{\left[}
\renewcommand{\]}{\right]}

\counterwithout{table}{section}

\newcommand{\bigx}[1]{\bBigg@{#1}}

\def\3pt#1#2#3{{\langle{#1}\vert{#2}\vert{#3}\rangle}}

\newcommand\doi[2]        {\href{http://dx.doi.org/#1}{#2}}

\newcommand{\EQ}{\begin{equation}}
\newcommand{\EN}{\end{equation}}
\usepackage{epsfig}
\usepackage{psfrag}
\usepackage{amsmath}
\usepackage{graphicx}
\usepackage{amsfonts}
\usepackage{amssymb}

\def\tilde{\widetilde}

\def\hat{\widehat}
\def\*{\star}
\def\[{\left[}
\def\]{\right]}
\def\({\left(}      

\def\){\right)}

\def\frac#1#2{\dfrac{#1}{#2}}
\def\inv#1{\dfrac{1}{#1}}
\def\half{\tfrac{1}{2}}

\def\2pi{\hbox{$2\pi i$}}

\def\dsl{\raise.15ex\hbox{/}\kern-.57em\partial}
\def\Dsl{\,\raise.15ex\hbox{/}\mkern-.13.5mu D}
\def\be{\beta}

   \def\CE{{\cal E}}

   \def\CN{{\cal N}}

\def\2pi{\hbox{$2\pi i$}}

\def\dsl{\raise.15ex\hbox{/}\kern-.57em\partial}
\def\Dsl{\,\raise.15ex\hbox{/}\mkern-.13.5mu D}

\def\barray{\begin{eqnarray}}
\def\earray{\end{eqnarray}}
\def\beq{\begin{equation}}
\def\eeq{\end{equation}}

\def\AA{\leavevmode\setbox0=\hbox{h}
\dimen0=\ht0 \advance\dimen0 by-1ex\rlap{\raise.67\dimen0\hbox{\char'27}}A}

\def\blue#1{{\color{blue}{#1}}}

\def\blue#1{{\color{blue}{#1}}}

\def\half{\tfrac{1}{2}}

\def\prob{{\bf Pr}}

\def\probbracL{\Bigl[}
\def\probbracR{\Bigr]}

\def\U{M}

\def\Primes{\mathbb{P}}
\def\Primesp{\mathbb{P}'}

\def\sqf{\mathfrak{f}}

\def\half{\tfrac{1}{2}}

\def\s{\ell} 

\def\muhat{\hat{\mu}}

\def\limsupn{{\underset{n\to \infty} {\rm lim ~ sup}}}

\def\limsupL{{\underset{L\to \infty} {\rm lim ~ sup}}}

\def\liminfn{{\underset{n\to \infty} {\rm lim ~ inf}}}

\begin{document}
\bibliographystyle{plainnat}

\title{{\Large {\bf Randomness of  M\"obius coefficents and Brownian Motion: \\
 Growth of the  Mertens  Function and the Riemann Hypothesis}
 }} 
\author{Giuseppe Mussardo}
\affiliation{SISSA and INFN, Sezione di Trieste, via Bonomea 265, I-34136, 
Trieste, Italy}
\author{ Andr\'e  LeClair}
\affiliation{Cornell University, Physics Department, Ithaca, NY 14850} 


\begin{abstract}
\noindent
\noindent
The validity of the Riemann Hypothesis (RH) on the location of the non-trivial zeros of the Riemann $\zeta$-function is directly related to the growth of the Mertens function 
$M(x) \,=\,\sum_{k=1}^x \mu(k)$, where $\mu(k)$ is the M\"{o}bius coefficient of the integer $k$: the RH is indeed true if the Mertens function goes asymptotically as 
 $M(x) \sim x^{1/2 + \epsilon}$, where $\epsilon$ is an arbitrary strictly positive quantity. We {argue} that this behavior can be established on the basis of a new probabilistic approach based on the {\em global} properties of Mertens function, namely, {based on reorganizing globally in distinct blocks the terms of its series.}    
To this aim, we focus the attention on the square-free numbers and we derive a series of probabilistic results concerning the prime number distribution along the series of square-free numbers, the average number of prime divisors, the Erd\H{o}s-Kac theorem for square-free numbers, etc. These results {point to} the conclusion that the Mertens function is subject to a normal distribution as much as any other random walk. {We also present an argument in favor of the thesis that the validity of the Riemann Hypothesis also implies the validity of the Generalized Riemann Hypothesis for the Dirichlet $L$-functions. }
Next we study the {\em local} properties of the Mertens function, {i.e. its variation induced by each M\"{o}bius coefficient restricted to the 
square-free numbers.} Motivated  by the natural curiosity to see how close to a purely random walk is any sub-sequence extracted by the sequence of the M\"{o}bius coefficients for the square-free numbers, we perform a massive statistical analysis on these coefficients, applying to them a series of randomness tests of increasing precision and complexity: together with several frequency tests within a block, the list of our tests include those for the longest run of ones in a block, the binary matrix rank test, the Discrete Fourier Transform test, the non-overlapping template matching test, the entropy test, the cumulative sum test, the random excursion tests, etc. {for a total number of eighteen different tests}. 
The successful outputs of all these tests  (each of them with a level of confidence  of $99\%$  that all the sub-sequences analyzed are indeed random) can be seen as impressive ``experimental" confirmations of the brownian nature of the restricted M\"{o}bius coefficients and the probabilistic normal law distribution of the Mertens function analytically established earlier. In view of the theoretical probabilistic argument and the large battery of statistical tests, we can conclude that while a violation of the RH is strictly speaking not impossible, it is however extremely improbable. 
\end{abstract}

\maketitle

.

\newpage 

. 

\newpage 
\tableofcontents

\newpage
\vspace{5mm}

\section{Introduction}

\noindent

\noindent
In Physics there are innumerable examples of a successful scientific strategy which consists of the following steps:  
\begin{description}
\item{A} Start from a consolidated set of empirical data.
\item{B} Elaborate a theory which explains the data: besides its elegance and beauty, the strength of the proposed theory 
is related to the diversity of phenomena it can explain, together with its simplicity and parsimony.
\item{C} Make repeated testable and thorough attempts to falsify the theory. When theories are falsified by such observations, 
scientists can respond by revising the theory, or by rejecting the theory in favor of  another one. In either case, however, 
this process must aim at the production of new, falsifiable predictions. 
\end{description}
Elaborating a theory is of course something different than proving a theorem, even though a theory may be based on theorems or may give rise 
to several theorems: $ \vec{F} = m \vec{a}$, i.e. Newton's law, expresses a {\em theory} of mechanics rather than a {\em theorem} that nobody can prove it.  
However, it is true that, assuming the validity of Newton's theory, there are of course several theorems which follow from it, e.g. the conservation of energy for conservative forces. 
The success of Newton's theory is certified, for instance, by showing the validity to  Kepler's laws while its limit turns out to be the existence of relativistic phenomena. 
Another striking example of this successful strategy in Physics comes from atomic phenomena:  based on a huge and consolidate set of spectroscopic data, Heisenberg, Schr\"{o}dinger and 
many others set up and developed the theory of Quantum Mechanics. While it does not make any sense to ask whether one can {\em prove} the Schr\"{o}dinger equation, 
it is on the contrary extremely crucial to check that all its consequences are not only theoretically self-consistent but also never contradicted by any 
experiment.

These considerations are particularly useful in order to put in the right perspective the approach we are going to present for facing an interesting scientific problem which does not 
come however from the realm of physical phenomena but, on the contrary, directly from the realm of mathematics. It concerns  an  infinite class of  functions of the complex variable $s$, the so-called Dirichlet $L$-functions,  to which  belongs Riemann's zeta function. Postponing to later chapters the discussion of all  relevant details, a guideline for the present paper can be summarized as follows: 
\begin{enumerate}
\item  For the problem at hand, the consolidate set of empirical data of the problem consists of the explicit computation of a huge number of the so-called non-trivial zeros of these functions in the complex plane 
of the variable $s$. As a matter of fact, all these computed zeros are found to be {\em always} on the line ${\mathbb Re}\, s = 1/2$. So far there are no theorems which make such a property of 
these functions obvious and this is precisely why the problem is interesting and challenging for a curious mind. The possibility that all the non-trivial zeros of the Dirichlet $L$-functions 
are indeed on the line ${\mathbb Re}\, s = 1/2$ is known in the literature as the {\em Generalized Riemann Hypothesis} (GRH). 
\item As shown in the following, we are able to argue that there exists a very simple and elegant theory which explains at once why {\em all} non-trivial zeros of {\em all} the Dirichlet $L$-functions are always on the line ${\mathbb Re}\, s = 1/2$: 
according to this theory, any Dirichlet $L$-function is associated to a particular realization of a Brownian motion which -- it can be shown -- rules the location of its zeros. The ubiquitous value $1/2$ for the real part 
of all the zeros of these functions finds then a natural explanation in terms of the universal critical exponent $1/2$ associated to the growth of any Brownian motion and the validity of a central limit theorem behind any 
Dirichlet $L$-function. 
\item  Since what is proposed in this paper is an overall {\em theory} for the zeros of the Dirichlet functions rather than a {\em theorem}, in the last part of the paper we have performed a large number of testable attempts to falsify it. These consist of 
a huge battery of ``numerical experiments" designed to bring under a severe scrutiny the many different aspects of the supposed Brownian motion behind the Dirichlet $L$-functions. 
As discussed in detail later, the successful outputs (with an extremely high level of confidence, as certified by standard statistical tests) of {\em all} the numerical analyses 
performed can be seen as impressive experimental confirmation 
 of the theoretical framework proposed for the GRH. In this respect, we invite the reader to simply look at the various figures
presented in the last part of the paper, showing spectacular agreement between the Brownian motion theory and the numerics. 
\end{enumerate}

As will  become more and more clear in the following pages, in order to set up and support the Brownian motion picture behind the Dirichlet $L$-functions we had to embark in a long and detailed tour of 
many fascinating fields in physics and mathematics coming, in particular, from the common area of probability, stochastic phenomena and alike. The ideas which guide our work will be often illustrated with elementary arguments which have the clear advantage to present, in an essential way, the results we obtained: it is worth stressing that the style adopted in the presentation is the one of a theoretical physicist rather than a mathematician. In this respect, we do not claim to have any rigorous proof of the results presented here but we like to think that our work has been guided by the famous and beautiful Feynman's comment: {\em A great deal more is known than has been proved}. Moreover, if this work will have the effect to stimulate further rigorous studies by genuine mathematicians on the subject, it has already reached the scope to draw the attention to a possible way to tackle a long standing problem such as  Riemann's hypothesis.   

\begin{figure}[t]
\centering
\includegraphics[width=0.50\textwidth]{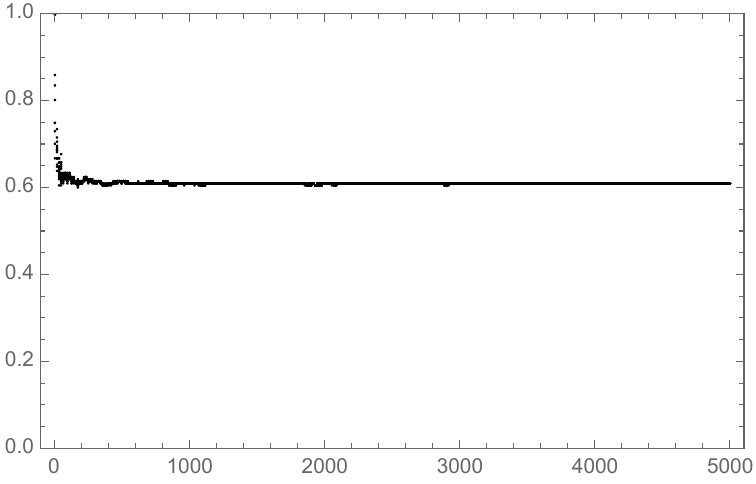}
\caption{Fraction ${\mathcal P}$ of square-free numbers in a sample of $n$ consecutive integers which very quickly converges to the predicted value ${\mathcal P} \sim 6/\pi^2$.
} 
\label{fffff}
\end{figure}

To make crystal clear the style adopted in this paper, let's present a paradigmatic and important example of a result which will be extremely useful later. Such an example concerns
 the square-free numbers. {An integer is said to be square free if no prime factor divides it more than once and these square free numbers will play an important role in the sequel of the paper.}
 Let us pose  the following question: {\em What is the probability ${\mathcal P}$ that an integer is square-free?} This is  equivalent to asking: {\em What is the density of square-free numbers among the integers?} To answer this question, our way of arguing will go as follows\footnote{This argument is presented in the book by M. Schr\"{o}der \cite{Schroeder} listed in our references, a text which is a true gem in the relations between Number Theory and Physics.}:
assuming the statistical independence of primes, the probability that a number $n$ is not divisible twice by the same prime $p$ is approximatively $(1 - 1/p^2)$ and, multiplying on all the primes, we get 
\beq
{\mathcal P} \sim \prod_{p} \left(1 - \frac{1}{p^2}\right) = \frac{6}{\pi^2} \,=\, 0.607927...
\label{probabilitysquarefree}
\eeq
Obviously, the first part of the line above plays the role of a {\em theory}, while the second part is a mathematical identity, i.e. a {\em theorem}. But, can we check the validity of this theory? Yes,  by simply counting the fraction of square free numbers within the integers! Doing this numerical experiment, with the corresponding plot shown in Figure \ref{fffff}, one can see that, quite rapidly, the measured density of square-free numbers within the integers indeed tends to the predicted theoretical value $6/\pi^2$. 
\vspace{1mm}
This example follows paradigmatically the A-B-C scheme mentioned above, since it elaborates a theory based on some assumption, it leads to a theorem (the infinite product on primes given above {\em is} indeed equal to $6/\pi^2$) and it is {\em not} falsified by a direct experiment. Of course the argument misses any mathematical rigor since, for instance, it never specifies how the probability is defined, if there is a proof of the independence of the divisibility of an integer by two prime numbers $p$ and $q$, what is the actual meaning of the symbol $\sim$, and finally the value to assign to a ``numerical check" of a mathematical property. We appeal once more to Feynman to go on with our presentation. 

A final disclaimer: this paper is rather long,  however, by no means should it be regarded as a review of the (Generalised) Riemann Hypothesis and the associated literature: the reader interested in  the history of the GRH may find satisfaction in reading some of the large literature specialized to this subject such as, for instance 
\cite{Riemann,Edwards,Titcmarsh,Davenport,Bombieri,Sarnak,Conrey,Polya,Borwein,Broughan,reviewRiemann,Wolf}.

Let's now go to the next section which provides a general overview of the topic and a presentation of the various parts of the paper. 
 
\section{A bird's eye view of the problem}

 An interesting aspect of the studies presented in this paper is the connection with the ``experimental" sides of mathematics, in particular those concerning with time series and the pseudorandom nature of their coefficients. Let's explain  briefly what  this is about.  Among other things, this paper deals indeed with an infinite binary sequence $\{{\mathcal S}_n\}$ made up of $\{0,1\}$'s  which looks like
\beq
\ldots \ldots 0, 0, 1, 0, 1, 0, 0, 1, 1, 1, 0,  0,  0,  0,  1, 1,  0,  1,  1,  0,  0,  0,  1,  1,  0,  1,  0,  1,  0,  0, 0,  1, 0,  1,  0, 0, \ldots\ldots 
\label{sequencebinary}
\eeq
It is  assumed that we can have access to arbitrarily large parts of this sequence, although always finite. One of our purposes will be to establish whether $\{\mathcal {S}_n\}$ is a random sequence and whether 
its properties are captured by the normal law distribution. As we will see, a relevant prior information on our sequence $\{\mathcal {S}_n\}$, which will be quite important in our analysis, is that such a sequence has on average an equal number of $0$'s and $1$'s. There is  also additional information on the sequence  $\{\mathcal {S}_n\}$ which will be presented later and which help us in better pinpointing the pattern of the appearance of the $0$'s and $1$'s. We will discuss below the origin of this problem, as well as its deep connection with an old-standing question of Number Theory: the Riemann Hypothesis  
\cite{Riemann,Edwards,Titcmarsh,Davenport,Bombieri,Sarnak,Conrey,Polya,Borwein,Broughan,Apostol,Iwaniec,Bombieri2,Steuding,Sarnak2,BK0,BK00,BK000,BK1,BK2,BK3,BK4,Bost,
Connes,KeatingSnaith,Sierra1,Sierra,Sierra2,Srednicki,Bender,reviewRiemann,Wolf,GriffinZagier,RodgersTao}.  It is however important to discuss some interesting methodological issues in these numerical studies.

\subsection{Time series}

If one does not have any prior information on an infinite sequence of numbers,  the problem of  establishing its nature can only be addressed statistically, i.e. in terms of a statistical study which  generically  goes under the name of 
{\em Time Series Analysis}, a well established and very powerful branch of Probability Theory  
 \cite{Knuth,Nist,Marsaglia,Feller,Jaynes,timeseries1,timeseries2,timeseries3}. It is worth stressing the power of this approach: using the Time Series Analysis to study sequences of numbers such as $\{\mathcal{ S}_n\}$, one can reach very robust and significant conclusions on whether the sequence under scrutiny is random or not. Indeed, there is no limit to the number of tests that can be developed and applied to an infinite sequence of numbers like the ones in (\ref{sequencebinary}) to check its level of randomness. Some of these tests are rather simple while others, as we will see, can  instead be very elaborate and sophisticated. 
 
 In pursuing this  kind of analysis, of course there is always a lurking {\em regressus ad infinitum} objection to face: namely,  if a sequence appears to be random under the tests $T_1, T_2, \ldots T_n$, how can we be sure that it will also appear to be  random under a new test $T_{n+1}$? 
 We cannot, of course. Assume, however, that  no matter how we increase the number and the quality of the statistical tests,  {\em any} arbitrarily large parts of $\{{\mathcal S}_n\}$  {\em always} passes them  successfully.  Well, at this point, our level of confidence of its randomness starts then to increase considerably, and we can quantify it in terms of a probability that becomes closer and closer to 1 as we increase the number and the refinement of our tests. In Donald Knuth's words, \textquotedblleft the sequence is presumed innocent until proven guilty" \cite{Knuth} and it is more and more innocent  by increasing the level of our screenings. One may even push such an argument forward and argue that this assumption is at the root of our general way of getting knowledge:  to form a judgment about the likely truth or falsity of any proposition A, the correct procedure is to calculate the probability that A is true, $P(A | E_1,E_2,\ldots)$ conditional on all the evidences $E_i$ at hand \cite{Jaynes}. In Physics these statements are, of course, very familiar: for instance, the discovery of the Higgs boson \cite{Higgs}, where the ATLAS and CMS experiments at CERN's Large Hadron Collider announced they had each observed a new particle in the mass region around 125 GeV, is after all a nice example of knowledge acquired by ``statistics'', given that this new particle owes its existence only to the satisfaction of stringent statistical analysis of the events collected at LHC. In Mathematics, on the other hand, knowledge reached by probabilist arguments seems to be almost flawed by an inherent uncertainty. But is it always the case? To make some progress on this question, let's discuss now in more detail the origin of our sequence.

\subsection{Probability in Number Theory} Our sequence $\{\mathcal {S}_n\}$ has an interesting origin and is deeply connected to one of the most famous problems in Number Theory, i.e. the Riemann conjecture about the location of the zeros of the Riemann zeta-function $\zeta(s)$ in the complex plane of the variable $s$ 
\cite{Riemann,Edwards,Titcmarsh,Davenport,Bombieri,Sarnak,Conrey,Polya,Borwein,Broughan,Apostol,Iwaniec,Bombieri2,Steuding,Sarnak2,BK0,BK00,BK000,BK1,BK2,BK3,BK4,Bost,Connes,
KeatingSnaith,Sierra1,Sierra,Sierra2,Srednicki,Bender,reviewRiemann,Wolf,GriffinZagier,RodgersTao}. 
As we will see, the problem is also related to a famous topic of Statistical Physics, i.e. the random walk problem \cite{Yuval,Mazo,Rudnick,Levy,diffusion} and, in particular, 
 how it is possible to establish the random nature of a Brownian motion when one has access to {\em only one} single trajectory rather then a collection of trajectories \cite{Perrin,Nordlund,Kappler,brow1,brow2,brow3,brow4}.   Although it is not surprising that probability concepts are applied to the study of Brownian motion,  one may wonder what probability has to do instead with Number Theory, a world dominated by the rigid rules of integer numbers and the like. However, some of the most remarkable progress of the latest decades is that many properties of Number Theory can be strikingly captured by ideas, methods and results which come indeed from the realm of probability, that branch of mathematics which describes aleatoric events 
\cite{Schroeder,Kac,erdosbook,ErdosKac,Kubilius,Tao,Cramer,Billingsley,Dyson,Montgomery,Odlyzko,Rudnick-Sarnack,Grosswald,EPFchi,Franca1,ML,LM, 
Denjoy, Churchhouse}.  Among the  pioneers of this probabilistic approach to Number Theory  are Mark Kac \cite{Kac} and Paul Erd\H{o}s \cite{erdosbook}.  
A nice example of what can be achieved by 
adopting such a point of view is provided by their famous Erd\H{o}s-Kac theorem \cite{ErdosKac} which states that, if we denote by $\omega(m)$ the number of distinct prime factors of the integer $m$, then in the large $m$ limit, the probability distribution of the variable 
\beq
x\,=\, \frac{\omega(m) - \log\log m}{\sqrt{\log\log m}} 
\eeq
is the standard normal distribution ${\mathcal N}_{0,1} (x)$, with mean $0$ and variance $1$, where  ${\mathcal N}_{\mu,\sigma^2} (x)$ stands for the 
normal distribution of a random variable $x$ with mean $\mu$ and variance $\sigma^2$ given by 
\beq
{\mathcal N}_{\mu,\sigma^2} (x) \,=
\,
\frac{1}{\sqrt{ 2 \pi \sigma^2}} \, \exp\left[{- \frac{(x-\mu)^2}{2 \sigma^2}}\right] \,\,\,. 
\eeq
It is important to emphasize that  regardless of the fact that $\omega(m)$ is a completely deterministic  arithmetic function, the distribution of its values follow the probabilistic normal distribution. 
Both the normal distribution and the Erd\H{o}s-Kac theorem will accompany us in the course of our discussion and we will have the opportunity to 
comment more on their importance. 
{Here,  the ``randomness" is related to the fact that the event of an integer being divisible by a prime $p$ is independent of being divisible by another prime $q$.  
This is in essence the origin of randomness in this paper. }

\def\Re{{\mathbb Re}}

\subsection{Riemann-zeta function}
Referring  all the details to the coming Sections, let's see the context in which  our sequence (\ref{sequencebinary}) emerges.  Consider the Riemann zeta function 
$\zeta(s)$, given by its series and infinite product representation on the prime numbers $p_k$ ($k=1,2,3,\ldots$), hereafter ordered wrt their increasing values 
\beq
\zeta(s) \,=\,\sum_{m=1}^{\infty} \frac{1}{m^s} \,=\,\prod_{k=1}^\infty \frac{1}{1 - \frac{1}{p_k^s}} \,\,\,. 
\label{riemannzeta}
\eeq
The above definition converges for $\Re \, s > 1$,  and the function can be analytically continued to the whole complex $s$ plane.  
This function has a pole at $s=1$, an infinite number of (trivial) zeros at $s = - 2 \,n $ (where $n=1,2,\ldots)$ and an infinite number of other zeros in the critical strip $0 \leq {\mathbb Re}\, s \leq 1$ which are supposed to satisfy the Riemann Hypothesis  \cite{Riemann,Edwards,Titcmarsh,Davenport,Bombieri,Sarnak,Conrey,Polya,Borwein,Broughan,Apostol,Iwaniec,Bombieri2,Steuding,Sarnak2,BK0,BK00,BK000,BK1,BK2,BK3,BK4,Bost,Connes,
KeatingSnaith,Sierra1,
Sierra,Sierra2,Srednicki,Bender,reviewRiemann,Wolf,GriffinZagier,RodgersTao}: 

\vspace{3.5mm}
\fbox{
\parbox{15cm}{
{\bf Riemann Hypothesis (RH)}:  {\em in the critical strip $0 \leq {\mathbb Re}\, s \leq 1$, the zeros of $\zeta(s)$} \\
{\em are 
{\underline{\em all}} and {\underline{\em only}} 
on the infinite line $ {\mathbb Re} \,s = 1/2$.}
}}
\vspace{3.5mm}

\noindent
Equivalently, if the RH is true, the multiplicative inverse function of $\zeta(s)$, defined as well as by its series and infinite product representation 
\beq
\inv{\zeta(s)} \,=\,\sum_{m=1}^{\infty} \frac{\mu(m)}{m^s} \,=\,\prod_{k=1}^\infty \left(1 - \frac{1}{p_k^s}\right)\,\,\, ,
\label{mufunction}
\eeq
must necessarily have {\em all} its poles along the line $ {\mathbb Re} \,s = 1/2$. The arithmetic function $\mu(m)$, known as the M\"obius function, has values 
\beq
\label{mudef}
\mu(m)  \equiv
 \begin{cases}
 -1 \qquad \mbox{if $m$ is square free and has an odd number of prime factors} \\
~0 \qquad  ~~\mbox{if $m$ has a squared prime factor}\\
 ~1 \qquad ~~
\mbox{if  $m$  is square free and has  an even number of prime factors.} 
\end{cases}
\eeq
Using the Stieltjes measure, we can express $1/\zeta(s)$ in terms of its inverse Mellin transform 
\beq
\inv{\zeta(s)} \,=\, s\, \int_1^{\infty} \frac{M(x)}{x^{s +1}} \, dx \,\,\,,
\label{Mellin}
\eeq
where $M(x)$ is the so-called Mertens function, given by 
\beq
\label{Mertenfunction}
M(x) \,= \, \sum_{1 \leq m \leq x} \mu(m) \,\,\,.
\eeq
It is then simple to see that, if asymptotically $M(x)$ goes as $M(x) \sim x^{1/2 + \epsilon}$, for any arbitrarily small positive $\epsilon$,   then the RH is indeed true: in this case, in fact, the integral (\ref{Mellin}) diverges at $ {\mathbb Re} \,s = 1/2$, making clear the presence of a singularity on this axis
\footnote{{Throughout this paper, the symbol $\sim$ signifies  the ``big O" notaton,
$M(x) = O(x^{1/2 + \epsilon})$ for any arbitrarily small and positive $\epsilon$.  We will not always indicate the $\epsilon$ and may write simply
$M(x) \sim x^{1/2}$.}}.   In the following we will call this approach to the RH the {\em statistical/horizontal approach}, to be compared later with the {\em quantum/vertical approaches} to the RH.

\subsection{Law of iterated log's} Note that all terms with  $\mu (m)=0$ do not contribute to the sum in $M(x)$ and that the non-zero values of $\mu(m)$ only come from {\em square-free numbers}. As we have 
seen in eq.(\ref{probabilitysquarefree}), these numbers are a {\em finite fraction} of all integers.  Their first representatives are 
\beq
\label{fn}
\{ \sqf_1, \sqf_2, \sqf_3, \sqf_4,  \sqf_5, \ldots \} = \{1,2, 3, 5, 6, 7, 10,11, 13,\ldots \} 
\eeq
and the $n$-th square-free number $\sqf_n$ scales as 
\beq
\sqf_n \,\sim\, \frac{\pi^2}{6} \, n\,\,\,.
\label{scalefree}
\eeq 
Let's then introduce the map $f : {\mathbb N} \rightarrow \sqf$ such that $f(n) = \sqf_n$ and, restricting the attention only to these square-free numbers, 
let's also define 
\beq
\label{muhat}
\muhat (n) \,=\, (\mu \circ f)(n)\,=\,  \mu (\sqf_n )\,\,\,.
\eeq
 The relation of these coefficients with our initial ${\mathcal S}_n$ is simply
 \beq
 \boxed{
 \muhat(n) \,=\, 2 \, {\mathcal S}_n -1 \,\,\,. 
 }
 \label{relationmuands}
 \eeq
 In the following we will also denote these coefficients as $\hat\mu_n$. 
For our future considerations, in particular those concerning with the statistical analysis, sometimes we find  it useful to deal with our original sequence $\{\mathcal {S}_n\}$ while in other cases we find it more useful to focus our attention on the sequence $\{\hat\mu_n\}$, keeping in mind that the relationship between the two quantities is expressed by eq.\,(\ref{relationmuands}). It is important to stress that, focusing the attention only on the values taken by the M\"{o}bius coefficients on square-free numbers is the way to avoid all ``periodicities"\footnote{We have, for instance, that $\mu(n) =0$ whenever $n$ is multiple of $4$, $9$, $25$, $36$, etc., and, 
in general, a multiple of any $p_n^2$ or $\sqf_n^2$. These considerations imply that $\mu(n)$ has many {multiplicative}  periodicities while $\hat\mu(n)$ is expected to have none.} 
 present in the original arithmetic function $\mu(n)$ and, therefore, to better evaluate the random properties of this function which emerge by considering its restriction $\hat\mu(n)$ to the square-free numbers.{ In other words, we take seriously the observation on $\hat\mu(k)$ made originally by Denjoy \cite{Denjoy} (see also the book by Edwards \cite{Edwards} p.\,268) and we develop the main heuristic of this paper based precisely on the pseudo-random properties of the restricted M\"{o}bius coefficients $\hat\mu(n)$ to the square-free numbers.}

In the following the final object of our study is the {\em restricted Mertens function} $\hat M(n) = M \circ f$,  i.e. the Mertens function restricted to the square-free numbers only
 \beq
 \hat M(n) \equiv (M \circ f)(n) \,=\,\sum_{k=1}^n \muhat(k) \,\,\,.
 \label{hatmertens}
 \eeq
 Notice that, in light of the scaling law (\ref{scalefree}), we have the scaling relation 
 \beq
 \hat M(n) \,\sim \,  M\left(x=\frac{\pi^2}{6} n\right) \,\,\,.
 \eeq 
If we could show that the sequence $\{\hat\mu_n\}$ follows the stochastic laws of a random motion, then for the 
asymptotic growth of $\hat M(n)$ we would have  
\beq 
\hat M(n) \sim n^{1/2 + \epsilon}\,\,\,, 
\label{aymptotic}
\eeq
for an arbitrarily small positive $\epsilon$, which,  as explained above,  leads to an understanding of the origin and the validity of the RH.

\begin{figure}[t]
\centering
\includegraphics[width=0.45\textwidth]{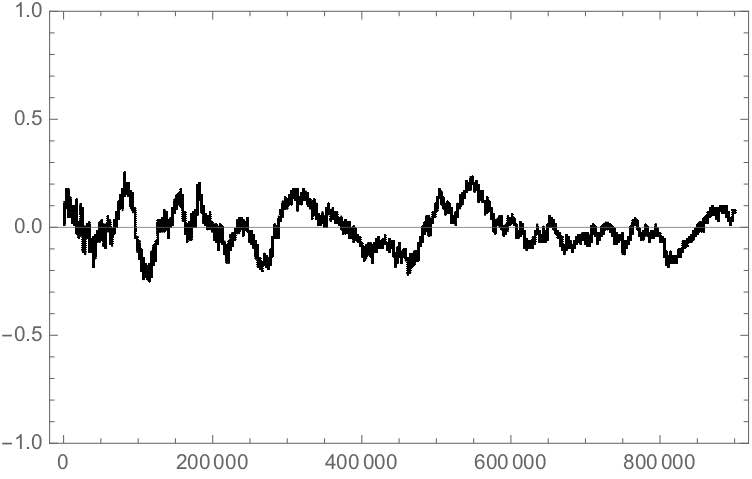}
\,\,\,
\includegraphics[width=0.45\textwidth]{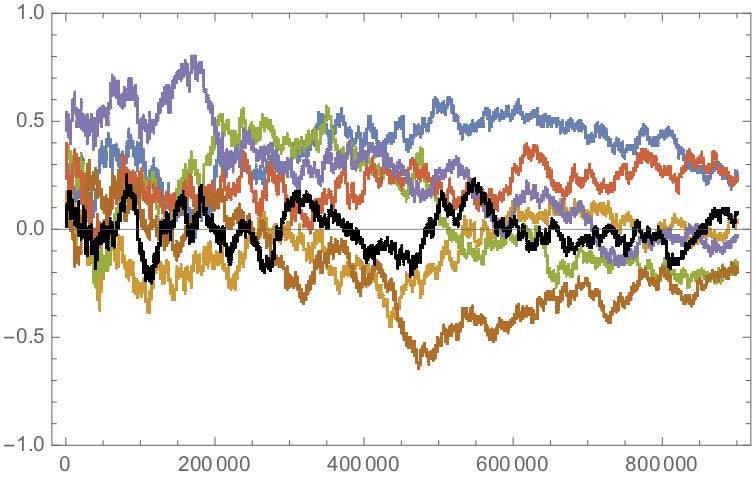}
\caption{(a) Left-hand side: the plot of $\hat M(n)/\sqrt{2 n \log\log n}$ vs $n$, for $n \leq 10^6$; (b) Right-hand side: the plot in colors of several random instances of 
$R(n)/\sqrt{2 n \log\log n}$  while, in black, the curve of $\hat M(n)/\sqrt{2 n \log\log n}$. 
} 
\label{Mertenseq}
\end{figure}

To show that $\hat M(n)$ has indeed such a behaviour, in the Part B of this paper we will present a {\em global approach}\footnote{As it will become clear in 
Sections \ref{globalapproach1} and \ref{globalapproach2}, the terminology {\em global approach} refers to the possibility of studying the function $\hat M(n)$ 
by grouping globally and at once large sets of its coefficients $\hat\mu(k)$. In contrast, we will refer to the {\em local property} of the function $\hat M(n)$ when 
we focus our attention just on the individual coefficients $\hat \mu(k)$.} to this function which involves the prime number distribution along the series of square-free numbers, the average number of prime divisors, the Erd\H{o}s-Kac theorem for square-free numbers, etc. 
Moreover, we will rely on the law of iterated logarithm which describes the magnitude of the fluctuations of a random walk \cite{Feller} to argue that, asymptotically, for the fluctuations of the 
function $\hat M(n)$ we have\footnote{For the original Mertens function $M(x)$ such a behavior (with different coefficient) was also conjectured in \cite{Churchhouse} on the basis of some computer studies: the analysis pursued here not only significantly enlarges the set of theoretical tools which leads to understand the validity of the formula (\ref{sqrtloglog}) but also provides for the first time a very robust support to the random behavior of any sufficiently large subsequence of $\{\mathcal{S}_n\}$ on the basis of many and very sophisticated statistical tests.}
\beq
\limsupn  \, \frac{| \hat M(n)| }{\sqrt{2 n \log\log n}} \,= \, 1 
\hspace{3mm} 
, 
\hspace{3mm}
{\rm a. s. } 
\label{sqrtloglog}
\eeq
where ``a.s.'' stands for ``almost surely'', in a probabilistic sense,  i.e. with probability equal to $1$ \cite{Feller}.
 As is well known, this law of iterated logarithm applies to all Brownian motion random sequences 
\beq
R(n) \,=\, \sum_{k=1}^n a_k 
\hspace{3mm}
, 
\hspace{3mm}
a_m \,=\,\{\pm 1\}\,\,\,, 
\eeq
where $a_k$ are random uncorrelated variables. The quantity $\hat M(n)/\sqrt{2 n \log\log n}$ is plotted in Figure \ref{Mertenseq}.a for $n \leq 10^6$ and 
compared in Figure \ref{Mertenseq}.b with the behavior of an analogous quantity $R(n)/\sqrt{2 n \log\log n}$ computed for several random uncorrelated sequences.  
These figures are  purely illustrative  of the brownian motion behavior of the restricted Mertens sequence which, as we will explicitly check in great detail later, persists to all scales, i.e. no matter how  much we enlarge the interval in which we study this sequence. 

It is worth stressing that there is nothing which prevents a deterministic sequence from being a realization of a random process, and the restricted M\"{o}bius coefficients together with the associate Mertens function seem indeed to be a very good example of this fact.  
{For instance,   imagine the random process of  a person  flipping a coin a large number of times,  and recording the sequence of the outcome of the flips 
 as a list in a sacred book kept in 
a secure chamber in a church,  or the National Bureau of Standards.       
The existence of the book now makes the sequence completely deterministic,  since anyone can  access the book,  and everyone agrees on its precise content.    It's a rather boring book and, unlike the gambler,  an  individual feels no nervous anticipation of 
what will be the next flip,  since it is already all known to anyone.}

\subsection{Generalized Riemann Hypothesis} To put the problem in a proper perspective, it is very useful to remind (see Part A for further details) that the Riemann $\zeta$-function is just a particular case of a more general class of analytic functions known as Dirichlet L-functions $L(s,\chi)$, where $\chi$ is an arithmetic character which we will discuss in more detail in Section \ref{Dirichelet} of Part A. There are  infinitely many of these functions in relation to the infinitely many characters $\chi$ and, for ${\mathbb Re} \,s > 1$, they admit both an infinite series and product representation 
\beq
L(s,\chi)\,,=\,\sum_{m=1}^\infty \frac{\chi(m)}{m^s} 
\,=\, 
 \prod_{k=1}^\infty  \( 1 -  \frac{\chi (p_k)}{p_k^s} \)^{-1} \,\,\,,     ~~~~ {\mathbb Re} \,s > 1\,. 
\label{Lfunctionsfirst}
\eeq 
For all these functions,  the Generalized Riemann-Hypothesis has been  conjectured to hold:

\vspace{3.5mm}
\fbox{
\parbox{15cm}{
{\bf Generalized Riemann Hypothesis (GRH)}:  {\em in the critical strip $0 \leq {\mathbb Re}\, s \leq 1$, the zeros of  all the infinitely many L-functions $L(s,\chi)$} 
{\em are 
{\underline{\em all}} and {\underline{\em only}} 
on the infinite line $ {\mathbb Re} \,s = 1/2$.}
}}

\vspace{3.5mm}
\noindent
Notice that 
\beq
\frac{1}{L(s,\chi)} \, =\, \prod_{k=1} \left(1 -  \frac{\chi (p_k)}{p_k^s} \right) \,=\, \sum_{m=1}^\infty \frac{\mu(m) \chi(m)}{m^s} \,\,\,, 
\eeq
and this function can be expressed in terms of its inverse Mellin transform 
\beq
\frac{1}{L(s,\chi)} \,=\, s \,\int_1^{\infty} \frac{M_{\chi}(x)}{x^{s+1}} dx \,\,\,,
\label{mellinl-1}
\eeq
where the analog of the Mertens function for the Riemann  zeta function is played in this case by
\beq
M_{\chi}(x) \,=\, \sum_{1 \leq n \leq x} \mu(n) \chi(n) \,\,\,,
\label{generalisedMertens}
\eeq 
which we call the {\em Generalized Mertens function}.  As in the case of the Riemann zeta function, if it can be shown that $M_\chi(x)$ goes asymptotically as $M_\chi(x) \sim x^{1/2 + \epsilon}$, for any arbitrarily small $\epsilon$, then the GRH is indeed true. We will discuss later the relation which links together $M_\chi(x)$ to $M(x)$ (see Section \ref{GRHRH}) and, in view of this relation, how  it is possible to unify the two hypothesis, in particular, how it is possible to argue that the validity of the GRH can be considered as a consequence of the RH. 

\vspace{1mm}

Let us stress that one of the key features of the statistical approach to the RH pursued in this paper is its level of naturalness and universality which helps in clarifying at once why {\underline {\em all}} the infinite number of Dirichlet $L$-functions (including the Riemann $\zeta$-function) have {\underline {\em all}} their non-trivial zeros on the line\footnote{
Denoting by $\rho = \sigma + i t$ the position of a generic zero of a Dirichlet function, it is worth to point out that the statistical approach presented here is extremely powerful 
in identifying the abscissa of all these zeros, leading to the conclusion that for all of them $\sigma = 1/2$, but it is literally unable to address the positions of their imaginary part $t$. On the other hand, one should keep in mind that the imaginary parts of these zeros are not at all universal since they change by changing the $L$-functions.}   
 $ {\mathbb Re} \,s = 1/2$. As it will become clear later, the reason is indeed quite interesting: the location of all the zeros of all these functions are ruled by sequences of numbers which, albeit strictly deterministic, behave though as brownian random walks! In other words, the universal location on the line $ {\mathbb Re} \,s = 1/2$ of all the zeros of these functions can be traced back to the universal statistical properties of the brownian random walks, such as the existence of the central limit theorem or the validity of the law of iterated log's. 

\subsection{Organization of the paper} The paper is divided into four Parts, each of them addresses different aspects of the problem and therefore has a different style and length. 

\vspace{1mm}
{\bf Part A} presents the general framework of the Dirichlet $L$-functions and explains why we can consider the Riemann zeta-function as a particular example thereof. We discuss in particular the reason why the argument which we previously used to show the validity of the GRH for the $L$-functions of non-principal characters \cite{ML,LM} cannot be used to show the RH for the Riemann zeta-function. However, we argue that the validity of the RH would imply the validity of the GRH.  In Part A we also discuss some important results about the properties of the Dirichlet $L$-functions which justify the use of a probabilistic approach for the GRH.  We recall, in particular, both the theorem by Grosswald and Schnitzer \cite{Grosswald} concerning a set of random functions which share the same zeros of the Dirichlet (Riemann) functions in the critical strip, and our previous results on the location of the zeros of the Dirichlet L-functions of non-principal characters. So, besides the original section on the relation between the RH and the GRH, Part A contains well known but also less known properties of the Dirichlet $L$ functions and is meant to be just a reference part for the rest of the paper, where the reader can easily find the most relevant definitions, the discussion of some key features of these functions, as well as some key points of our approach.

\vspace{1mm}
{\bf Part B} is the core original theoretical part of this paper, where we present our new approach to the restricted Mertens function, which we call the {\em global approach}.  It concerns  a series of new results relative to the sequence of the square-free numbers which lead us to argue positively about the validity of the RH on the basis of a probabilistic argument. In this part of the paper, we address important issues of Number Theory such as the distribution of prime numbers along the sequence of the  square free-numbers, the important role of the {\em primorial} (the equivalent of the factorial for the prime numbers) in controlling the growth of the Mertens functions, the average  number of prime divisors, the Poissonian distribution satisfied by the prime divisors and the Erd\H{o}s-Kac theorem for square-free numbers. All these results are instrumental for us to {\it argue} that the moments of the restricted Mertens function asymptotically behave as those of a random walk, i.e. 

\beq
\langle (\hat M(n))^{2 k+1}\rangle \,=\, 0 
\,\,\,\,\,\,\,\,\,\,\,\,\,\,\,\,
;
\,\,\,\,\,\,\,\,\,\,\,\,\,\,\,\,
\langle (\hat M(n) )^{2 k} \rangle\,=\, n^k \, (k-1)!! \,\,\,.
\eeq

\vspace{1mm}
\noindent
The careful reader may have noticed that the expressions above refer to some average quantities, i.e.  $\langle \cdots \rangle$. We will explain 
the meaning of taking such an average for a deterministic sequence as the one given by the M\"{o}bius coefficients, and this leads us, in particular, to define the proper statistical 
ensemble in which this average makes sense. As we will see, the problem is remarkably close to the so-called Single Trajectory Random Walk problem
\cite{Nordlund, Kappler, brow1,brow2, brow3, brow4}, which consists in establishing the random nature of a brownian motion when one has access to {\em only one} single trajectory rather then a collection of trajectories. In a nutshell, we will see that all the considerations related to the average of the restricted Mertens function have a precise analog in statistical mechanics, when one substitutes {\em phase space average} with {\em time average}, under the assumption of the ergodicity and time translation of the system under scrutiny. 

\vspace{1mm}
{\bf Part C} presents a large battery of statistical tests which we have performed on very long sub-sequences of the restricted M\"{o}bius coefficients.  Altogether 18 statistical tests were performed,   each of them repeated thousands of times in very large intervals along the infinite sequence of the restricted M\"{o}bius coefficients $\hat\mu(n)$. All these studies aim to probe the {\em local properties} of $\hat\mu(n)$ and the associated restricted Mertens function $\hat M(n)$ and they were triggered by our natural curiosity to see how good such finite sub-sequences behave randomly, as predicted on the basis of the global theoretical considerations on the restricted Mertens function that we presented in Part B. In this part we extensively use the block variables associated to the sub-sequences of $\{{\mathcal S}_n\}$'s for reaching  very robust probabilistic conclusions on a deterministic sequence as the one given by the restricted Martens sequences $\hat M(n)$. The statistical tests applied to our sequence $\{{\mathcal S}_n\}$ are of increasing order of complexity and sophistication. Let's anticipate and, at the same time emphasize, that {\underline {\em  all}} the sub-sequences analyzed passed successfully our statistical tests, with an overall level of confidence of $99\%$,  as quantified by the associated $\chi^2$ distribution. To the best of our knowledge, never before  has such a massive statistical analysis  been performed 
on the restricted M\"{o}bius coefficients and Mertens function, and the outputs of this analysis may be regarded as a very robust and striking  ``experimental" confirm of the RH. . 

\vspace{1mm}
{\bf Part D} contains our conclusions.  In this section we come back to the relationship between GRH and RH and we also make a comparison of our result versus other conjectured results concerning the growth of the Mertens function. We will comment, in particular, on some unpublished conjectured formulas due to Gonek which were obtained, though, assuming the validity of the RH. We will show that our prediction on the growth of the restricted Mertens function, obtained {\em without} assuming the validity of the RH\footnote{As a matter of fact, our logic has been the reverse! Namely,  the Riemann Hypothesis holds true on the basis of the behavior we have established for the restricted Mertens function.}, is perfectly compatible with the estimated growth of the restricted Mertens function obtained assuming the Riemann Hypothesis.

\begin{figure}[t]
\centering
\includegraphics[width=0.80\textwidth]{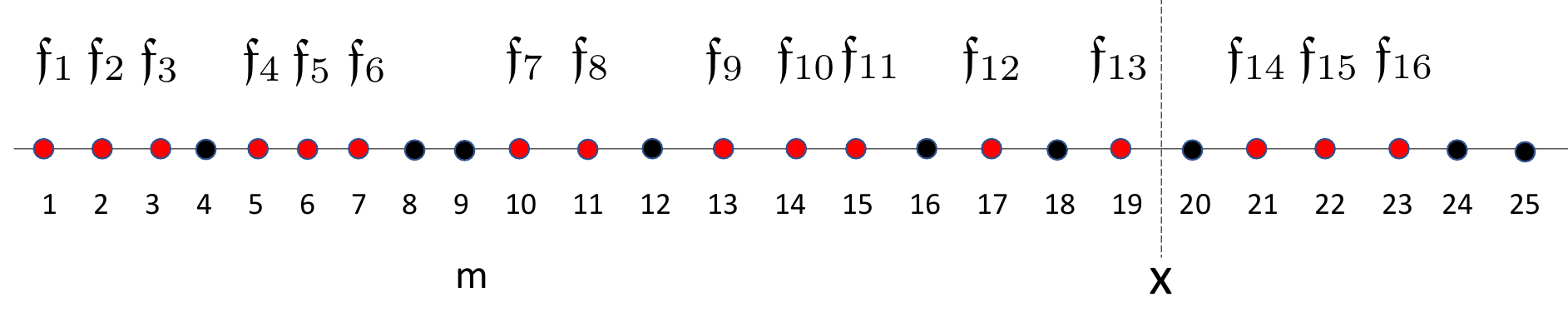}
\caption{Sequence of integers $m$ (black dots) and square-free numbers $\sqf_n$ (red dots). $x$ denotes a generic real number. Square-free numbers have 
a density $6/\pi^2$ among the integers.
} 
\label{scalingfiguresf}
\end{figure}

\vspace{3mm}
\noindent
{\bf Notation}. In the following we often use the index ${n}$ to refer to the square-free numbers $\sqf_n$ and quantities related to them.    
For an arithmetic function $F(n)$ defined over the integers,   the same function over the square free numbers $\sqf_n$ will be denoted as 
$\hat F (n) \equiv  F(\sqf_n)$.   This applies to  the restricted  Mertens  function $\hat M(n)$ given in eq.\,(\ref{hatmertens}), the restricted M\"obius function, prime number counting function
and number of divisors functions  $\hat \mu$,  $\hat \pi$ and $\hat \omega$ defined below. 
 On the other hand, we use ${m}$  and ${k}$ to label a generic integer and ${x}$ to denote a generic real number, whose integer part $[x]$ are related to the integers. It is of course simpler to use $x$ rather than $n$ for  denoting intervals but one has to keep in mind that $n$ is the true number of terms in a sum as the one in (\ref{hatmertens}). The scaling relation between $x$ and $n$ is given in eq.\,(\ref{scalefree}) and illustrated in Figure \ref{scalingfiguresf}. 
Moreover, as usual, in the following we use the notation 
\beq
s \,=\, \sigma + i t
\eeq
to denote the real and imaginary parts of the complex variable $s$.

\newpage

\addcontentsline{toc}{section}{Part A}
\begin{center}
{\Large {\bf PART A}}
\end{center}

 The aim of this part is to clarify a few important points of the subject that help in illuminating the nature of the RH. As we are going to see,  a crucial role in our future discussions is played by probabilistic arguments, characterised by their economy and simplicity.  Moreover, they provide a natural way to provide a concrete reason  for the validity of the Riemann and the Generalized Riemann Hypothesis.  Indeed, the probabilistic approach pursued in this paper has the advantage to disclose, in a very natural way, the reason why the non-trivial zeros of the Riemann zeta-function have to be along the line $ {\mathbb Re} \,s = 1/2$ in the complex plane and, at the same time, why this property is also shared by the infinitely many Dirichlet $L$-functions. In other words, if arguments based on probability provide on one side a very robust ground of naturalness for the Riemann and the Generalised Riemann Hypothesis, on the other side they also point to the very high level of universality behind these hypotheses, ultimately related to the properties of the random walk, i.e. to the ubiquitous appearance of the normal law distribution and the central limit theorem.  

In this part of the paper we show that the Riemann zeta functions belongs to an infinitely large class of functions known as Dirichlet $L$-functions: these may be regarded as generating functions constructed in terms of local data associated with an arithmetic object.  Dirichlet $L$-functions are particular examples of so-called {\em Dirichlet series}, which provide very useful tools in analytic number theory. One of the main properties of the Dirichlet series is their half-plane absolute convergence. In the case of the $L$-functions, another important property is that, besides their series representation, they also admit an infinite product representation over prime numbers. In the following, we will give a very short overview of the $L$-functions, in particular focusing on their origin from number theory and their analytic structure in the complex plane, referring to classical texts of the literature for an extended discussion of their properties \cite{Apostol,Iwaniec,Bombieri2,Steuding,Sarnak2}.

\vspace{3mm}
\noindent
{\bf Setting the stage}. The main purpose of this Part is to set the stage for the results presented later in this paper. In the next sections, in particular, we intend to clarify the 
following topics: 
\begin{itemize}
\item Why is the Riemann zeta function a particular example of the Dirichlet $L$-functions?
\item What is so special about the Riemann zeta function? In particular, why the approach \cite{ML,LM} successfully applied to show the validity 
of the Generalized Riemann Hypothesis for a generic $L$-function of non-principal characters cannot be simply extended to the Riemann zeta function? 
\item Would  it be possible to show that the validity of RH implies the validity of the GRH?
\item What is the importance of the Grosswald-Schnitzer theorem on the zeros of a set of random functions which are the analogue of the 
Dirichlet $L$-functions? 
\item What is the difference between the Statistical Approach and the Quantum Approach for establishing the validity of the Generalized Riemann Hypothesis?
\item What do we know so far about the properties of the Mertens function and where  does this knowledge come from? 
\end{itemize}

\vspace{3mm}
\noindent
{\bf Playing with dice}. As we will see shortly, the answer to the first two questions is very straightforward: any Dirichlet $L$-function is associated to a natural number $q \geq 1$, known as the modulus, and the Riemann zeta-function simply corresponds to the $q=1$ case. Moreover, we will show that, behind all these functions, there is a stochastic process which may be regarded as the outputs of a dice made of $q$ faces. Hence, for the Riemann zeta-function, we are formally dealing with a dice having only $q=1$ face, which is what makes it difficult to apply our previous approach  \cite{ML,LM} directly to $\zeta(s)$.  {A signature of this issue is the pole at $s=1$ for Dirichlet L-functions based on principal characters, which is absent for those based on non-principal characters.}
As a matter of fact, this difficulty is what  motivated us to develop the alternative approach to the Riemann zeta function discussed in this article. However, interestingly enough,  the story has had a compelling turn, in the sense that  the approach developed here to deal with the RH may be also used to establish the validity of the GRH. 
Let's see in more detail how all this comes about.

\section{Dirichlet $L$-functions and the Riemann zeta-function}\label{Dirichelet}

The Dirichlet L-functions of the complex variable $s = \sigma + i t$ admit a series and an infinite product representations given by   
\beq
L(s,\chi) \,=\,\sum_{m=1}^\infty \frac{\chi(m)}{m^s} 
\,=\, 
 \prod_{m=1}^\infty  \( 1 -  \frac{\chi (p_m)}{p_m^s} \)^{-1} \,\,\,,     ~~~~ {\mathbb Re} \,s > 1\,,
\label{Lfunctions}
\eeq 
where $\chi(m)$ is a Dirichlet character and $p_m$ is the $m$-th prime in ascending order. 
Comparing with eq.\,(\ref{riemannzeta}), it is easy to see that the Riemann zeta function corresponds to $\chi(m) =1$, for all natural numbers $m$. 
It is convenient to briefly discuss these Dirichlet characters for better appreciating the nature of the Dirichlet $L$-functions, in particular to 
discover that everything starts from a classical problem in number theory. 


\vspace{3mm}
\noindent
{\bf Primes in Arithmetical Progressions}. A problem which attracted the attention of Dirichlet in 1837  was to prove that there is an infinite number of primes in 
arithmetic progressions such as 
\beq
A_m \,=\, q \, m + h \,\,\,\,\,,\,\,\,\, m=0, 1, 2, \ldots 
\,\,\,\,\,\,\,\,\,\,\,\, 
q, h \in \mathbb{N} 
\label{sequence}
\eeq
where the number $q$ is known as the {\em modulus} while the number $h$ as the {\em residue}. Dirichlet proved that if $q$ and $h$ have no common divisors, namely they are {\em coprime}, a condition expressed as $(q,h) = 1$, this is a sufficient and necessary condition for finding indeed infinitely many primes in the sequence (\ref{sequence}). 
Notice that, taking $q=1$, the statement is equivalent to say that there are infinitely many primes in the sequence of natural numbers. 
 
\subsection{Characters} Given an integer $q$, we consider all integers $m$ coprime with $q$, i.e. $(m,q) =1$. The set of these integers $m$, called the prime residue classes modulo $q$, 
under the multiplication mod $q$ forms an abelian group, denoted as   
\beq
(\mathbb{Z}/q\mathbb{Z})^* := \{m\, {\rm mod} \,q \,:\, (m, q) = 1\} \,\,\,. 
\label{groupab}
\eeq
The dimension of this group is given by the Euler totient arithmetic function $\varphi(q)$, defined as the number of positive integers less than $q$ that are coprime to $q$: its value 
is given by 
\beq
\varphi(q) \,=\, q \, \prod_{p | q} \left(1 - \frac{1}{p}\right) \,\,\,,
\label{Eulertotient}
\eeq
where the product is over the distinct prime numbers dividing $q$. Notice that $\varphi(q)$ is an even integer number for $q \geq 3$. 

Being an abelian group, all its Irreducible Representations are one-dimensional and coincide with their characters. A Dirichlet character $\chi$ of modulus $q$ is an arithmetic function from the finite abelian group $(\mathbb{Z}/q\mathbb{Z})^*$ onto $\mathbb{C}$ satisfying the following properties: 
\begin{enumerate}
\item 
$\chi(m+q) \,=\, \chi(m) $. 
\item 
$\chi(1) =1 $ and $\chi(0) = 0$. 
\item 
$\chi( n \, m ) \,=\, \chi(n) \, \chi(m)$. 
\item 
$\chi(m) = 0$ if $(m,q) > 1$ and $\chi(m) \neq 0$ if $(m,q) =1$. 
\item 
If $(m,q) =1$ then $(\chi(m))^{\varphi(q)} =1$, namely $\chi(m)$ have to be $\varphi(q)$-roots of unity. 
\item 
If $\chi$ is a Dirichlet character so is its  complex conjugate $\overline\chi$. 
\end{enumerate}
From property $5$, it follows that for a given modulus $q$ there are $\varphi(q)$ distinct Dirichlet characters that can be labeled as $\chi_{j}$ where $j = 1, 2, . . . , \varphi(q)$ denotes an arbitrary ordering. We will not  display the arbitrary index $j$ in $\chi_j$, except for explicit examples. For values of $m$ coprime with $q$, the character $\chi(m)$ mod $q$ may have a period less than $q$. If this is the case, $\chi$ will be called a {\em non-primitive} character, otherwise  $\chi$ is {\em primitive}. Obviously if $q$ is a prime number, then every character mod $q$ is primitive. 

Moreover, it is noteworthy that there is an important difference between {\em principal} versus {\em non-principal} characters.

\vspace{3mm}
\noindent
{\bf Principal character}. It is important to notice that for any $q$ there always exists  the {\em principal} character, usually denoted $\chi_1$, defined as 
\beq
\chi_1(m) \,=\, 
\left\{ 
\begin{array}{cl}
 1 & \, {\rm if} \,\,(m, q) = 1 \\
 0 & \, \, {\rm otherwise} 
 \end{array}
 \right.
 \label{principalcharacter}
 \eeq
The principal characters take only  the values $1$ or $0$ and satisfy 
\beq
\sum_{m=1}^{q-1} \chi_1(m) \,=\,\varphi(q) \,\neq \,0 \,\,\,.
\label{summm}
\eeq
Notice that, when $q = 1$, we have only the {\em trivial} principal character $\chi(m) = 1$ for every $m$. This case corresponds to the Riemann zeta function and therefore this remark clarifies why the Riemann zeta function is just a particular case of the $L$-functions. 

\vspace{3mm}
\noindent
{\bf Non-principal characters}. Contrary to the principal characters, which are made of the real numbers $0$ and $1$, the non-principal characters are in general  
complex numbers expressed in terms of phases $\theta_m$ 
\beq
\chi (m) \,=\,  e^{i \theta_m} 
\eeq
related to the $\varphi(q)$ roots of unity. These non-principal characters satisfy
\beq
\sum_{m=1}^{q-1} \chi(m) \,=\,0 \,\,\,.
\label{sumtozero}
\eeq
We will see below that the different results associated to the two sums given above, eqs.\,(\ref{summm}) and (\ref{sumtozero}), 
have far-reaching consequences on the analytic structure of the corresponding $L$-functions. 

\subsection{$L$-functions of principal characters and  the Riemann $\zeta$ function}\label{pppLLL}
 Notice that the principal character of modulus $q$ satisfies 
eq.\,(\ref{principalcharacter}) and therefore the relative $L$-functions can be expressed as 
\beq
 \boxed{
L(s,\chi_1) \,=\,\prod_{p \nmid q} \left(1 - \frac{1}{p^s}\right)^{-1} \,=\, \zeta(s) \, \prod_{p\mid q} 
\left(1 - \frac{1}{p^s}\right)^{-1}}\,\,\,,
\label{identitylprincipal}
\eeq
where $d | n$ denotes the integer $d$ which divides the integer $n$, and $d \nmid n$ otherwise. Since the finite product involving the primes which divide $q$ in the the right hand side never vanish, the zeros of the Dirichlet $L$-functions of principal characters coincide exactly with the zeros of the Riemann $\zeta$ function. Hence, establishing the GRH for these functions is equivalent to prove the original RH for the $\zeta$ function.

\subsection{Analytic structure of the $L$-functions}  As previously mentioned, there is an important distinction between the $L$-functions based on non-principal verses  principal characters which will be very
important for our purposes.  
\begin{itemize}
\item
The $L$ functions for  non-principal characters are {\em entire} functions,  i.e. analytic everywhere in the complex plane with no poles. 
\item The $L$-functions  $L(s,\chi_1)$ for  principal characters, on the contrary, are analytic everywhere except for a {\em simple pole} at $s=1$ with residue $\varphi(q)/q$. 
\end{itemize}
To show this result, let us  first  express any $L$-function in terms of a {\em finite} linear combination of the Hurwitz zeta function defined by the series  
\beq
\zeta(s,a) \,=\, \sum_{m=0}^\infty \frac{1}{(m+a)^s} \,\,\,,
\label{Hurwitz}
\eeq
whose domain of convergence is ${\mathbb Re} \,s >1$. Since we can split any integer $m$  as  
$$
m \,=\, q\, k +r \,\,\,\,\,\,\,,\,\,\,\,\,\,\,\,
{\rm where} \,\, 1 \leq r \leq q \,\,\,
{\rm and} \,\,\, k=0,1,2, \ldots 
$$
we have 
\begin{eqnarray}
L(s,\chi) &\,=\,& \sum_{m=1}^\infty \frac{\chi(m)}{m^s} \,=\, 
\sum_{r=1}^q \sum_{k=0}^\infty \frac{\chi(q k + r)}{(q k +r)^s} \,=\, 
\frac{1}{q^s} \, \sum_{r=1}^q \chi(r) \, \sum_{k=0}^\infty 
\frac{1}{\left(k+\frac{r}{q}\right)^s} \\
&=& \frac{1}{q^s} \, \sum_{r=1}^q \chi(r) \, \zeta\left(s, \frac{r}{q}\right) 
\,\,\,.\nonumber
\label{LHurwitz}
\end{eqnarray} 
The Hurwitz $\zeta$-function has a simple pole at $s=1$ with residue 1 and therefore the residue at this pole of the $L$-function is 
\beq
{\rm Res} \, L(s,\chi) \,=\, \frac{1}{q} \sum_{r=0}^q \chi(r)\,=\, 
\left\{
\begin{array}{cll}
\frac{\varphi(q)}{q} & & {\rm if }\, \, \chi = \chi_1 \\
0 & & {\rm if }\,\, \chi \neq \chi_1 \,\,\,.
\end{array}
\right.
\eeq

\subsection{Functional equation} The $L$-functions associated to the primitive characters satisfy a functional equation similar to  that of the Riemann $\zeta$-function. This functional equation strongly constrains the position of the zeros of these functions. To express such a functional equation, let's define the index $a$ as 
\beq\label{order}
a \equiv \begin{cases}1 \qquad
\mbox{if $\chi(-1)= -1$ ~~\,\,\,\,(odd)} \\
0 \qquad \mbox{if $\chi(-1) = \,\,\,\,\,1$  ~~~~(even)}
\end{cases}
\eeq
Moreover let's also introduce the Gauss sum 
\beq
G(\chi) \,=\,\sum_{m=1}^q \chi(m) \, e^{2 \pi i m/q}
\,\,\,,
\label{Gausssum}
\eeq
which satisfies $|G(\chi)|^2 = q$ if and only if the character $\chi$ is primitive. With these definitions, the functional equation for the primitive $L$-functions can be written as 
\beq
\label{FELambda}
L(1-s, \chi) \,=\, i^{-a} \,\frac{q^{s-1} \, \Gamma(s)}{(2\pi)^s} 
\, G(\chi) \, \left\{\begin{array}{c}\cos(\pi s/2)\\ \sin(\pi s/2) 
\end{array}\right\} \, 
L(s, \overline\chi)\,\,\,. 
\eeq
where the choice of cosine or sine depends upon the sign of $\chi(-1) = \pm 1$.  An equivalent but a more symmetric version of the functional equation (\ref{FELambda}) can be given in terms of the so-called {\em completed $L$-function} $\tilde L(s,\chi)$ defined by  
\beq
\tilde L(s,\chi) \equiv \left(\frac{q}{\pi}\right)^{(s+\delta)/2} \, 
\Gamma\left(\frac{s+\delta}{2}\right) \, L(s,\chi)\,\,, 
\label{completedLfunction}
\eeq
where $\delta = \frac{1}{2} (1 - \chi(-1))$. The completed $L$-function satisfies the functional equation 
\beq
 \boxed{
\tilde L(s,\chi) \,=\, \epsilon(\chi) \, \tilde L(1-s, \overline\chi) }\,\,\,,
\label{equivFunctional}
\eeq
where the quantity $\epsilon(\chi)$ 
\beq  
\epsilon(\chi) \,=\, \frac{G(\chi)}{i^\delta \sqrt{q}} \,\,\,,
\label{epsilonchi}
\eeq
is a constant of absolute value $1$. For the Riemann zeta-function $\zeta(s)$, the functional equation can be expressed in terms of the function $\xi(s)$ 
\beq
 \boxed{
\xi(s) \,=\, \xi(1-s)} \,\,\,,
\label{functionaeqriemann}
\eeq
where 
\beq
\xi(s) \,=\, \pi^{-s/2} (s-1) \,\Gamma\left(1+\frac{s}{2}\right) \zeta(s) \,\,\,.
\label{xiversuzeta}
\eeq
Notice that $\xi(s)$ is an {\em entire function} whose only zeros are inside the critical strip and, if the RH is true, all of them are on the critical axis ${\mathbb Re} \,s =1/2$. 

\subsection{Trivial Zeros} Using the Euler product representation of the $L$-function it is easy to see that these functions have no zeros in the half-plane ${\mathbb Re} (s) > 1$, in particular $\log L(s,\chi)$ is finite in this region since the series converges there.  Examining the functional equation (\ref{FELambda}) one sees that, analogously to the Riemann $\zeta$-function, the trivial zeros of the $L$-functions are those in correspondence with the zeros of the trigonometric functions present in the expression. Therefore 
\begin{enumerate}
\item If $\chi(-1) = 1$, then the trivial zeros are along the negative real axis located at $\sigma = - 2 k$, with $k=0,1,2,\ldots$. This is also the case of the Riemann zeta-function.
\item If $\chi(-1) = -1$, then the trivial zeros are along the negative real axis  but  now located at $\sigma = - 2 k -1$, with $k=0,1,2,\ldots$.
\end{enumerate} 

\subsection{Non-trivial Zeros and Generalized Riemann Hypothesis}
All other non-trivial zeros of the $L$-functions must lie in the critical strip $ 0 < \sigma < 1$. First of all, 
it is known that there are an infinite number of zeros in the critical strip and, to leading order,  the number of them with ordinate $0< t< T$ is given by\footnote{This result holds for $L$-functions relative to 
primitive characters mod $q$.}   
\beq
N(T,\chi) \,=\,\frac{T}{\pi} \,\log\frac{q \,T}{2 \pi e} + {\mathcal O}(\log (q T)) 
\,\,\,.
\label{asymptoticzeross}
\eeq
When the character is real (as it is also the case for the Riemann $\zeta(s)$) if 
\beq 
\rho_* = \sigma + i t
\label{rhoo}
\eeq 
is a zero of $L(s,\chi)$ then, from the duality relation (\ref{equivFunctional}) 
\beq
\hat\rho_* = (1- \sigma) - i t
\label{rhooo}
\eeq
 is also a zero of the same $L$-function. Hence, if $\sigma = 1/2$, the two zeros are then complex conjugates of  each other.  When the character $\chi$ is instead complex, a zero $\rho_*$ as in (\ref{rhoo}) corresponds to a zero $\hat\rho_*$ as in (\ref{rhooo}) 
 of $L(s,\overline\chi)$: in this case, if $\sigma =1/2$, the zeros of the $L$-functions associated to complex characters are not necessarily complex conjugates, since they refer to different characters. 

\begin{figure}[t]
\centering
\includegraphics[width=0.75\textwidth]{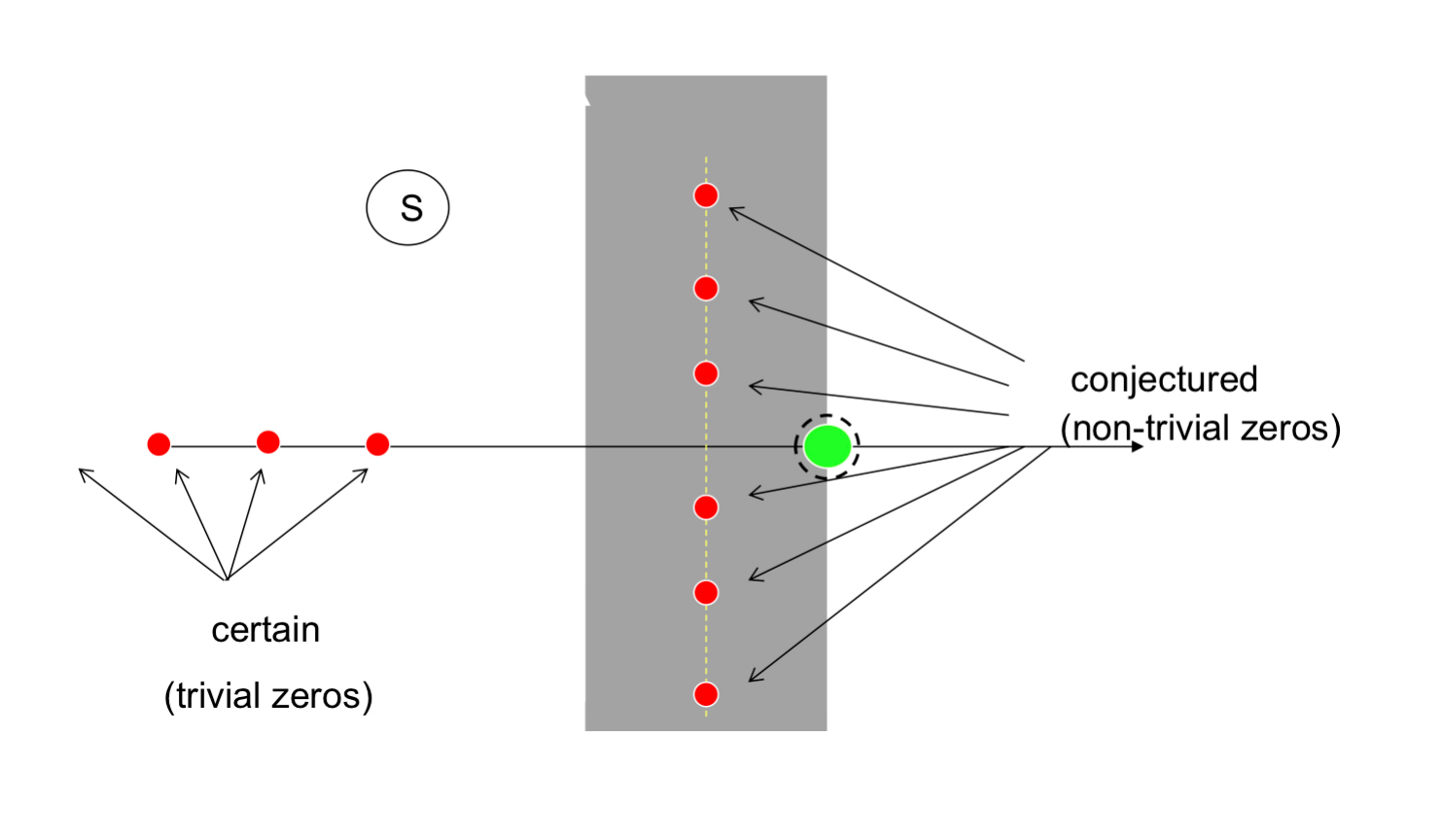}
\caption{Analytic structure of the $L$-functions, where the green circle at $s=1$ stays for a simple pole ($L$-functions corresponding to principal charcater) 
, while the dashed black circle at $s=1$ stays for the absence of this pole ($L$-functions corresponding to non-principal characters).
} 
\label{AnalyticstructureL-functions11}
\end{figure}

According to the Generalized Riemann Hypothesis, all non-trivial zeros of the primitive\footnote{It is important to refer to primitive characters in order to exclude the zeros of the factors $\prod_{p | {\hat q}} \left(1 - \chi(p) \,p^{-s}\right)$ present in the non-primitive characters which are all along the line $\sigma = 0$.} $L$-functions lie on the critical line $ \sigma = \half$, i.e. they have the form 
\beq
\rho_n \,=\, \frac{1}{2}  + i \gamma_n \,\,\,.
\label{conjecturedzeros}
\eeq
The conjectured analytic situation is summarised in Figure \ref{AnalyticstructureL-functions11}.  An explicit formula for the $n$-th zero of the Riemann zeta-function (and Dirichlet $L$-functions as well) as the solution of a transcendental equation was proposed in \cite{Transcendental}. Moreover, as followed by eq.\,(\ref{asymptoticzeross}), the imaginary part of these zeros for the Riemann zeta-function is expected to 
scale as 
\beq
\gamma_n \sim \frac{2 \pi n}{\log n} \,\,\,. 
\label{scalingimag}
\eeq
A lot is known about the non-trivial zeros of the $L$-functions.  Concerning the zeros of the Riemann function and therefore of all the $L$-functions with principal character, among many theoretical results (see \cite{Riemann,Edwards,Titcmarsh,Davenport,Bombieri,Sarnak,Conrey,Polya,Borwein,Broughan,Apostol,Iwaniec,Bombieri2,Steuding,Sarnak2,BK0,BK00,BK000,BK1,BK2,BK3,BK4,Bost,Connes,KeatingSnaith,Sierra1,Sierra,Sierra2,
Srednicki,Bender,reviewRiemann,Wolf,GriffinZagier,RodgersTao} for more details)   it is worth stressing that in 1914 G.H. Hardy \cite{HardyRiemann}  proved that there are infinitely many zeros along the critical line 
$ {\mathbb Re} \,s = 1/2$. Notice that this result does not imply that {\em all} the zeros are on the critical line. In 1974 N. Levinson \cite{Levinson} showed that more 
than one-third of the zeros of the Riemann zeta-function are on  the critical line, a bound further improved in 1989 by B. Conrey \cite{ConreyRiemann} who proved that at least two-fifths of the zeros of this function are on the critical line. These results have also been obtained for the generic $L$-functions, see \cite{Selberg1,Fujii,IwaniecS,Hughes,Conrey2}. 
Interestingly enough, it is also possible to prove that most of the nontrivial zeros of $\zeta(s)$ and also of the $L$-functions cannot lie too far from the critical line $\sigma = 1/2$, an observation due to Bohr and Landau \cite{BohrLandau}, and also Littlewood \cite{Littlewoodzeros}. The derivation of this theorem, an example of the so-called {\em density theorems}, is shown in detail in \cite{Steuding}: denoting by $N(\sigma,T,\chi)$ the number of zeros of the $L$-function of character $\chi$ with ordinate $0< t< T$ but at 
$\sigma \neq 1/2$, it holds   
\beq
N(\sigma, T,\chi) \,=\, o(N(T,\chi)) \,\,\,,
\label{density}
\eeq
as $T$ tends to infinity, i.e. all but an infinitesimal proportion of the zeros of $L(s,\chi)$ lie in the strip $\frac{1}{2} - \epsilon < \sigma < \frac{1}{2} + 
\epsilon$,  no matter how small $\epsilon$ may be. 

There is also a long tradition of studies for the numerical determination of the Riemann zeros built on previous results such as \cite{Titchmarsh,Turing,Riele}. In 1992 
A. M. Odlyzko  \cite{Odlyzkozeros}  computed 175 million zeros of heights $T$ around $10^{20}$ not only to verify that they alined along 
the axis $ {\mathbb Re} \,s = 1/2$ but also to check the Montgomery-Dyson pair-correlation conjecture as well as other conjectures that state that the zeros of the Riemann zeta-function behave like eigenvalues of random matrices coming from the Gaussian Unitary Ensemble \cite{Dyson,Montgomery,Odlyzko}. The present record on the computation of the zeros is due to D. Platt and T. Trudgian \cite{Platt} who verified that the first $12.363.153.437.138$ zeros up to the height $T=3.000.175.332.800$ are all along the line $ {\mathbb Re} \,s = 1/2$.

\section{On the zeros of random function analogues of $L$-functions}\label{surprising}
There is a rather surprising result concerning the zeros of the Riemann $\zeta$ function (and all other Dirichlet  $L$-functions) in relation to the zeros 
of a family of random functions. These results are the content of two theorems due to Grosswald and Schnitzer \cite{Grosswald}.

\begin{theorem} (Grosswald and Schnitzer) 
\label{GrosswaldSchnitzer}  
Let $p_m$ be the $m$-th prime and select an integer number $p_m'$ so that  
\beq
p_m \leq p_m' \leq p_{m+1} \,\,\,.
\eeq
With these $p_m'$ form the infinite product 
\beq
\zeta'(s) \,=\, \prod_{m=1}^\infty \left(1 -\inv{ (p_m')^{s}}\right)^{-1} \,\,\,.
\eeq
The function $\zeta'(s)$ possesses the following properties: (i) $\zeta'(s) \neq 0$ for $\sigma > 1$; (ii) $\zeta'(s)$ can be continued as a meromorphic function in $\sigma >0$; 
(iii) in $\sigma > 0$, $\zeta'(s)$ has a single pole at $s=1$ with residue $r$, $1/2 \leq r \leq 1$; (iv) in $\sigma >0$, $\zeta'(s)$ has the {\em same} zeros with the same multiplicities of the Riemann zeta function $\zeta(s)$. 
\end{theorem}

\begin{theorem} (Grosswald and Schnitzer) 
\label{GrosswaldSchnitzer1}
Let $L(s,\chi)$ be the Dirichlet $L$ function based on any  Dirichlet character of modulus 
$q$.  Let $\Primes = \{p_1, p_2, \ldots \}$  denote the set of primes while $\Primesp = \{ p'_1, p'_2, \ldots \}$ a set of  integers $p'_n$ satisfying 
\beq
 \label{ppn}
 p_n \leq p'_n < p_n + K, ~~~~~~p'_n = p_n ~~ ({\rm mod} ~ q)
\eeq
where $K \geq q$ is an arbitrary integer,  and define the modified $L$-function according to the infinite product 
\beq
\label{Lp}
L' (s, \chi) \,=\,  
\prod_{n=1}^\infty  \( 1 -  \frac{\chi(p'_n)}{(p'_n)^s} \)^{-1} \,\,\,.
\eeq
Then $L'(s,\chi) $ can be analytically continued to the half plane $\sigma > 0$ and, in this domain, it has the same zeros as the Dirichlet $L$-function $L(s,\chi)$. Moreover, if $\chi$ is a non-principal character then $L'(s, \chi)$ has no poles for $\sigma>0$, as does $L(s,\chi)$. 
\end{theorem}

\begin{figure}[b]
\centering\includegraphics[width=.5\textwidth]{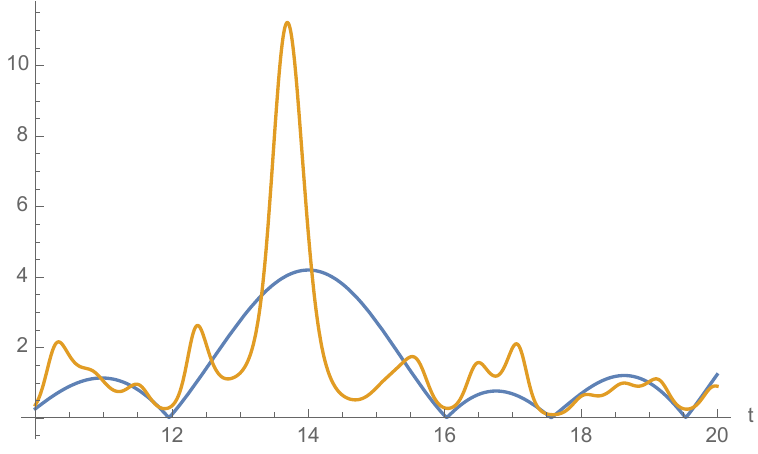}
\caption{Numerical illustration of Theorem \ref{GrosswaldSchnitzer} for the non-principal character $\chi$ mod $5$ indicated in \eqref{char}.    The blue line is $|L(\half + i t)|$ where $t$ is the x-axis.  
The erratic orange line is $|L'(\half + i t)|$ for some randomly chosen state $\Primesp$.  We chose $N =5 \times 10^3$ and $K=Mq=10$  for $M =2$.}
\label{LpDir}
\end{figure}

{\bf Remark:} What is surprising in these theorems is the emergence of the following scenario: if in the Euler product representation of the $L$-functions or the Riemann zeta-function we use another set of random numbers which shares with the primes the same residue wrt the modulus $q$ and the same rate of growth, then all the non-trivial zeros of the original $L$-functions or the Riemann zeta function remain  exactly at the same location  in the critical strip!  In particular, Theorem \ref{GrosswaldSchnitzer} suggests that the validity of the Generalized Riemann Hypothesis may not depend on the detailed properties of the primes and this further justifies the probabilistic considerations presented later in this paper.
Figure  \ref{LpDir} shows the results of a numerical check of Theorem \ref{GrosswaldSchnitzer1}, where we have chosen as example the non-principal character $\chi$ mod $5$, whose values in the first period are given by 
\beq
\label{char}
\{ \chi (1), \chi(2),  \chi (3), \chi(4), \chi(5) \} = \{ 1, -1, -1, 1, 0 \} \,\,\,, 
\eeq
In Figure  \ref{LpDir} we plot $|L(\half+it)|$ and $|L'(\half+i t)|$  for a randomly chosen set of the integers $p_n'$ as a function of $t$ in the region of  the first 3 zeros. Whereas $|L'(\half+i t)|$ is erratic due to the randomness of the integers $p'_n$ and changes its shape if we change the set of these random numbers, the  validity of Theorem \ref{GrosswaldSchnitzer} is nevertheless clear, i.e.  the two functions share the same zeros.  

Amazingly enough, there is another theorem due to Chernoff that leads to a completely different outcome if, in the infinite product representation of the Riemann zeta-function, we substitute the primes $p_n$ with their average behavior $p_n \sim n \log n$, namely

\begin{theorem} ({\em Chernoff}) \label{Chernoff} \cite{Chernoff}.  
Consider the Euler infinite product representation of the Riemann $\zeta$-function. Substitute the primes $p_n$ in such a formula with their  approximation $p_n \sim n \log n$ and define the modified  function $\zeta"(s)$ according to the infinite product 
\beq 
\zeta''(s) \,=\,\prod_{n=1}^\infty  \( 1 -  \frac{1}{(n\,\log n)^s} \)^{-1}\,\,\,.
\label{CsC}
\eeq 
The  function $\zeta''(s)$can be analytically continued into the half-plane $\Re(s) > 0$ except for an isolated singularity at $s=0$.  Furthermore it no longer has any zeros in this region. 
\end{theorem}

\vspace{1mm}
The first  two theorems will help us in shaping our later considerations, although our constructions are based on the true prime numbers.  

\section{Quantum vs Statistical}

There are many different ways to formulate the RH and a very good selection of these different formulations of the problem can be found e.g. in \cite{Borwein,Broughan}. Here we focus our attention on two of these approaches which are closer than the others to a physicist's background and sensibility. The first we call hereafter 
the {\em quantum approach}, the second one the {\em statistical approach}.  

\subsection{ Quantum Approach}  
This approach has attracted  for many years the attention of  theoretical physicists, for the simple reason that it is deeply related to the spectral theory of quantum mechanics. Originally stated by P\'olya and Hilbert around 1910,  this approach has given rise to an important series of works on the Riemann  $\zeta$-function 
by Berry, Keating, Bost and Connes, Sierra, Srednicki, Bender and many others \cite{BK0,BK00,BK000,BK1,BK2,BK3,BK4,Bost,Connes,KeatingSnaith,Sierra1,Sierra,Sierra2,Srednicki,Bender} (for a more complete list of references, see the reviews \cite{reviewRiemann,Wolf}).  In a nutshell, this approach can be iconically rephrased as 
\beq
\xi\left(\frac{1}{2} + i \hat H_R\right) \,=\, 0 \,\,\,,
\label{Polyaeq}
\eeq
where the function $\xi(s)$ is defined in (\ref{xiversuzeta}) while $\hat H_R$ is a hermitian operator whose spectrum coincides with the imaginary parts $\gamma_n$ of the non-trivial zeros of the Riemann zeta-function along the 1/2 axis  
\beq
\hat H_R \,|\psi_n\rangle \,=\, \gamma_n \, |\psi_n\rangle \,\,\,.
\label{polyaspectrum}
\eeq
Finding such a hermitian operator (i.e. a quantum Hamiltonian) has been the focus of an intense research activity for decades. The task is particularly challenging in view of the 
random properties exhibited by the known zeros of the Riemann zeta function, which behave like eigenvalues of large random Hermitian matrices \cite{Dyson,Montgomery,Odlyzko}. We refer to this quantum approach as the {\em vertical approach} because it aims at proving that all the non-trivial zeros of the Riemann zeta-function are vertically  aligned along the axis 1/2 on the basis of 
the spectral property of hermitian operators.

\vspace{3mm}
\noindent
{\bf Spectral determinant}. 
It is better to stress the key point of this approach and the true meaning of eq.\,(\ref{Polyaeq}): the function $\xi(s)$ is an entire function whose only zeros are within the critical strip and admits an infinite product representation \cite{Edwards} 
\beq
\xi(s) \,=\, \xi(0) \, \prod_{\rho} \left(1 - \frac{s}{\rho}\right) \,\,\,,
\label{xiproduct}
\eeq
where $\rho$ ranges over all the roots $\rho$ of $\xi(\rho) = 0$ and the infinite product is understood to be taken in an order which pairs each root $\rho$ with the corresponding $1-\rho$. Hence, the meaning of eq.\,(\ref{Polyaeq}) is to find a non-trivial  hermitian operator $H_R$ such that $\xi(s)$ emerges as its {\em exact} spectral determinant {\em without} assuming the validity of the RH. In fact, on the contrary, assuming the RH there is a tautological way of solving eq.\,(\ref{Polyaeq}) and showing the existence of a quantum mechanical Hamiltonian $H_R$ which satisfies (\ref{Polyaeq}). Indeed, given the scaling behaviour of the $\gamma$'s given in (\ref{scalingimag}), one can either find such a Hamiltonian using the semi-classical method previously applied to find the potential of the prime numbers \cite{primenumberpotential} (as it has been done in \cite{zerosemiclassical}) or, using methods of supersymmetric quantum mechanics \cite{susyqm}. In both cases, denoting by $|\gamma_n\rangle$ the eigenfunction of such a 
Hamiltonian, and assuming the RH one can explicitly show that  
\beq
H_R \,=\, \sum_n \gamma_n\, | \gamma_n \rangle \, \langle \gamma_n| 
\eeq
satisfies eq.\,(\ref{Polyaeq}) but, unfortunately, nothing is learned about  the RH by this tautological construction.  

\vspace{3mm}
\noindent
In summary, the quantum approach is a very appealing way to prove the RH but so far the sought after Hamiltonian has remained elusive.

\subsection{Statistical Approach}  

\label{StatApp}
As shown originally in the papers \cite{ML,LM} (elaborating on previous approaches of the same type discussed in \cite{EPFchi, Franca1}), the GRH can be 
addressed in a different way. The starting point of this new approach comes from a simple remark: if all the infinitely many Dirichlet $L$-functions have their non-trivial zeros along the axis $\sigma = \half$, behind this fact there should be some universal and very robust reason which transcends the details of the characters entering their definition and which rely instead on some of the general properties of these quantities. Such a reason can be nailed down to the existence of a random walk and its diffusive universal scaling law $N^{1/2+ \epsilon}$ after $N$ steps, in the sense that the value $\sigma = \half$ can be identified with the critical exponent of a random walk process which exists for all these functions. 

\begin{figure}[t]
\centering\includegraphics[width=.5\textwidth]{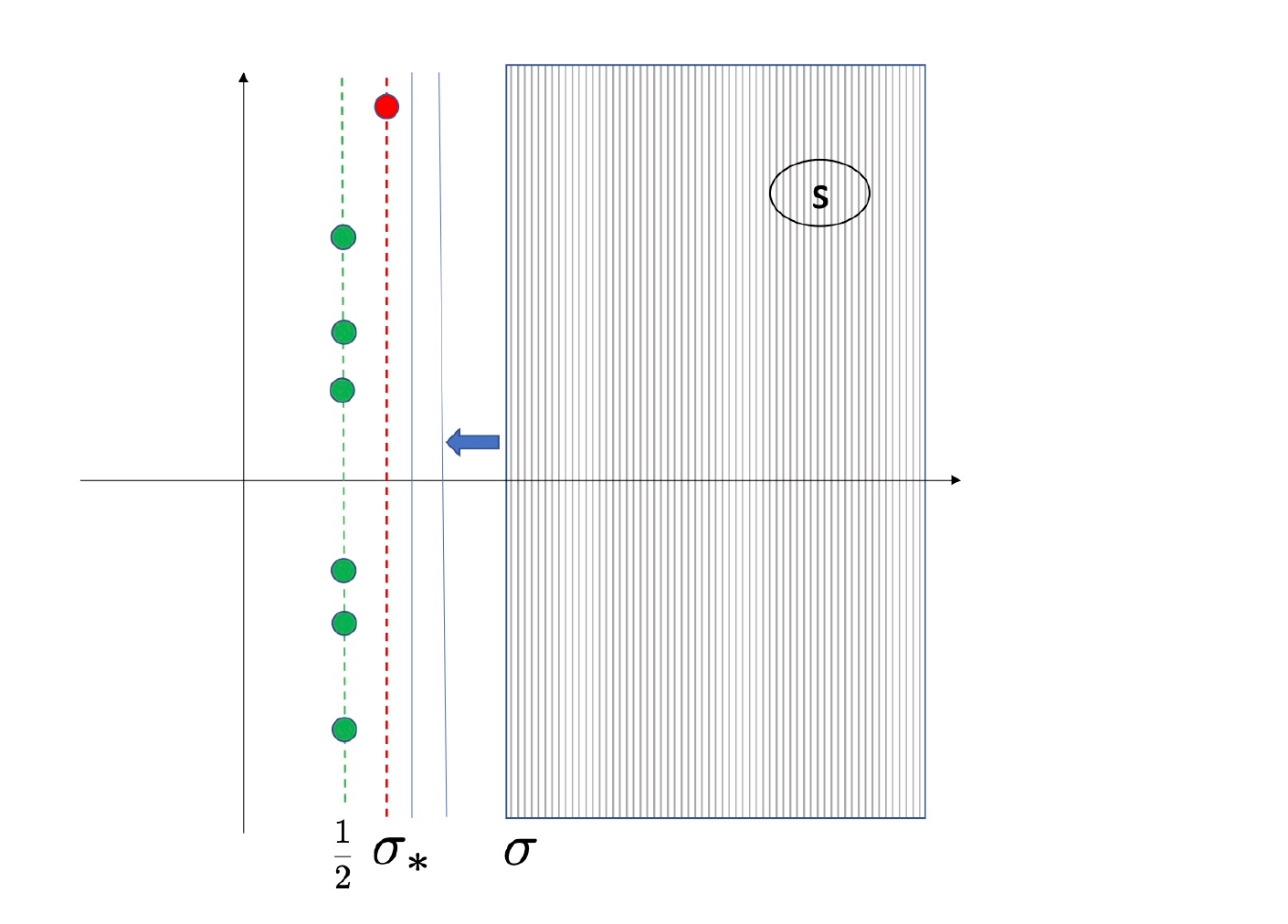}
\caption{In the statistical approach the problem is to find the abscissa $\sigma_*$ of convergence of a Dirichlet series  $D(s)$ originally defined for $\sigma \geq 1$. 
For the Riemann zeta-function, such a critical value $\sigma_*$ is dictated by the abscissa of the zero more distant from the origin, and the same for a generic 
$L$-function. If $\sigma_* = 1/2$, then {\em all} zeros are on the critical axis.}
\label{horizontalgraph}
\end{figure}

From the mathematical point of view, we will see that in this statistical approach the problem consists in general to determine the abscissa of convergence $\sigma_*$ 
of a inverse Mellin transform $D(s)$ of a weight function $\mathcal{W}(x)$ such as   
\beq
D(s) \,=\,\int_{1}^\infty \frac{\mathcal{W}(x)}{x^{s+1}} \,dx \,\,\,. 
\label{genericDirichseries}
\eeq
The specific nature of the weight function ${\mathcal W}(x)$ will change according to the case under scrutiny, i.e. whether we consider the Riemann zeta function or the $L$-functions of non-principal characters. In both cases, however, the corresponding $D(s)$ will have an original domain of convergence given by $\sigma \geq 1$ and the goal will consist of showing that it can be extended down 
 to $\sigma_* = 1/2$. Clearly $\sigma_*$ cannot be less than $1/2$, since we know that there are infinitely many zeros on the axis $1/2$, although we do not know if they are {\em all} on this axis. However, using the duality properties of the Riemann and Dirichlet functions, the RH and the GRH could  be proved to be true as far as we are able to show that the abscissa of convergence of the relative inverse Mellin transform (\ref{genericDirichseries}) associated to the Riemann zeta function or to the $L$-functions is precisely $\sigma_* = 1/2$.  This approach can be  called the {\em horizontal approach} since its aim is to establish how far we can move towards the origin the abscissa of the half-plane  domain of convergence of the function $D(s)$, as shown in Figure \ref{horizontalgraph}.  

So far it was important to distinguish two cases: 
\begin{enumerate}
\item $L$-functions relative to non principal characters. 
\item $L$-functions  relative to principal characters.
\end{enumerate}
These cases gave rise to inverse Mellin transform of two different origins. Let's first discuss how these two cases come about and, later, how it can be argued on their unification 
under a single mathematical umbrella.

\begin{enumerate}
\item {\bf  $L$-functions relative to non principal characters.} For these functions, the corresponding expression $D(s)$ originates from the infinite product representation on the {\em prime numbers} of the $L$ functions. Indeed, we know that the values of the characters $\chi(n)$ of modulus
$q$ (see eq.\,(\ref{Lfunctions})) are phases which are related to the $q$-roots of  unity. This means that, in this case, we are dealing with a dice of $q$ faces, as becoming immediately evident taking the logarithm of the infinite-product representation of these functions, where we have \cite{EPFchi}
\beq 
\log L(s,\chi) \,=\,  X(s,\chi) + R(s,\chi)\,\,\,,
\eeq
where
\beq\label{PDir}
X(s,\chi) = \sum_{n=1}^{\infty} \dfrac{\chi(p_n)}{p_n^{\, s}}\,\,\,\,\,\,\, , \qquad
R(s,\chi) = \sum_{n=1}^{\infty} \sum_{m=2}^{\infty}
\dfrac{\chi(p_n)^{m}}{m p_n^{\, ms}}\,\,\,.
\eeq
Since $R(s,\chi)$ is absolutely convergent for $\sigma > \half$, the possibility to enlarge the convergence of the original Euler product depends only on properties of $X(s,\chi)$ and 
therefore we can write 
\beq
\label{logsum2}
\log L(s,\chi) \, \sim\,X(s,\chi) + O(1) \,\,\,.
\eeq
The singularities of $ \log L(s,\chi)$ are determined by the zeros and poles of $L(s,\chi)$ but, since the $L$-functions of non-principal characters do not have poles, $X(s,\chi)$ {\em is the diagnostic quantity which directly locates their non-trivial zeros}. Taking now the real part\footnote{Analogous arguments apply to the imaginary part of $X(s, \chi)$.} of $X(s, \chi)$ in \eqref{PDir}, i.e. 
$S(\sigma, t,\chi) = Re (X(s,\chi))$, and focusing on its expression at $t=0$ (for the half-line convergence of this kind of series), we end up  considering  this series 
\begin{equation}
\label{SDef}
S(\sigma,\chi)\,=\, \sum_{n=1}^\infty \dfrac{\cos( \theta_{p_n})}{p_n^{\, \sigma}} \,\,\,. 
\end{equation}
Defining
\begin{equation}
\label{CxDef}
B(x, \chi) \,=\, \sum_{p \le x} \cos\( \theta_p\)\,\,\,,  
\end{equation}
we have  
\beq
B(p_n;,\chi) - B(p_{n-1}, \chi) \,=\, \cos (\theta_{p_n})\,\,\,,
\eeq
and then  
\begin{equation}
S(\sigma, \chi) = \sum_{n=1}^{\infty} B(p_n,  \chi) \left( \dfrac{1}{p_n^{\,\sigma}} - 
\dfrac{1}{p_{n+1}^{\,\sigma}}\right) = \sigma \sum_{n=1}^{\infty} 
B(p_n ,  \chi) \int_{p_n}^{p_{n+1}} \dfrac{1}{u^{\sigma+1}} du\,\,\,.
\end{equation}
Given that $B(x, \chi) = B(p_n,  \chi)$ is a constant for $x \in (p_n, p_{n+1})$, we finally arrive to 
\begin{equation}
\label{importantintegral}
S(\sigma, \chi) = \sigma \int_{2}^{\infty} \dfrac{B(x,  \chi)}{x^{\sigma+1}} dx\,\,\,. 
\end{equation}
Hence, looking at eq.\,(\ref{genericDirichseries}), the function $S(\sigma, \chi)$ plays in this case the role of the aforementioned function $D(s)$ while Re$(s B(x,\chi))$ plays the role of the weight function $\mathcal{W}(x)$. 
Hence, the convergence of the integral is dictated by the behavior of the function $ B(x ; \chi)$ at $x \rightarrow \infty$: if $B(x ; \chi) = O(x^{\alpha})$ for $x \rightarrow \infty$, then the integral converges for $ \sigma > \alpha$ and diverges precisely at $\sigma =\alpha$. In Refs. \cite{ML,LM} it was argued that one should expect to have 
\beq
B(x,\chi) \sim x^{1/2+\epsilon} \,\,\,,
\label{scalingbhat}
\eeq
where $\epsilon$ is arbitrary and  strictly positive, on the basis of 
\begin{enumerate}
\item the Dirichlet theorem \cite{Dirichetresidues} on the equidistribution of the residues (mod q) along the sequence of the prime numbers, which implies the equi-distribution of the 
angles $\theta_{p_n}$ moving along the sequence of the prime numbers; 
\item the very weak correlation between a pair or more k-plet of successive angles $\theta_{p_n}, \ldots \theta_{p_{n+k}}$, as result of the analysis of Lemke Oliver and Soundararajan \cite{OliverSoundararajan} on the basis of the Hardy-Littlewood prime k-tuples conjecture.
\end{enumerate}
\item {\bf  $L$-function relative to principal characters.} For these $L$-functions, which are all proportional to the Riemann zeta function (see 
eq.\,(\ref{identitylprincipal})), we cannot  use the previous approach, because in this case there is not a q-plet of residues to play with: for the Riemann zeta-function, the values of the principal character is in fact identically equal to $\chi(n) = 1$. Hence, in this case, instead of studying the convergence of the infinite product representation on the 
Riemann function, it is much more useful to study the convergence of the inverse Mellin transform which originates from  the infinite series representation of the 
(multiplicative) inverse of the Riemann zeta function 
\beq
1/\zeta(s) \,=\, \sum_{n=1}^{\infty}\frac{\mu(n)}{n^s} \,=\,
s\, \int_1^{\infty} \frac{M(x)}{x^{s +1}} \, dx \,\,\,,
\label{Mellina}
\eeq
where $M(x)$ is the Mertens function, given by 
\beq
\label{Mertenfunction222}
M(x) \,=\, \sum_{m=1}^x \mu(m) \,\,\,.
\eeq
As discussed in the Introduction, this is indeed the main object of our study, further refined to be the restricted Mertens function $\hat M(x)$ based on the square-free numbers. 
 Therefore, if we are able to show that asymptotically $\hat M(x)$ goes as $\hat M(x) \sim x^{1/2 + \epsilon}$, for any arbitrarily small $\epsilon$, then the RH is indeed true, since the integral (\ref{Mellina}) diverges at $ {\mathbb Re} \,s = 1/2$, making clear the presence of a singularity on this axis. The duality of the Riemann zeta function, expressed by the eq.\,(\ref{functionaeqriemann}), finally ensures that {\em all} the zeros are along the axis 
$ {\mathbb Re} \,s = 1/2$. 
\end{enumerate}

Let us  finally return to the random $L$ functions in the Grosswald-Schnitzer theorems described in section III in light of the above results. For the  $\zeta$ case,  not much can be said since there is  no analog of Mertens function,  since the $p'_n$ don't lead to an arithmetic  M\"obius function $\mu'$ which is simple to handle.  However  for non-principal characters, we can make a connection. In Theorem 2 of section III   we are free to select only random $p'_n$ that are ordered:
$$ p'_1 < p'_2 <  p'_3  \ldots$$  of course still subject to  \eqref{ppn}.   In such a case,  the inverse Mellin transform goes through and we have
 $$ S'(s, \chi) =  s \int dx  ~ B'(x,\chi)/x^{s+1}$$
 where
 $$B'(x, \chi) = \sum_{n<x}  \chi (p'_n).$$ 
 The main point is that we argued that $B$ behaved like a random walk, but $B'$ is even {\em more random}, meaning with these words that in this case even the numbers $p_n'$ on which the characters are evaluated are random as well.  Thus this provides an even stronger argument that 
 the random Dirichlet functions $L'$,  which share the same zeros as $L$,  have no zeros to the right of the critical line.

 \section{Some known results about the Mertens function}

Let's briefly discuss some known facts about the Mertens function which are also helpful for better understanding 
our subsequent  analysis carried out in Parts  B and C. 

\subsection{Mean of the Mertens function}

The first result concerns the mean of the Mertens function: as shown in Appendix \ref{AA}, the average of the Mertens function vanishes 
\beq
\lim_{x\rightarrow \infty} \frac{1}{x} \sum_{m=1}^x \mu(m) \,=\, 0 \,\,\,. 
\label{averageMertens}
\eeq
This implies that our original sequence $\{{\mathcal S}_n\}$, given in eq.\,(\ref{sequencebinary}), has a perfect balance\footnote{
The arithmetic function $\mu(m)$ has values $0, \pm 1$.  Given that its average vanishes, the number of $+1$ must balance the number of $-1$: these are the two values 
present in the M\"{o}bius coefficients $\hat \mu(n)$ restricted to the square-free numbers and are those mapped to the values $0$ and $1$ of the sequence $\{\mathcal{S}_n\}$, 
as stated in eq.\,(\ref{relationmuands}).}  in the numbers of $0$'s and $1$'s.

An equivalent way to arrive to the same result is to return to the formula
\beq
\label{PNT1}
\inv{\zeta(s)} =  s \int_0^\infty  \frac{M(x)}{x^{s+1}}   ~ dx
\eeq
which was presented earlier. It is well known that the Prime Number Theorem (PNT) follows if $\zeta(s)$ has no zeros along the line ${\mathbb Re} \,s =1$ and, in fact, the first proof of the PNT was based on this 
connection.  If there are no zeros along this line, then $1/\zeta(s)$ is necessarily finite along this line. This implies that the above integral converges for ${\mathbb Re} \,s =1$. This is guaranteed if 
$M(x) = O(x^{1-\delta})$  where $\delta$ is small and positive. Since the PNT has been proven to be true, we can assume 
\beq
\label{PNT2}
M(x) = O(x^{1-\delta})
\eeq
The above result implies
\beq
\label{PNT3}
\lim_{x\to \infty}  \frac{M(x)}{x} =0
\eeq
which implies that the mean of $M(x)$ is zero.   
Of course,  $M(x) = O(x^{1/2+\epsilon})$ is a stronger bound than \eqref{PNT2},  and the PNT would follow from it.  

\subsection{Titchmarsh and Kotnik-van de Lune expansion} In his book on  the Riemann zeta-function \cite{ Titcmarsh}, Titchmarsh derived a trigonometric series in $\log x$ for the 
function $t(x)$, defined as
\beq
t(x) \,=\, \frac{M(x)}{\sqrt{x}}\,\,\,,
\label{definitionofq1(x)}
\eeq
based on the following assumptions 
\begin{itemize}
\item The Riemann hypothesis is true. 
\item All non-trivial zeros of the Riemann zeta-function are simple. 
\end{itemize}
A nice re-rewriting of Titchmarsh's trigonometric series for $t(x)$ has been given by Kotnik and van de Lune \cite{Kotnik} and in the following we use the expression found by these authors for 
the truncated version of this series for $t(x)$ up to a cut-off $K$ of terms, given by 
\beq
t_K(x) \,\sim \, 2 \sum_{k=1}^K |a_k| \,\cos(\gamma_k \,\log x + {\rm arg}\, a_k) \,\,\,,
\label{q_k}
\eeq
with   
\beq
a_k\, =\, (\rho_k \,\zeta'(\rho_k))^{-1} 
\hspace{5mm}
,
\hspace{5mm}
\rho_k = \frac{1}{2} + i \,\gamma_k
\label{tivan}
\eeq
\begin{flushleft}
\begin{figure}[t]
\hspace{-10mm}
\includegraphics[width=1.0\textwidth]{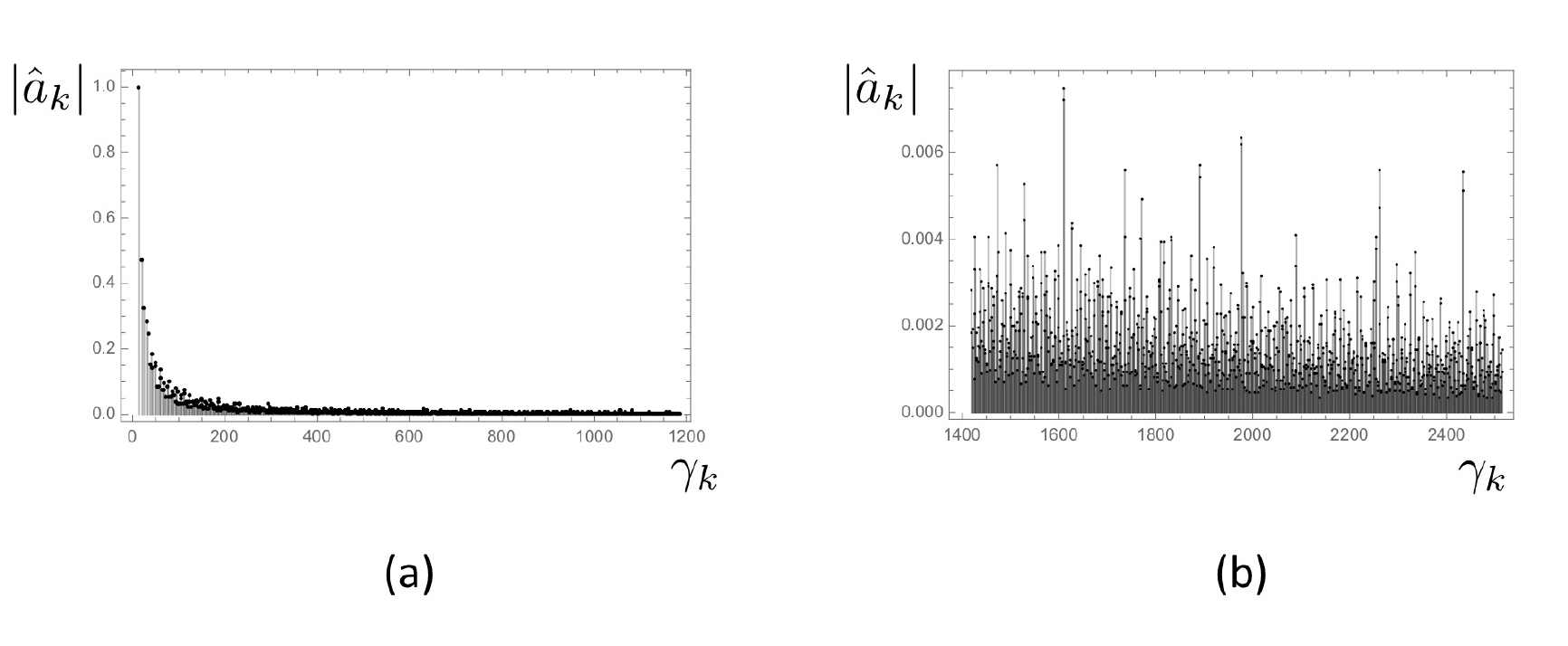}
\caption{(a) Left-hand side: ratios of the absolute Fourier coefficients $|\hat a_k| =|a_k/|a_1|$ of Titcmarsh's function versus the imaginary part of the Riemann zeta function zeros;
(b) Right-hand side: zoom on the ratios of the absolute Fourier coefficients $|\hat a_k|$  for $1400 <\gamma_k  < 2500$.} 
\label{Titchmarsh2}
\end{figure}
\end{flushleft}
\begin{figure}[b]
\centering
\includegraphics[width=0.7\textwidth]{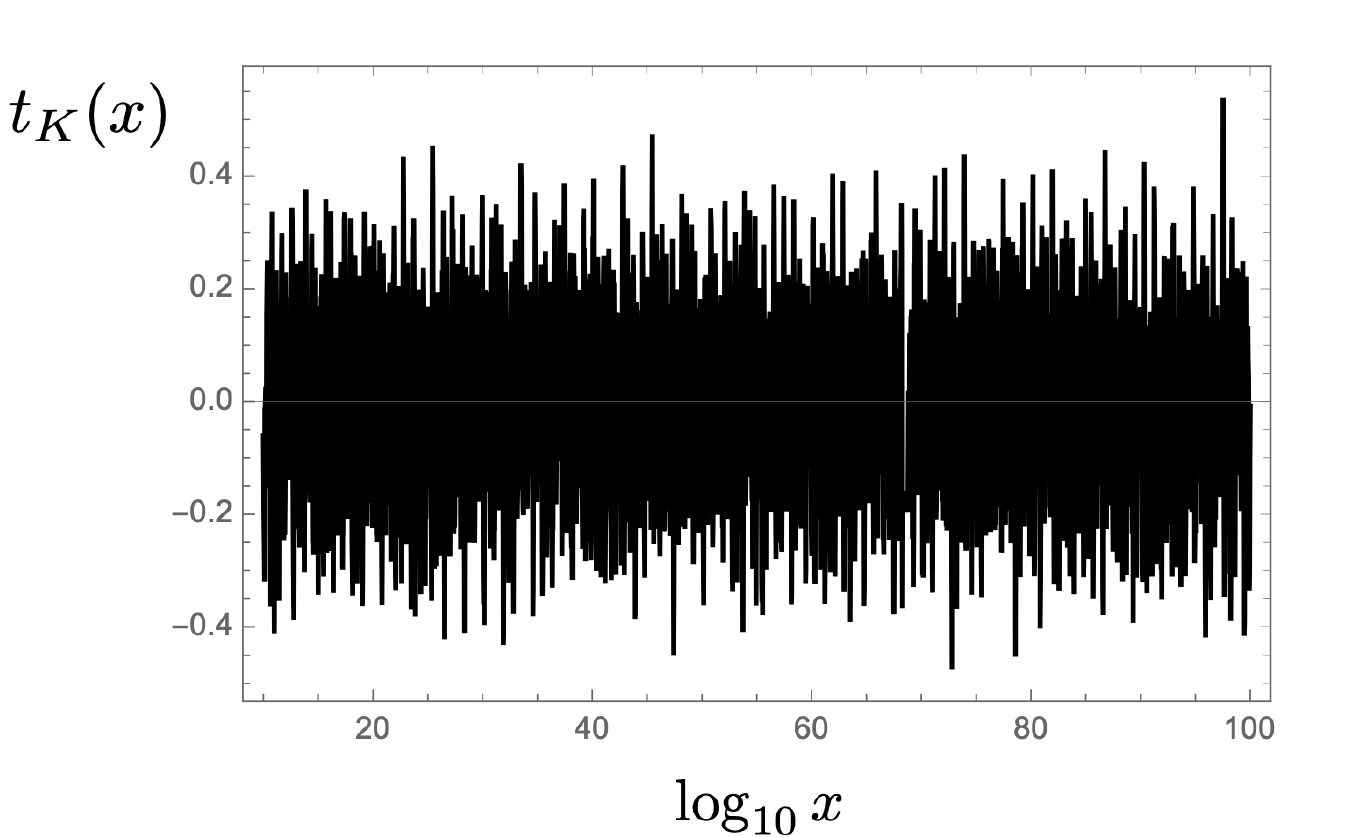}
\caption{Plot of Titcmarsh's function $t_K(x)$ in the range $1 \leq x  \leq 10^{100}$ with $K=1200$. } 
\label{Titchmarsh1}
\end{figure}
where $\rho_k$ denotes the $k$-th zero of the Riemann zeta function in the upper half-plane, counted in increasing order. The ratio $|\hat a_k| \equiv |a_k/a_1|$ of the coefficients of this series is shown in Figure \ref{Titchmarsh2}, from where one can see that the coefficients $a_n$ do not form a monotonically decreasing sequence but behave instead quite irregularly.
The plot of the function $t_K(x)$ (for $K=1200$) is shown in Figure \ref{Titchmarsh1} and, at the first sight, it looks to be the plot of a completely random function. From a celebrated theorem by Mark Kac 
\cite{Kac} this would be indeed the case for a trigonometric series as the one in (\ref{q_k}) (properly normalised though) if the frequencies $\gamma_n$'s were linearly independent over the rationals\footnote{ This seems to be just the  case. In other words, there are no known reasons why the $\gamma_n$'s should be instead linearly dependent over the rationals. The fact that, on a large scale, the $\gamma_n$'s, their gaps and their correlations are well described by random matrix theory \cite{Dyson, Montgomery, Odlyzko} gives an additional support to this hypothesis.}.  
However, in this case, since it is known that the series\footnote{The conjectured behaviour is 
$\sum_{0 < \gamma < T} |\rho \zeta'(\rho)|^{-1} \sim (\log T)^{5/4}$, see \cite{Ng}.}
\beq
\sum_{\rho} \frac{1}{|\rho \zeta'(\rho)|} 
\eeq
diverges, the sum of the coefficients of $\cos(\gamma_k \,\log x + {\rm arg}\, a_k)$ can be made arbitrarily large by choosing $K$ large enough. Hence, if one can find a value 
of $y=\log x$ such that all the arguments ($\gamma_k y + {\rm arg}\,a_k)$ of the cosines are close to the integer multiples of $2\pi$, the function $t_k(x)$ could become arbitrarily large. Indeed this condition could be satisfied if the $\gamma_n$'s would be linearly independent on the rationals. 
Our expectation, on the statistical analysis performed in the remaining parts of the paper, is however that $t_K(x)$ can go {\em at most} as $\sqrt{\log\log x}$. For further details on the  expansion (\ref{q_k}) and other properties of the Mertens function we refer the reader to the papers \cite{Ng,Kotnik,Riele2,Riele3,Pinz} and references therein.

\section{Random Dirichlet Series and Mertens Function}

For the purpose of determining the abscissa of convergence $\sigma_*$ of $1/\zeta(s)$ 
 (i.e. the location of its pole with the largest positive real part), it is interesting to see what happens for a {\em purely random} Dirichlet series \cite{Kahane}
\beq
\tilde R(s) \,=\, \sum_{n=1}^\infty \frac{a_n}{n^s} \,\,\,,
\label{kahareseries}
\eeq
with $a_n$ truly {\em independent random variables} taking values $\pm 1$. Once we express this function in terms of its inverse Mellin transform 
\beq
\tilde R(s) \,=\, s\, \int_1^\infty \frac{R(x)}{x^{s+1}}\,dx \,\,\,,
\label{kah}
\eeq
the corresponding function ${ R}(x)$ is given in this case by the canonical displacement sum of a pure random walk 
\beq
R(x) \,=\, \sum_{n=1}^ x a_n \,\,\,.
\eeq
Therefore, the probability distribution of the random variable 
\beq
\frac{R(n)}{\sqrt{n}} 
\eeq
is clearly in this case the standard normal distribution ${\mathcal N}_{0,1}$ and for the moments of the random variable $R(n)$, in the large $n$ limit, we have 
\beq
\langle (R(n))^{2k} \rangle \,=\, n^k \, (k-1)!! 
\hspace{3mm}
,
\hspace{3mm}
\langle (R(n))^{2k+1} \rangle \,=\, 0 \,\,\,. 
\eeq
In this case we are fully entitled to refer to the law of iterated logarithms \cite{Feller} in order to conclude that 
\beq
\limsupn ~ \frac{| R(n)|}{\sqrt{2 n \log\log n}} \,= \, 1 
\hspace{3mm} 
, 
\hspace{3mm}
{\rm a. s. } 
\label{sqrtloglog22}
\eeq
where ``a.s.'' stays for ``almost surely'', in  the probabilistic sense \cite{Feller}. Hence, with probability 1, the abscissa of convergence  $\sigma_*$ of the random Dirichlet series 
(\ref{kahareseries}) is exactly $\sigma_* = 1/2$ (see also \cite{Kahane}). This is an encouraging result on the road to prove that also the restricted M\"{o}bius function $\hat\mu(s)$ has its abscissa of 
convergence equal to $\sigma_* = 1/2$. .

\subsection{Mertens Conjecture and its disproval} The random Dirichlet series just analyzed allows us to make an important comment on one of the first attempts to prove the 
RH that, however, was finally proven wrong. We are referring to the famous conjecture about the Mertens function $M(x)$. The story is well known 
(see, for instance, \cite{Edwards,Ng,Kotnik,Riele2,Riele3,Pinz}): it was conjectured by Thomas Joannes Stieltjes\footnote{This conjecture was written down in a letter by 
Stieltjes sent to Charles Hermite in 1885  (reprinted in Stieltjes (1905) and also printed later by Franz Mertens (1897)).}, that the Mertens function satisfies the bound 
\beq
| M(x) | < \sqrt{x} \,\,\,.
\label{mertensconj}
\eeq
In 1885 Stieltjes claimed to have proven a weaker result \cite{StieltjesC}, namely that the quantity 
\beq
t(n) \,=\, \frac{M(n)}{\sqrt{n}}\,\,\,,
\label{definitionofq(x)}
\eeq
previously introduced, is always bounded. However, he never published a proof of this statement. Of course if eq.\,(\ref{mertensconj}) was true, the RH would be true as well. 
In 1985, Andrew Odlyzko and Herman te Riele proved however that the strong version of the Mertens conjecture is false using the Lenstra-Lenstra-Lov\'{a}sz 
lattice basis reduction algorithm \cite{Odlyzyko2}, in the sense that they were able to show that 
\beq
\liminfn  ~ t(n) < -1.009 
\hspace{3mm}
,
\hspace{3mm}
\limsupn ~  t(n) > 1.06 \,\,\,. 
\eeq
Of course, one can  hardly be surprised of this result, since the order of growth of a function made of a random sequence of $\pm 1$ (as the function $R(n)$ above) 
is, with probability 1, of the order $\sqrt{2 n \log\log n}$. This simply comes from the law of iterated logarithm \cite{Feller}. This means that, sooner or later, the hypothetical bound (\ref{mertensconj}) of Stieltjes was  expected to be violated, as indeed it is.

\section{From RH toward (non-principal) GRH}\label{GRHRH}

Let us now discuss how the two cases discussed above (i.e. the GRH for $L$-functions of non-principal characters and the RH for $L$-function of principal characters) can be related.  In particular, here we want to show how it is possible to argue that the GRH is true as far as the RH is true. We are obviously referring to the non-principal cases of the $L$-functions since the the validity of the GRH for the principal case is simply related to the RH, as discussed in \ref{pppLLL}.

For the sake of the argument  presented in this section, let us assume we know that 
\begin{itemize}
\item the restricted Mertens function goes as $|\hat M(x) |\sim x^{1/2 + \epsilon}$ (see Part B)
\item the restricted M\"{o}bius coefficients $\hat \mu(n)$ behave as random independent variables (see Part C).
\end{itemize} 
(Recall that the restriction is to square-free integers.)
These two points are, of course, the main objects of analysis of this work. Assuming both to be true,  our purpose  here is to show how, in addition to ensuring the validity of the RH, they also imply the validity of the GRH.

Let us then consider the Generalized Mertens function $M_\chi(x)$ for a non-principal character $\chi$ relative to modulus $q$. For simplicity we consider the case where the modulus $q$ is a prime, therefore $\varphi(q) = (q-1)$. Since $\chi( k \,q) =0$ (for any integer $k$), we can decompose the sum on the index $m$ into the sum on the $(q-1)$ residue classes with residue $r$ ($r=1,2,\ldots, q-1$). Moreover, using the periodic property of  the characters (see Section \ref{Dirichelet}), we have 
\beq
M_\chi(x) \,=\,\sum_{1 \leq m \leq x} \chi(m) \mu(m) \,=\, 
\sum_{r=1}^{q-1} \sum_{k=0}^{k_r^*} \chi( k q +r) \mu(k q +r) \,=\, 
 \sum_{r=1}^{q-1} \chi(r) M_{r}(x) \,\,\,,
 \label{firstidentity1}
 \eeq
 where 
 \beq
 M_{r}(x) \,\equiv \, \sum_{k=0}^{k_r^*} \mu( k q + r) \,\,\, ,
 \label{firstidentity2}
 \eeq
 and the upper indices $k_r^*$ are given by 
 \beq
 k_r^* \,=\, \left[\frac{x-r}{q} \right] \sim \left[\frac{x}{q}\right] \,\,\,.
 \eeq
 Since the characters are pure phases, we have the inequalities 
 \begin{eqnarray}
 | M_\chi(x) | &\leq & |\chi(1) M_{1}| + |\chi(2) M_{2}| + \cdots + |\chi(q-1) M_{q-1}(x) | 
 \label{ineqq}\\
& \leq & | M_{1}| + | M_{2}| + \cdots + | M_{q-1}(x) | \,\,\,. \nonumber 
 \end{eqnarray}
As in the case of the Riemann zeta function,  we can restrict the M\"{o}bius coefficients along each arithmetic sequence $ k q +r$ only to the square-free numbers in order to restrict to the non-zero values and to filter out their obvious periodicities. How many square-free numbers are in the arithmetic progression $x_{k,r} =k q + r$ (with $q$ and $r$ comprime, i.e. $(q,r)=1$), and with $x_{k,r} \leq X$? Let's call this number $\hat n$: it turns out to be the same for any residue class $r$ and can be estimated probabilistically, generalising the argument which will be presented in Section \ref{ssqqff}, with the result\footnote{To estimate the number of square-free numbers which are divisible by the $q$, i.e. with residue $r=0$, one must keep in mind that 
 the probability that a square-free number is divisible by a prime $q$ is $1/(q+1)$, see Section \ref{squarefreeeeee}.}
 \cite{Nunes} 
\beq 
\hat n \,=\,\sum_{m \leq X, (q,r) = 1} \mu^2(m) \sim \frac{6}{\pi^2}\,\frac{X}{q} 
 \,\prod_{p | q} \left(1 - \frac{1}{p^2}\right)^{-1} \,\,\,.
\label{sqarithm}
\eeq
An explicit check of this probabilistic prediction is shown in Table \ref{tablessqqff2}. 

\begin{table}[t]
\begin{center}
\begin{tabular}{| c| c|c|c|}\hline
 r  &  
  \#\, square-free numbers & probabilistic estimate & relative error \\\hline
$0\, $  & $3.799.542$ & $3.799.544$ & $5\times 10^{-5}$\\ 
$1\,  $  & $4.432.807 $ & $4.432.801$ & $ 1 \times 10^{-4}$\\
$2 \,  $  & $ 4.432.777$ & $4.432.801 $ & $5 \times 10^{-4}$\\
$3\, $  &  $ 4.432.811 $ & $ 4.432.801$ & $2 \times 10^{-4}$\\
$4 \, $ &   $4.432.800$ & $4.432.801$ & $ 2\times 10^{-5}$  \\
$5 \,  $ &   $4.432.822$ & $4.432.801$ & $5 \times 10^{-4}$\\
$6 \, $  & $4.432.784$ & $4.432.801$& $ 4 \times 10^{-4}$\\ 
\hline
\end{tabular}
\end{center}
\caption{
  Number of square-free numbers along the arithmetic sequences $x_{k,r} = k \, q + r$, where here the modulus $q$ is chosen to be $q=7$, 
  for $x_{k,r} \leq X$, with $X = 50.000.000$, compared with the probabilistic 
  prediction (\ref{sqarithm}). The final column reports the percent error of the theoretical estimate, namely $(N_{pr} - N_{count})/N_{pr}$, where 
  $N_{pr}$ is the expected number coming from probability argument while $N_{count}$ is the actual number coming from counting.
}
\label{tablessqqff2}
\end{table}

Let us now denote by $\hat M_\chi(x)$ and $\hat M_{r}(x)$ the corresponding quantities of the generalised Mertens functions restricted however to the square-free numbers. Making now the assumption that the restricted M\"{o}bius coefficients are largely independent of each other, as thoroughly shown in Part C of this paper, the various $M_{r}(x)$ in eq.\,(\ref{ineqq}) behave statistically in the same way and therefore we have 
\beq
{\rm max} \, |\hat M_\chi(x) | \leq (q-1)\,{\rm max}\, |\hat M_{r}|\,\,\,.
\label{Mqqr}
\eeq
Relying once again on random independence of the restricted M\"{o}bius coefficients, the ${\rm max} \,|\hat M_{r}(x)|$ is obviously related to the maximum of the familiar restricted Mertens function opportunely rescaled: indeed, the only thing which matters is the number of terms present in their sum, which can be estimated using (\ref{sqarithm}). Hence, assuming that ${\rm max}\, |\hat M(x)| \sim x^{1/2 +\epsilon}$, we have then 
\beq
{\rm max}\, |\hat M_{r} (x) | \sim \sqrt{\frac{1}{q}
  \prod_{p | q} \left(1 - \frac{1}{p^2}\right)^{-1}} \, x^{1/2 + \epsilon} \,\,\,,
\label{indepr}
\eeq
and therefore\footnote{It is important to stress that, assuming that the restricted M\"{o}bius coefficients $\hat \mu(n)$ are independent random variables, as all checks discussed in Part C unequivocaly show, 
the right hand side of eq.(\ref{indepr}) is independent on the residue index $r$.}
\beq
{\rm max}\, |\hat M_\chi(x) | \sim  A_q \, x^{1/2 +\epsilon}\,\,,
\label{basicinequality}
\eeq
where 
\beq
A_q \,=\, \frac{(q-1)}{\sqrt{q}} \, \sqrt{ \prod_{p | q} \left(1 - \frac{1}{p^2}\right)^{-1}}\,\,\,. 
\eeq
The inequality (\ref{basicinequality}) shows that the GRH will be true as long as the RH holds true, the only difference being the presence of the overall constant $A_q$. 
In the following our efforts will  then be focused on arguing  that the restricted Mertens function 
behaves as $\hat M(x) \sim x^{1/2 + \epsilon}$ and to show that the restricted M\"{o}bius coefficients behave as random independent variables. We will come back to the relation (\ref{Mqqr}) between GRH and RH in our conclusions. 
 
\section{Summary}
 
In this Part A we have collected some well-known but also less-known properties of the Dirichlet and Riemann functions. We have clarified why the Riemann function may be considered as a particular case of the general Dirichlet $L$-functions, although we have shown that it is a case that needs a different approach than the one used previously for discussing the GRH for a generic Dirichlet $L$-function of non-principal character  \cite{ML,LM}: indeed, while for a generic Dirichlet $L$-function of non-principal character of modulus $q$ one can employ their infinite product representation and work directly on the statistical distribution of the residues (mod q) of the prime numbers\footnote{The obvious advantage of working, in this context, with prime numbers rather than natural numbers is that for prime numbers there are theorems, such as the Dirichlet theorem for the equidistribution of residues (mod q), or robust conjectures, such as the Hardy-Littlewood result for the correlations of these residues, which are extremely helpful in the statistical approach to the Generalized Riemann Hypothesis. On the other hand, for the M\"{o}bius coefficients on square-free numbers, their correlations and the corresponding behaviour of the Mertens function are much less known.}, for the Riemann zeta function it is instead necessary to employ the natural numbers and the statistical distribution of the Mertens function and its relative M\"{o}bius coefficients coming from the multiplicative inverse function of the Riemann's, i.e. $1/\zeta(s)$. However, in this part we have also argued that, establishing the validity of the RH through the asymptotic behavior of the Mertens function, could be enough to show the validity of the GRH as well. 

We have also recalled the duality properties of the Riemann/Dirichlet $L$ functions and a theorem of Grosswald and Schnitzer about random functions defined by infinite product representation on pseudo-primes which share exactly the same zeros in the critical strip of the Riemann/Dirichlet functions: this theorem offers a particular perspective on the (Generalized) Riemann Hypothesis and the role of randomness in pursuing its proof. Finally, we have identified as key object of our study the asymptotic behaviour of the Mertens function $\hat M(n)$ restricted to the square-free numbers: this will be our main  focus of the next Part B and Part C of this paper.   

\newpage

\addcontentsline{toc}{section}{Part B}
\begin{center}
{\Large {\bf PART B}}
\end{center}

In this Part B of the paper we present an alternative view of Mertens function which allows us to clarify many of its properties, in  particular the fluctuations of this function. A key tool of our analysis will be the set of square-free numbers, since they are the only ones for which $\mu(n) \neq 0$. In the next sections we will establish a prime number theorem restricted to the set of square-free numbers and, in addition, we will also address the counting of the number of divisors for the square-free numbers and their Poisson distribution. Using these results, we will formulate the Erd\H{o}s-Kac theorem relative to the square-free numbers. More importantly, we will be able to study the mean and variance of the restricted Mertens function arriving in this way to the important result that $\langle (\hat M(n))^{2k}\rangle \sim n^k $, quite relevant for the validity of the RH. 

Let's start our analysis discussing in more detail the square-free numbers and their properties.

\section{Square-free numbers and their M\"{o}bius coefficients}\label{ssqqff}

 The square-free numbers $\sqf_1,\sqf_2, \ldots,\sqf_n$ are those integers which are divisible by no perfect square other than 1. That is, their prime factorization has exactly one factor for each prime that appears in it. A generic square-free number is given by 
\begin{equation}
\sqf \,=\,p_1^{\alpha_1}\,p_2^{\alpha_2} \cdots \,p_k^{\alpha_k} 
\hspace{3mm}, 
\hspace{3mm} 
\,\,\,\alpha_a \in \{0,1\} 
\label{ssff}
\end{equation}
and their first representatives are $\{ \sqf_n \} \,= \{\,2, 3, 5, 6, 7, 10,  11, 13, 14, 15, 17, 19, 21,..\}$. Remarkably, these numbers are a finite fraction of all the integers, as we have seen in the Introduction.
Hence, denoting by Q(x) the number of square-free integers between $1$ and 
$x$, we have
\beq
Q(x)\, \sim \, \frac{6}{\pi^2} \, x \,\,\,.
\label{numbersquarefree}
\eeq
 In view of (\ref{numbersquarefree}), notice that the $n$-th square-free number has an approximate value 
\beq
\sqf_n \sim \frac{\pi^2}{6} \,n \,= 1.64493.. \, n 
\label{nthsquarefree}
\eeq
Square-free numbers have been the subject of several papers, also in relation to the M\"{o}bius function \cite{Cellarosi1,Cellarosi2,Granvillesqfree}.  A physical realization of square-free numbers has been recently proposed in terms of hard-core bosons and a ladder system made of coupled quantum spin chains \cite{MTZ}. As we will see below, square-free numbers acquire a suggestive interpretation by promoting the prime numbers to quantum energies of a fermionic system \cite{Julia, Spector}. 

\begin{figure}[t]
\centering\includegraphics[width=.7\textwidth]{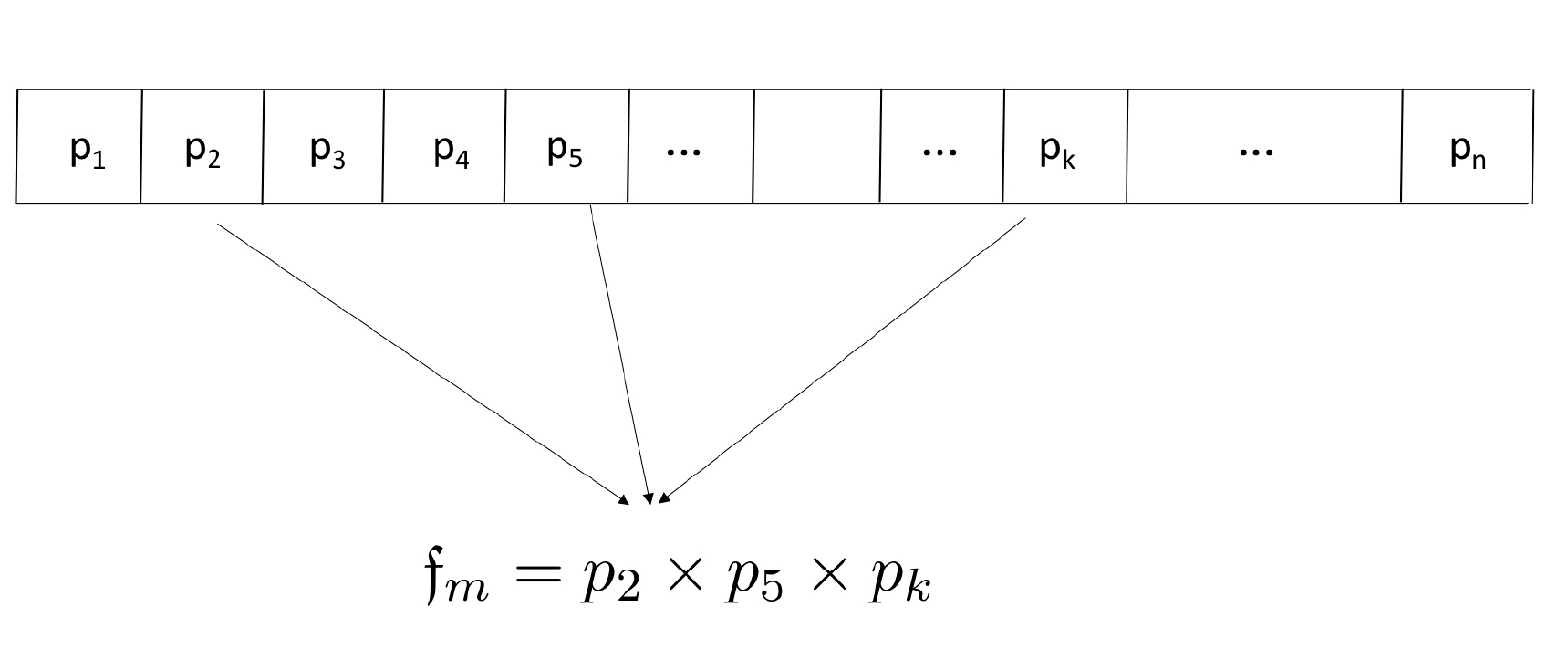}
\caption{Boxes associated to the prime numbers and draw of a square-free number.}
\label{Primeboxes}
\end{figure}

\vspace{3mm}
\noindent
{\bf Randomness and correlations of M\"{o}bius coefficients}. By using an argument which can be called the {\em random box lottery draw}, it is rather simple to understand the random nature of the M\"{o}bius coefficients as well as the origin of their eventual correlations. The argument consists of the following. Imagine we put each prime, in an increasing order, in a box along a row, so that the entire sequence of primes is associated to the infinite sequence of boxes shown in Figure \ref{Primeboxes}. Let's associate to each prime a fermion degree of freedom, so that we can implement a draw protocol as follows: we are allowed to select randomly any number of boxes, but never choosing the same prime twice  (respecting the Pauli principle of the fermions). In this way, after each of these draws of the boxes, we end up with a square free number. Each of these selections has a definite parity, corresponding to the even or odd number of boxes chosen. Let's denote the generic element of these two even and odd classes as $T_{\pm}$. Given the infinite sequence of primes, it is obvious that to any given even selection $T_+$ we can immediately find another odd selection $T_-$ (either extracting one or an odd number of boxes and adding to the previous selection or eliminate randomly one or an odd number of the selected boxes), and vice versa. Any of these selections obviously corresponds to a square-free number of a given parity, to which we can associate the corresponding value of the M\"{o}bius coefficient according to the law (\ref{mudef}). Since the boxes are all equivalent and the choices of the boxes absolutely random, it is obvious that the values assumed by the M\"{o}bius coefficients have to be random as well. Moreover, from the perfect equivalence between even and odd number of boxes chosen, the mean value of the (restricted) Mertens function has to be zero (we  previously discussed that this is a rigorous result).

So, according to the {\em random box lottery draw} argument we just presented, it would seem natural to argue positively about the randomness of the restricted M\"{o}bius coefficients and the random walk behavior of their sum, given by the restricted Mertens function. However, we must be aware that a correlation between these values can arise if we start {\em ordering} the square free numbers which we extracted! Actually, this ordering is inherently present in the definition of the Mertens function, where we sum on the values of the M\"{o}bius coefficients in increasing order of the square-free numbers involved. Consider, for instance, that using our lottery draw we have extracted this random sequence of square-free numbers 
\beq
\sqf_{100}^{(1)} \longrightarrow \sqf_{3}^{(2)} \longrightarrow \sqf_{72}^{(3)}\longrightarrow  \sqf_{15}^{(4)} \cdots  \longrightarrow \sqf_{k}^{(n)}  \longrightarrow 
\label{trialsequence}
\eeq
where the upper index refers to the $n$-th draw while the lower index refers instead to the actual ordering index of the square-free numbers. As we said earlier, we expect that, on average, such a sequence contains square-free numbers with as many even and odd number of primes, i.e. on average an equal number of $\pm 1$ values for the M\"{o}bius function 
computed along the sequence (\ref{trialsequence}). But, when we order it according to the increasing order of the $a_k$ index, we end up in a {\em different} sequence 
\beq 
\sqf_3  \longrightarrow \sqf_{15}  \longrightarrow  \sqf_{72}  \longrightarrow  \sqf_{100}  \longrightarrow , \ldots 
\eeq
where now the sequence of $\pm1$ of the M\'{o}bius function computed for these numbers may show, in principle, a certain unbalance in the order in which they appear. 
Said differently, the only feature responsible for an eventual correlation of the restricted M\"{o}bius coefficients is the ordering of the drawing: {\em if} this would end up in a strong correlation, this could lead to a violation of their random walk behavior and {\em then} to a violation of the RH. In the following we are going to show that this seems to be {\em not} the case, namely if correlations in the sequence of M\"{o}bius numbers  exist, they are so weak that they do not spoil the central limit theorem and therefore the validity of the RH.

\section{Prime number theorem for square-free numbers}

In the last section we have seen that the square-free numbers are a finite fraction of the integers. It is then natural to ask how many primes there are in their sequence. Before facing this problem, let's however initially present the probabilistic derivation of the prime number theorem in the usual sequence of integers, since this leads us to the answer to the same question for the sequence of square free numbers.

\subsection{Prime number theorem for integers} 
The probabilistic argument for the prime distribution among the integers goes as follows (see, for instance \cite{Schroeder}): let $1/p_k$ be the probability that a generic integer is divisible by the prime $p_k$ since, after all, among $p_k$ consecutive integers there is one which is divisible by $p_k$. Assuming that divisibility for different primes is independent, the probability $U(x)$ that an integer $x$ is not divisible by any prime below it is given by 
\beq
U(x) \sim \, \prod_{p_k < x} \left(1 - \frac{1}{p_k} \right)\,\,\,.
\label{basicprob}
\eeq
Of course, if $x$ is not divisible by any smaller prime below it, $x$ is a prime as well and therefore $U(x)$ is the probability that $x$ is a prime. 
Taking the logarithm of this expression and expanding at the lowest order in $1/p$, we have 
\beq
\log\,U(x) \,\sim \, \sum_{p_k < x} \log\left(1 - \frac{1}{p_k}\right) \,\sim - \sum_{p_k < x} ~ \frac{1}{p_k} \,\,\,.
\label{leading}
\eeq
We can now convert the sum on the primes into a sum over  all integers by making use of the function $U(x)$ itself, namely
\beq
\log\,U(x) \,\sim \, - \sum_{m=1}^x \frac{U(m)}{m} \,\sim -\int_2^x \frac{U(m)}{m} \, dm \,\,\,.
\label{integral}
\eeq
Taking now the derivative wrt $x$ on both sides of this equation, we get the differential equation satisfied by $U(x)$ 
\beq
\frac{U'(x)}{U(x)} \,=\, \frac{U(x)}{x} \,\,\,,
\eeq
whose solution is 
\beq
U(x) \,=\, \frac{1}{\log x} \,\,\,.
\label{probabilitynatural}
\eeq
Therefore, calling $\pi(x)$ the number of primes smaller than $x$, we arrive at  the result that this function is asymptotically given by 
\beq
\pi(x) \,\sim \int_2^x \frac{dx}{\log x} \,\equiv {\rm Li}(x) \,\,\,,
\label{primenumbertheorem}
\eeq
where ${\rm Li}(x)$ is the ``logarithmic  integral function".  As well known, eq.\,(\ref{primenumbertheorem}) expresses the content of the 
Prime Number Theorem.

\begin{table}[t]
\begin{center}
\begin{tabular}{|l|l|l|}\hline
  primes  &  
  numerical & theoretical \\
$p_a$ & probability & probability  \\\hline
$2 $ & $0.333331$ & $0.333333$  \\
$3 $ & $0.249998$ & $0.250000$\\
$5$ & $0.166670 $ & $0.166666$\\
$7 $ & $0.1250001 $ & 0.125000\\
$11$ & $0.0833331$ & $0.0833333$ \\
$13 $ & $0.0714281$ & $ 0.0714286$\\
$17$ & $0.0555547$ &$0.0555556$ \\ 
\hline
\end{tabular}
\end{center}
\caption{
  Numerical vs theoretical probability of divisibility of square-free numbers by a prime. The 
  numerical data are given by the ratio $N_a/N$, where $N=10^7$ is the number of square-free numbers considered and $N_a$ is the number of them divisible by the prime $p_a$. The theoretical probability is $1/(p_a +1)$.}
 \label{tablep1} 
\end{table}

\subsection{Probability of divisibility by $p_k$ of a square-free number}\label{squarefreeeeee} 
We would like now to estimate the function  $\hat\pi(n)$ 
defined as 
\beq
\hat\pi(n) \,=\, \# \left\{ \sqf_a \leq \sqf_n \,\,: \sqf_a = {\rm prime}\right\}\,\,\,.
 \eeq
This function counts the number of primes which appears in the sequence of square-free numbers less than $\sqf_n$ (notice that $n$ refers to the index of the square-free number). Of course 
$\hat\pi(n) = \pi(\sqf_n)$ but let's see how we can determine $\hat\pi(n)$ self-consistently.  To proceed, we have initially to determine the probability $\tilde \Delta_k$ that a randomly chosen square-free number is divisible by a prime factor $p_k$. As shown in \cite{MTZ}, instead of being $1/p_i$, this probability is given instead by 
\beq
\tilde \Delta_k \,=\, \frac{1}{p_k +1} \,\,\,. 
\label{probsqf}
\eeq
Although asymptotically $p_k +1 \approx p_k$,  the shift by $1$ has some important consequences.    For instance,   the fraction of square free numbers $\sqf_n$ that are even
is $1/3$ rather than $1/2$.   
It is relatively simple to show that this is indeed the correct result by making a simple numerical check, as shown in Table \ref{tablep1}. For proving (\ref{probsqf}), 
one can use the inclusion-exclusion principle as follows. Let $Q(x)$ the number of square-free numbers less than $x$, whose behavior is given in eq.\,(\ref{nthsquarefree}). Using $Q(x)$, we can give the first estimate of the number of square-free numbers which are less than $x$ and multiples of the prime $p_k$. This number is approximatively equal to $Q\left(x/p_k\right)$. If we now multiply a square-free number $y$ (with $y \leq x/p_k$) by $p_k$, this yields (for sure) a multiple of $p_k$ which is $\leq x$. This multiplication usually gives rise to a number which is also square-free, for the only perfect square that could possibly divide the number $y p_k$, where $y \leq x/p_k$ and $y$ is square-free, is $p_k^2$. This implies that $Q(x/p_k)$ over counts the set of multiples of $p_k$ that are $\leq x$ and square-free. In order to correct this discrepancy, we must subtract approximately $Q(x/p_k^2)$, which almost counts how many numbers $\leq x$ are divisible by $p_k^2$ but which are otherwise square-free. But this time we have subtracted too much, since we have also subtracted the numbers $\leq x$ which are divisible by $p_k^3$ but which are otherwise square-free. So, we need to add back approximately $Q(x/p_k^3)$ and so on. In this way, we have to deal with the sum of the infinite series
\EQ
\tilde \Delta_k\,=\, \frac{1}{p_k} - \frac{1}{p_k^2} + \frac{1}{p_k^3} - \frac{1}{p_k^4} + \cdots \,=\,\frac{1}{p_k+1}
\EN 
This yields the sought for  probability for a randomly  chosen square-free number to be divisible by a prime $p_k$. Hence, with respect to the integers, in the case of square-free numbers there is a renormalization in the expression of the probability  of divisibility by a prime. 

\begin{table}[b]
\begin{center}
\begin{tabular}{|l|l|l|l|}\hline
 \, n  &  
  \#\, primes & $\hat\pi(n)$ & relative error \\\hline
$1\, \times10^6 $  & $124.281$ & $124. 419$ & $1.1\times 10^{-3}$\\ 
$2\, \times 10^6 $  & $236. 242 $ & $236. 344$ & $ 4.3 \times 10^{-4}$\\
$3 \, \times 10^6 $  & $  344.244 $ & $344. 408 $ & $4.7 \times 10^{-4}$\\
$4\, \times 10^6 $  &  $449. 850$ & $450.109$ & $ 5.7 \times 10^{-4}$\\
$5 \, \times 10^6$ &   $553. 878$ & $554. 120$ & $ 4.3 \times 10^{-4}$  \\
$6 \, \times 10^6 $ &   $ 656. 561$ & $656. 823$ & $3.9 \times 10^{-4}$\\
$7 \, \times 10^6 $  & $758. 165$ & $758. 463$& $ 3.9 \times 10^{-4}$\\ 
$8\, \times 10^6 $  & $848. 921$ & $849. 173$ & $ 2.9 \times 10^{-4}$\\ 
$9\, \times 10^6 $ &  $938. 967$ & $939. 251$ & $ 3.0 \times 10^{-4} $\\
$1\,\times 10^7$ &  $ 1.028.462$ & $1.028.770$ & $ 2.9 \times 10^{-4}$\\
\hline
\end{tabular}
\end{center}
\caption{
  Number of primes along the sequence of square-free numbers compared with the approximate theoretical expression $\hat\pi(n)$ given in eq.\,(\ref{primenumbertheoremsqf}). 
  The final column reports the percent error of the theoretical estimate.   }
\label{tablep2}
\end{table}

\subsection{Prime number theorem for square-free numbers} 
\begin{figure}[t]
\centering\includegraphics[width=.5\textwidth]{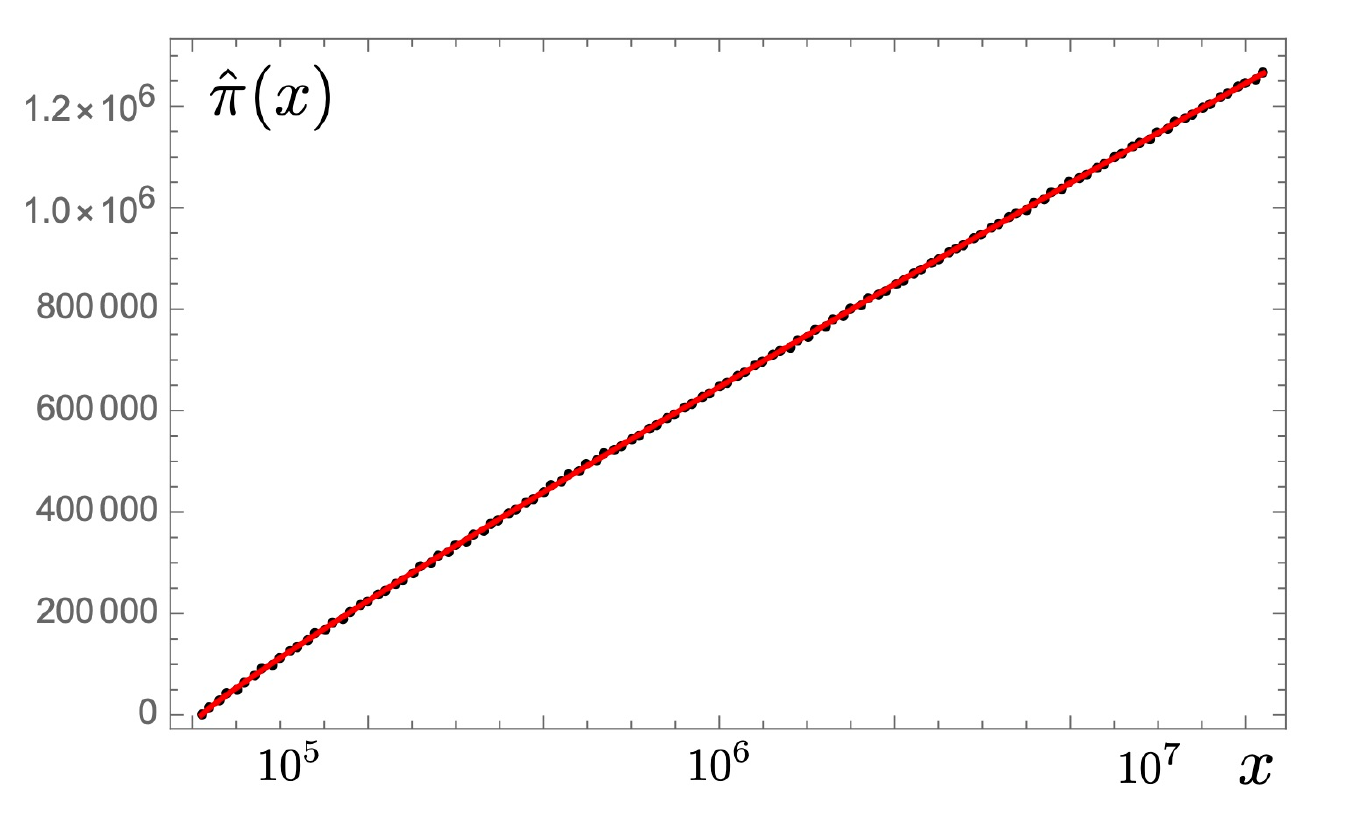}
\caption{Number of primes along the sequence of the square-free numbers (black dots) versus $\hat\pi(x)$, in the interval $x \in (2,3.2 \times 10^7$). The two curves are 
essentially indistinguishable.}
\label{comparisonnupr}
\end{figure}

We can now repeat the same steps which led us to the approximate expression $\pi(x)$ of eq.\,(\ref{primenumbertheorem}) employing though in the initial formula (\ref{basicprob}) the probability $\tilde \Delta_k$ determined above. Taking also into account the behaviour (\ref{nthsquarefree}) of the $n$-th square-free number we have
\beq
U_{sq}(x) \,\sim \frac{1}{\log(x+1)} \,\,\,,
\label{squarefreeW}
\eeq 
and therefore we arrive to 
\beq
\hat\pi(n) \,\sim \,\int_2^{\sqf_n} \frac{dx}{\log(x+1)}\,\equiv \, {\rm Li}_{sq}(\sqf_n) \,\,\,.
\label{primenumbertheoremsqf}
\eeq
We can check how good this approximation of the function $\hat{\pi}(x)$ is just by counting the primes in the sequence of the square-free numbers, as done in Table \ref{tablep2} and shown in Figure \ref{comparisonnupr}, with the relative error around $10^{-4}$, which is quite satisfactory.  Of course, for very large values of  $\sqf_n$, one recovers the expected identity $\hat\pi(n) = \pi(\sqf_n)$.

\section{A different look at Mertens function: I}\label{globalapproach1}
In this section and the next one, we are going to have a different look at the Mertens function, adopting a ``global'' point of view of this function rather than a ``local'' one. 
To this aim, let's first  define and discuss the {\em primorial}.

\subsection{Primorial}

The primorial $\mathbb{P}(q)$ is similar to the factorial function, but rather than successively multiplying positive integers, the function employs the multiplication 
of the first prime numbers up to the $q$-th prime 
\beq
\mathbb{P}(q) \,\equiv \prod_{k=1}^q p_k \,\,\,.
\label{primorial}
\eeq
As shown in Table \ref{TablePrimorial}, such a function grows exponentially fast and its scaling behaviour is captured by the following law 
\beq
\mathbb{P}(q) \sim e^{p_q} \sim e^{q \log q} \,\,\,, 
\eeq
which employs the $q$-th prime number \cite{Ruiz}. Any primorial $\mathbb{P}(q)$ is clearly a square-free number. 
\begin{table}[t]
\begin{center}
\begin{tabular}{| l | l |}\hline 
n & $\mathbb{P}(n)$ \\\hline
$ 1 $ & $2$  \\
 $2 $ & $6 $ \\
$ 3 $ & $ 30 $ \\
$ 4 $ & $210 $ \\
$ 5 $ & $2310 $ \\
$ 6 $ & $ 30030$ \\
$ 7 $ & $510510 $ \\
$ 8 $ & $9.69969\times 10^6$ \\
$ 9 $ & $2.23093\times 10^8 $ \\
$ 10 $ & $6.46969\times 10^9 $ \\
$ 11 $ & $ 2.00560\times 10^{11} $\\
 $12$ & $7.42074\times 10^{12} $\\
$ 13 $ & $3.04250\times 10^{14} $\\
 $14 $ & $1.30828\times 10^{16} $\\
 $15$ & $6.14890\times 10^{17} $\\
\hline
\end{tabular}
\end{center}
\caption{The first $15$-th values of the primorial $\mathbb{P}(n)$.  }
\label{TablePrimorial}
\end{table}

\begin{figure}[b]
\centering\includegraphics[width=.6\textwidth]{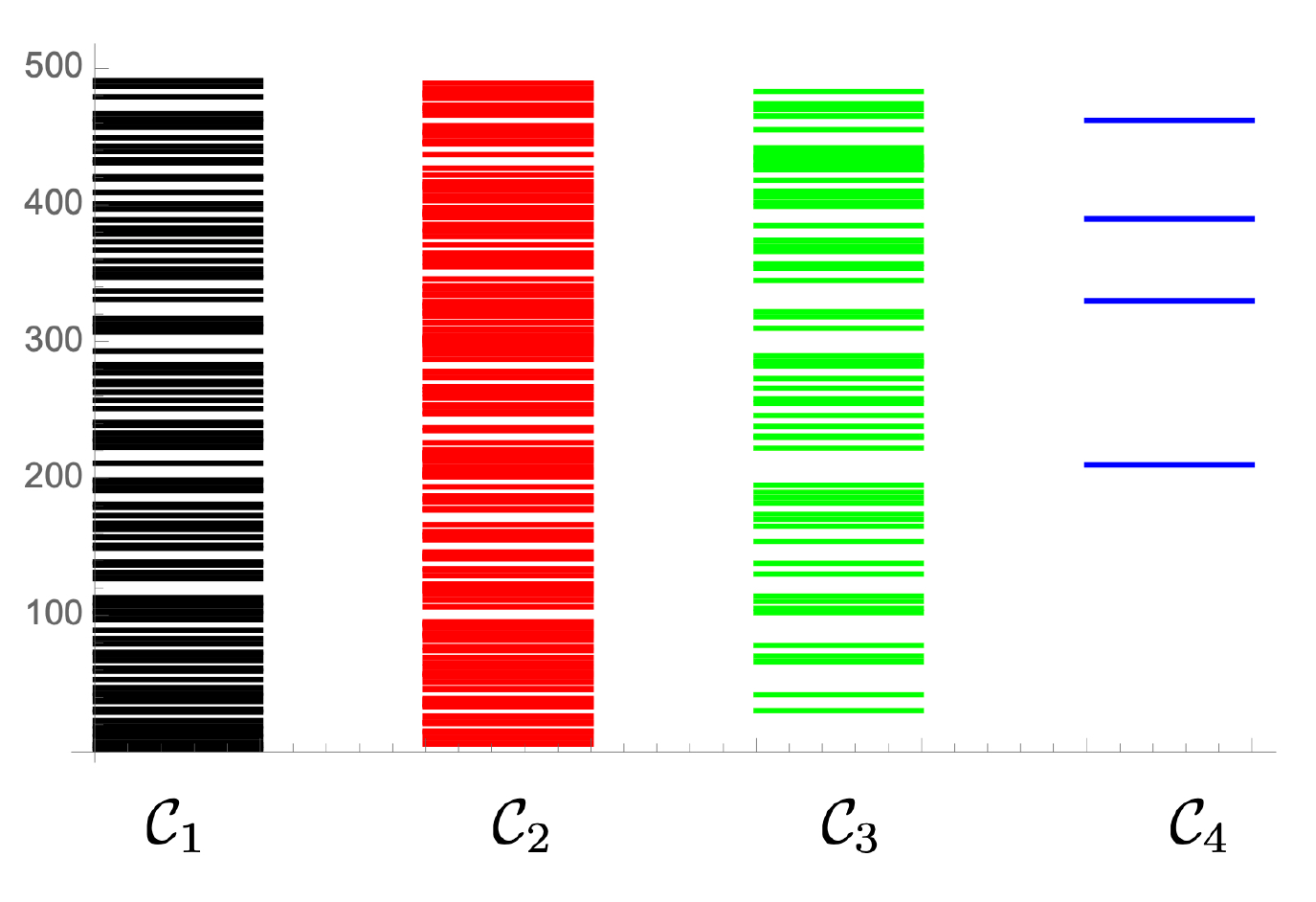}
\caption{The first $n=300$ square-free integers organized according to the number of their divisors, i.e. according to the classes ${\mathcal C}_k$.
 The figure shows the 
interpenetration of the various classes.}
\label{spectroscopy}
\end{figure}

\subsection{ Sets of square-free numbers made of $k$ primes}

Let's partition the square-free numbers into the following sets: 
\beq
\begin{array}{l}
{\mathcal C}_1 \,=\, \{\sqf_a \,: \, \sqf_a = p_{i_1}\} \\
{\mathcal C}_2 \,=\, \{\sqf_a \,: \, \sqf_a = p_{i_1} \cdot p_{i_2}\,\,\,\,\,, i_1\neq i_2\}\\
{\mathcal C}_3 \,=\, \{\sqf_a \,: \, \sqf_a = p_{i_1}\cdot p_{i_2} \cdot p_{i_3} \,\,\, \,\,, i_1 \neq i_2 \neq i_3\}\\
\ldots \hspace{15mm}\ldots \\
\ldots \hspace{15mm}\ldots\\
{\mathcal C}_k\,=\,\{\sqf_a \,: \, \sqf_a = p_{i_1}\cdot p_{i_2} \cdots p_{i_k} \,\,\,\,\,,  i_1 \neq i_2 \neq i_3 \cdots \neq i_k\} 
\end{array}
\eeq
Each of these sets ${\mathcal C}_k$ has a minimum, given by the corresponding primorial
\beq
{\rm min}\, {\mathcal C}_k =  \mathbb{P}(k)\,\,\,.
\label{minimum}
\eeq
Clearly the sequence made of the minima of ${\mathcal C}_k$ is monotonic
\beq 
{\rm min}\, {\mathcal C}_r \, < \, {\rm min}\, {\mathcal C}_{s} \,\,\,\,\,\,, \,\,\,\,\,\,\,r < s\,\,\,, 
\eeq
although the square-free numbers present in each class can be arbitrarly close to each other (for integers, their minimum difference is of course 1).  Thinking of the prime numbers as 
the elementary particles of arithmetic, we can interpret  the square-free numbers $\sqf_a$ of the class ${\mathcal C}_k$ as made of $k$ of these elementary particles, with the energy  given by the square-free numbers themselves.  We can then plot the various energies (alias the sequence of the square-free numbers), taking care of their "spectroscopy", as done in Figure \ref{spectroscopy} which makes pretty evident the relative {\em degeneracy} of the various classes.

\begin{figure}[b]
\centering\includegraphics[width=.7\textwidth]{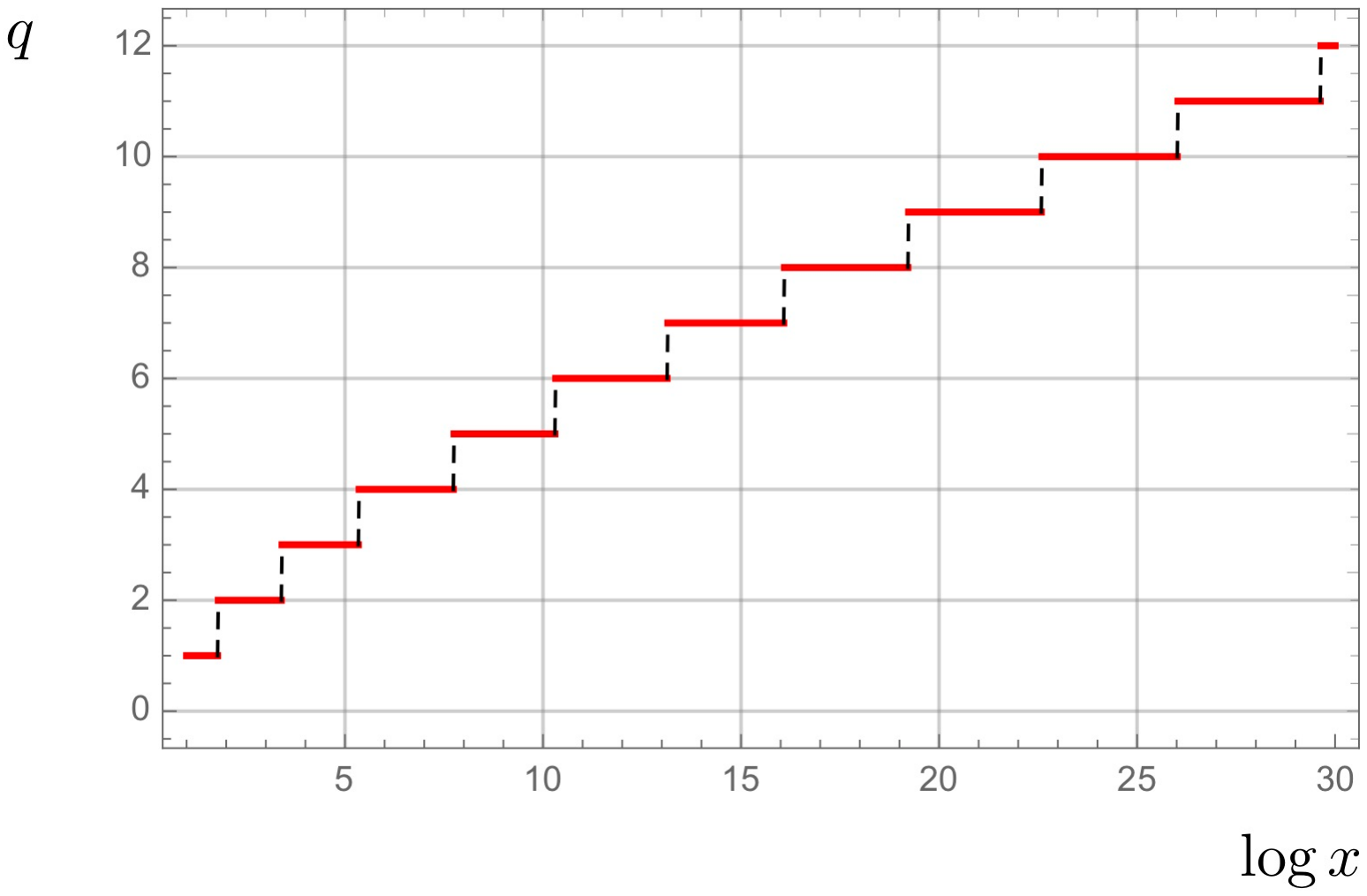}
\caption{Number $q$ of factors versus $\log x$.}  
\label{Primorialdivisor}
\end{figure}

\bigskip

\subsection{ Rewriting of the restricted Mertens function}

We can express differently the restricted Mertens functions adopting a ``global'' point of view of this function rather than a ``local'' one. By this we mean the following: so far (see eq.\,(\ref{hatmertens})), we have considered the Mertens function $\hat M(n)$ as the sum of the sequential series of the restricted M\"{o}bius coefficients $\hat\mu(n)$. 
But we can gain interesting information on the restricted Mertens function by fully exploiting the definition of the M\'{o}bius coefficients, and the definition of the sets ${\mathcal C}_k$ given above. Indeed, we can organize the expression of $\hat M(n)$ as follows\footnote{The sign $(-1)^k$ is equal to the so-called Liouville function $\lambda(m)$, given by $\lambda(m) = (-1)^{\Omega(m)}$. 
The Dirichlet series for the Liouville function is related to the Riemann zeta function by $$\frac{\zeta(2 s)}{\zeta(s)} \,=\,\sum_{m=1}^\infty \frac{\lambda(m)}{m^s}$$.}
\beq
\hat M(n) \,=\, \sum_{k=1}^{q(n)} (-1)^k \, N_k(n) \,\,\,,
\label{rewritingMertens}
\eeq
where $N_k(n)$ are the number of square-free numbers less than $\sqf_n$ made of just $k$ factors and therefore which belong to the set ${\mathcal C}_k$
\beq
N_k(n) \,\equiv \# \left\{\sqf_a \leq \sqf_n \,\, : \,\, \sqf_a = \underbrace{p_{a_1} \cdot p_{a_2} \cdots p_{a_k}}_\text{k\,prime\,factors}\right\}
\,=\, \# \left\{\sqf_a \leq \sqf_n \,\, : \,\Omega(\sqf_a) = k\right\} \,\,\,.
\eeq
We have introduced the arithmetic function $\Omega(b)$ which counts the total number of prime factors of the natural number $b$ and which, in general, also includes their multiplicities\footnote{
In number theory, the function $\omega(b)$  counts the number of distinct prime factor present in the natural number $b$ while 
$\Omega (b)$ count the number of prime factors of a natural number $b$, honoring their multiplicity. Of course, for a square-free number 
$\sqf_n$ we have $\omega(\sqf_n) = \Omega(\sqf_n)$.}. The $N_k(n)$'s are fluctuating quantities whose properties we will comment on in detail in the next section. 

The upper limit $q(n)$ of the sum in eq.\,(\ref{rewritingMertens}) is given by the integer $\eta_n$ which is the index of the primorials $\mathbb{P}(\eta_n)$ and $\mathbb{P}(\eta_n+1)$ satisfying the inequality 
\beq 
\mathbb{P}(\eta_n) \, \leq \, \sqf_n \, < \, \mathbb{P}(\eta_n+1) \,\,\,, 
\label{inequality}
\eeq 
namely 
\beq
q(n) \,=\, \eta_n\,\,\,. 
\label{eqsfond}
\eeq
Therefore, there is an important conclusion of this analysis:  for exponentially large intervals of the variable $n$, the number of terms $q(n)$ in $\hat M(n)$ is constant and relatively small!  So, for instance, up to $\sqf_n \sim 10^{10}$, the Mertens function $\hat M(n)$ involves square-free numbers made up of at most of $9$ different primes (see, for instance, 
Figure \ref{Primorialdivisor}). 

\section{A different look at Mertens function: II}\label{globalapproach2}
\begin{figure}[b]
\centering\includegraphics[width=.7\textwidth]{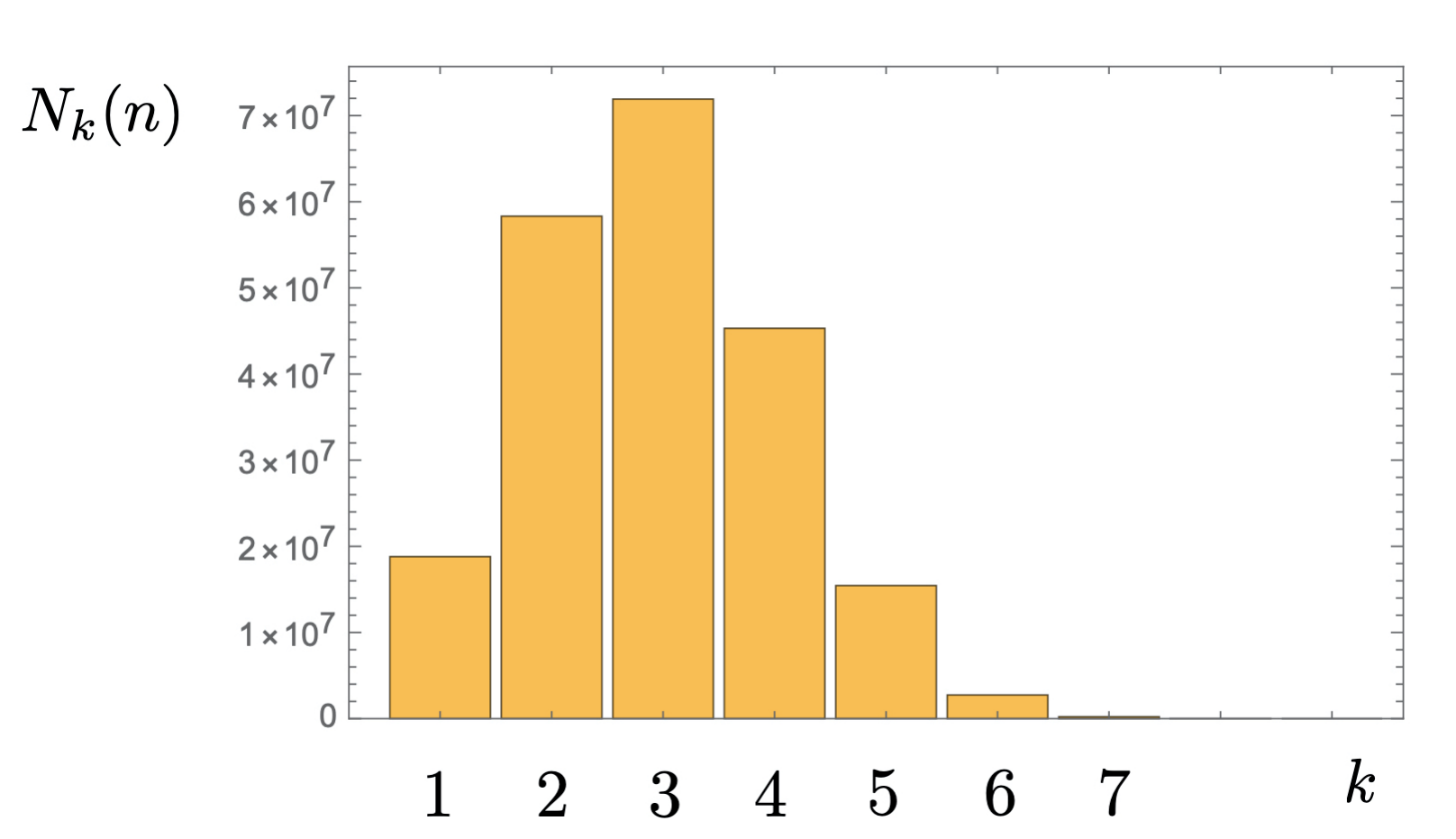}
\caption{The numbers $N_k(n)$ of square-free numbers less than $\sqf_n$ made of just $k$ factors 
versus $k$ for all $\sqf_a\leq \sqf_n = 182.378.126$, namely for $n=110.872.566$.
}

\label{Poisson1}
\end{figure}
Let's now analyze in more detail the quantities $N_k(n)$. From the previous section we know that each $N_k(n)$ becomes different from zero above a certain threshold, given by 
\beq
N_k(n) > 0 \,\,\,\,\,,\,\,\,\,\, {\rm if} \, \,\,\sqf_n \geq \mathbb{P}(k) \,\,\,.
\ee
This means that, in any exponentially large intervals $I_q = (\mathbb{P}(q),\mathbb{P}(q+1))$, the values taken by the Mertens function $\hat M(n)$ is determined only by the fluctuations of the $q(n)$ quantities $N_k(n)$ ($k=1,2,\ldots,q(n))$. A typical distribution of these quantities is presented in Figure \ref{Poisson1} and, as shown by this figure and as we are going to heuristically argue below, these quantities are distributed according to a Poisson distribution. To arrive to this result, we need to study initially the behavior of the number of prime divisors for square-free numbers. It is worth mentioning that, for the usual integers,  Wintner and Granville \cite{Wintner,GranvillePoisson} argued about the validity of the Poisson distribution for the relative frequency (``probability") of those natural numbers which possess $m$ divisors: Wintner's argument was made more quantitative by Granville \cite{GranvillePoisson} 
and also by Schr\"{o}der \cite{Schroeder} and, in the following, we take advantage of the method proposed by Schr\"{o}der to arrive at a quantitative estimate, {more specifically the distribution}   of the number of prime divisors for square-free numbers. 
\subsection{Average of the number of prime divisors}

 Let's define the function $\hat\omega(n)$ for the square-free numbers  according to 
 \beq
\hat{\mu} (n) = \mu(\sqf_n) = (-1)^{\hat\omega (n)}\,\,\,. 
\eeq
Equivalently 
\beq
\hat\omega(n) \,=\,\sum_{p_i | \sqf_n} 1 \,\,\,.
\label{omega}
\eeq
This function counts the number of (necessarily distinct) prime number divisors of the square-free numbers $\sqf_n$.
As in previous functions, we have $\hat\omega(n) = \omega(\sqf_n)$. 
 It is an additive function, because if two square-free numbers $\sqf_n$ and $\sqf_m$ are coprime, 
$(\sqf_n,\sqf_m) =1$, we have 
\beq
\omega(\sqf_n \, \sqf_m) \,=\, \sum_{p_i | \sqf_n \sqf_m} 1 \,=\, \sum_{p_i | \sqf_n} 1 + \sum_{p_i | \sqf_m} 1 \,=\, \omega(\sqf_n) + \omega(\sqf_m) \,\,\,, 
\eeq
so that, defining $\hat\omega(n \star m) \equiv \omega(\sqf_n \,\sqf_m)$, we have 
\beq
\hat\omega(n \star m) \,=\, \hat\omega(n) + \hat\omega(m) \,\,\,.
\eeq
Making the appropriate changes for the square-free numbers based on the familiar results for the  integers, (along the lines discussed in the previous Section), we can adjust the probabilistic arguments presented in  \cite{Schroeder} to get an estimate of $\hat\omega(n)$ for large $n$ and the argument goes as follows. 

First of all, we can convert the sum on those primes that divide $\sqf_n$ into a sum of {\em all} primes up to $\sqf_n$, using the probability factor $1/(p_i+1)$ that a prime $p_i$ occurs in a square-free number
\beq
\hat\omega(n) \sim \sum_{p_i \leq \sqf_n} \frac{1}{p_i +1} 
\,\,\,.
\label{firstestimate}
\eeq
The latter expression, using now the probability (\ref{squarefreeW}), can be further converted into a sum on all square-free numbers up to $\sqf_n$ and therefore used to compute the 
average of this quantity.
\begin{eqnarray}
&& \overline {\hat \omega(n})\,\sim \, 
 \, \int_2^{\sqf_n} \frac{dx}{(x+1) \log (x+1)} \label{secondestimate} \\
&& \,\approx \log\left(\log\left(\frac{\pi^2}{6} n +1 \right)\right) + \log\log 3\,\,\,.\nonumber 
\end{eqnarray}
where we have used $\sqf_n \sim \pi^2 n/6$. 
A better way to estimate $\overline{\hat\omega(n)}$ consists in summing up to a prime $p_m$ and convert the rest of the sum into an integral 
\beq
\overline{\hat\omega(n)} \,\sim \,\sum_{p_i=2}^{p_m} \frac{1}{p_i+1}+ \sum_{x=p_m}^{\sqf_n} \frac{1}{(x+1) \log (x+1)} \,\,\,.
\eeq
If we now approximate the second sum by an integral we have 
\beq 
\overline{\hat\omega(n)} \,\sim\,\sum_{p_i=2}^{p_m} \frac{1}{p_i+1} + \log\left(\log\left(\frac{\pi^2}{6} n +1 \right)\right)
- \log\left(\log\left(p_m +1 \right)\right)\,\,\,.
\eeq
Therefore the value of the constant $C$, defined as 
\beq
\overline{\hat\omega(n)} \equiv \log\left(\log\left(\frac{\pi^2}{6} n +1 \right)\right) + C \,\,\,,
\eeq
is obtained by taking this limit 
\beq
C \,=\, \lim_{p_m \rightarrow \infty} \left(\sum_{p_i=2}^{p_m} \frac{1}{p_i +1} - \log\left(\log\left(p_m +1 \right)\right)\right) \,\,\,.
\eeq
For estimating its value, first we expand the first term in powers of $1/p_i$ as 
\beq
\sum_{p_i=2}^{p_m} \frac{1}{p_i+1} \,=\,\sum_{s=0}^{\infty}  (-1)^s \sum_{p_i=2}^{p_m} \frac{1}{p_i^{s+1}}\,\,\,,
\eeq
and then we write $C$ as 
\beq 
C \,=\, {\mathcal A} - {\mathcal B}\,\,\,,
\eeq
where
\begin{eqnarray*}
{\mathcal A}\,&=& \,\left[\lim_{p_m \rightarrow \infty} \left(\sum_{p_i=2}^{p_m} \frac{1}{p_i } - \log\left(\log\left(p_m  \right)\right)\right)\right] \,\,\,;\\
{\mathcal B} \,&=&\,  \sum_{k=2}^{\infty}(-1)^k \sum_{p_i=2}^\infty \frac{1}{p_i^k}\,\,\,.
\end{eqnarray*}
The constant ${\mathcal A}$, also known as the  Kronecker constant, 
  can be expressed as  \cite{Schroeder}
\beq
{\mathcal A} \,=\, \gamma_E + \sum_{k=2} \frac{\mu(k)}{k} \log\zeta(k) \,\sim 0.261497...
\eeq
where the numerical value is easily extracted since the last series converges very quickly to its asymptotic value. The constant ${\mathcal B}$ can be numerically estimated 
in terms of the probability $U(x)$ (\ref{probabilitynatural}) of the primes along the integers as 
\beq
{\mathcal B} \,=\, \sum_{k=2}^{\infty}(-1)^k \sum_{p_i=2}^\infty \frac{1}{p_i^k}\,\,\sim \,\,\sum_{k=2} ^\infty \int_2^\infty \frac{1}{t^k\,\log t} \, dt \,\sim \,0.291479..
\label{formulaabove}
\eeq
Hence, put everything together, the theoretical estimate of the average of the number of prime divisors for square-free numbers is given by
\beq
\overline{\hat\omega(n)} \,\sim \,\log\left(\log\left(\frac{\pi^2}{6} n +1 \right)\right)  -0.029982 \,\,\,.
\label{theoreticalomega}
\eeq
{The estimate done above obviously leads, for the asymptotic behaviour of $\overline{\hat\omega(n)}$, to the prediction 
\beq \overline{\hat\omega(n)} \sim \log\left(\log\left(\frac{\pi^2}{6} n \right)\right)\,\,\,,
\eeq
but in the formula (\ref{theoreticalomega}) given above we prefer to keep explicitly present also the last constant term because, in comparing with some actual counting, for the very slow $\log\log$ behaviour of the first term, this constant can 
significantly improve the agreement with the data. This can be clearly seen in Table \ref{tableomega} for $x$ in the interval $x \in (1.0\times 10^7,1.0 \times 10^8)$, where the two quantities (the theoretical prediction (\ref{theoreticalomega}) and the actual counting agree one to the other for a relative error less than $10^{-3}$.  
} 

\begin{table}[t]
\begin{center}
\begin{tabular}{|l|l|l|l|}\hline
 \, n  &  
$ \overline{\hat\omega(n)}$ & $\overline{\hat \omega(n)}_{th}$ & relative error \\\hline
$1\, \times10^7 $  &\, $2.789$ & \, $2.780$ & $3.2\times 10^{-3}$\\ 
$2\, \times 10^7 $  & \, $2.828 $ & \, $2.821$ & $ 2.4 \times 10^{-4}$\\
$3 \, \times 10^7 $  & \, $2.850 $ &\, $2.844 $ & $2.1 \times 10^{-3}$\\
$4\, \times 10^7$  & \, $2.865$ & \, $2.860$ & $ 1.7 \times 10^{-3}$\\
$5 \, \times 10^7$ &  \, $2.877$ &\, $2.873$ & $ 1.3 \times 10^{-3}$  \\
$6 \, \times 10^7 $ & \,  $2.886$ &\, $2.883$ & $1.0 \times 10^{-3}$\\
$7 \, \times 10^7 $  & \, $2.894$ & \, $2.891$& $ 1.0 \times 10^{-3}$\\ 
$8\, \times 10^7 $  & \, $2.903$ & \, $2.901$ & $ 6.8 \times 10^{-4}$\\ 
$9\, \times 107 $ &  \, $2.907$ & \, $2.904$ & $ 1.0\times 10^{-3} $\\
$1\,\times 10^8$ & \, $2.913$ & \, $2.910$ & $ 1.0 \times 10^{-3}$\\
$2\,\times 10^8$ & \, $2.950$ &\, $2.949$ & $ 3.4 \times 10^{-4}$ \\
\hline
\end{tabular}
\end{center}
\caption{
 Numerical determination of the average of number of divisors for square free numbers, defined as $\overline{\hat \omega(n)} = 1/n \sum_{k=1}^n \hat\omega(k)$, vs its theoretical prediction, given in eq.\,(\ref{theoreticalomega}). In the last column the relative error of the theoretical prediction, i.e. $(\hat\omega(n) - \hat {\omega}(n)_{th})/\hat{\omega}(n)$.}
 \label{tableomega}
\end{table}

\subsection{Poisson distribution of ${\bf \hat\omega(n)}$}

It is well known (see, for instance, \cite{Wintner,GranvillePoisson}) that, denoting by $\pi_k(x)$ the number of integers less than $x$ having exactly $k$ prime factors (not necessarily distinct), these functions present a Poisson distribution
\beq
\pi_k(x) \sim \frac{x}{\log x} \frac{\tilde\lambda^{k-1}}{(k-1)!} \,\,\,,
\label{pikx}
\eeq
where the average value $\tilde\lambda$ is given by 
\beq
\tilde\lambda = \log\log x\,\,\,.
\label{avvvpoiss}
\eeq
Notice that $\pi_1(x)=\pi(x)$, where $\pi(x)$ is the familiar prime counting function. Moreover, summing on all values of $k$, we get correctly
\beq
\sum_{k=1}^\infty \pi_k(x) \,=\, x \,\,\,.
\label{normmmmpoiss}
\eeq
In view of these results, we expect to find a Poisson distribution also for number of divisors of the square-free numbers. Let's make some simple calculations which support 
such an expectation. 

 Let ${\mathcal O}_i(\sqf_n)$ be a binomial variable which takes values $1$ if the prime $p_i$ occurs in the prime decomposition of a square-free number $\sqf_n$ and $0$ otherwise. Since the probabilities that a prime factor $p_i$ occurs or not in a random square-free number are $\tilde q_i=1/(p_i+1)$ and 
$\overline q_i = 1-\tilde q_i = \left(1-1/(p_i+1)\right)$ respectively, the average of ${\mathcal O}_i$ is 
\beq
m_i \equiv \, \langle {\mathcal O}_i (\sqf_n)\rangle \,=\,\tilde q_i\,=\, \frac{1}{p_i +1} \,\,\,. 
\eeq
As a binomial variable, the variance of ${\mathcal O}_i(\sqf_n)$ is 
\beq
\sigma^2_i\equiv \langle ({\mathcal O_i(\sqf_n)} - m_i)^2 \rangle  \,=\, \tilde p_i \,\overline q_i \,=\, \frac{1}{p_i+1} \,\left(1-\frac{1}{p_i +1}\right) \,=\, 
m_i - \frac{1}{(p_i+1)^2} \,\,\,.
\label{varianceee}
\eeq
Assuming that the divisibility by different primes are independent and summing up to the $n$-th square-free number, for the overall average $\overline{\hat\omega(n)}$ we got 
the result of eq.\,(\ref{theoreticalomega}), while for the variance $\sigma^2$ we have 
\beq
\sigma^2\,\sim \,\overline{\hat\omega(n)} - \sum_{p_i=2}^{\infty} \frac{1}{(p_i+1)^2} \,\,\,,
\label{vvv}
\eeq
where the last sum was extended to infinity since it converges quickly. This series can be estimated by following the same steps as before, namely first we expand its coefficients in powers of $1/p_i$ and then we evaluate the sums on the inverse powers of the primes according to the the prime distribution along the integers
\begin{eqnarray}
&& \sum_{p_i=2}^\infty \frac{1}{(p_i+1)^2}\,=\,\sum_{k=1}^\infty (-1)^{k-1} k \sum_{p_i=2}^\infty \frac{1}{p_i^{k+1}} \,=\, 
\\
&&\,=\,
\sum_{k=1}^\infty (-1)^{k-1} k \int_2^\infty \frac{1}{t^{k+1}\,\log t} \, dt \,\sim 0.226978 \,\,\,.\nonumber 
\end{eqnarray}
Therefore 
\beq
\sigma^2 \sim \overline\omega(\sqf_n) - 0.226978.. \,\,\,.
\eeq
So, for very large $n$, we have that $\sigma^2 \sim \overline{\hat\omega(n)}$, one of the basic properties of a Poisson distribution, and based on previous results for a generic integer (see eq.\,(\ref{pikx}), it is natural to assume that $\hat\omega(n)$ is also approximately distributed 
according to the Poissonian distribution 
\beq
{\rm Prob}\,\left\{\hat\omega(n)\, =\,  k\right\}\,=\, {\mathcal P}(k,n)\,\sim \frac{\,\,\,\lambda^{k-1}}{(k-1)!} \, e^{-\lambda} \,\,\,,
\label{Poissondistribution}
\eeq
where\footnote{We have shifted by $1$ the distribution since each square-free number has always at least {\em one} prime factor.}
\beq 
\lambda \,=\, \overline{\hat\omega(n)} -1\,\,\,.
\eeq
According to the Poisson distribution, the most probable value of $k$ is given by 
\beq
\tilde k \,=\, [\overline{\hat\omega(n)}] +1 \,\,\,,
\label{massimoPois}
\eeq
where $[..]$ is the integer part, a prediction which is well satisfied even for relatively small number of square-free numbers employed. However, one should be aware that, for finite $n$, we should expect deviations of the actual distribution of the $N_k(n)$ from the Poisson distribution: first of all, for any finite $n$ there is only a {\em finite} number $q(n)$ of $k$ and not an infinite number of them, as required instead by the Poisson distribution; moreover, given the very slow $\log\log \pi^2/6 n$ dependence of $\overline{\hat\omega(n)}$, one must go to values of $n$ greater than, say, $n =10^{10^{50}}$ in order to satisfy the assumption we made that $\sigma^2 \sim \overline{\hat\omega(n)}$. Notice, however, that for $n \rightarrow \infty$ the value 
\beq
{\mathcal P}(1,n) \,=\,e^{-\lambda} 
\,\,
\underset{n \rightarrow \infty}{\sim}
\,\,
\frac{1}{\log\left(\frac{\pi^2}{6} n\right)}
\,\,\,
\label{limite}
\eeq
coincides with the large $n$ limit of the the correct probability (\ref{squarefreeW}) to find primes along the sequence of square-free numbers.

\subsection{Erd\H{o}s-Kac theorem for square-free numbers} Notice that in the limit $\sqf_n \rightarrow \infty$, the mean value of the Poisson distribution 
(\ref{Poissondistribution}) goes to infinity and, it is well known, that when this happens the Poisson distribution therefore tends to a gaussian (see, for instance, \cite{Feller}). Hence, we expect that the probability distribution of the number of distinct prime factors
 $\hat\omega(n)$ of the square-free number $\sqf_n$ is, in the large $n$ limit, the standard normal distribution ${\cal N}_{0,1}$, with mean $0$ and variance $1$ in the 
 variable 
\beq
\frac{\hat\omega(n) - \log\log \left(\frac{\pi^2}{6 }n\right)} {\sqrt{ \log\log \left(\frac{\pi^2}{6 }n\right)}}
\eeq
This result may be considered as the Erd\H{o}s-Kac theorem for square-free numbers (for a different derivation, see also Appendix A of the paper \cite{KLAGSBRUN}). 
It is of course reassuring that this expression coincides with what it could have been guessed on the basis of the density of square free numbers in the integers.

\subsection{Moments of the Mertens function} 
\begin{figure}[b]
\centering\includegraphics[width=1.\textwidth]{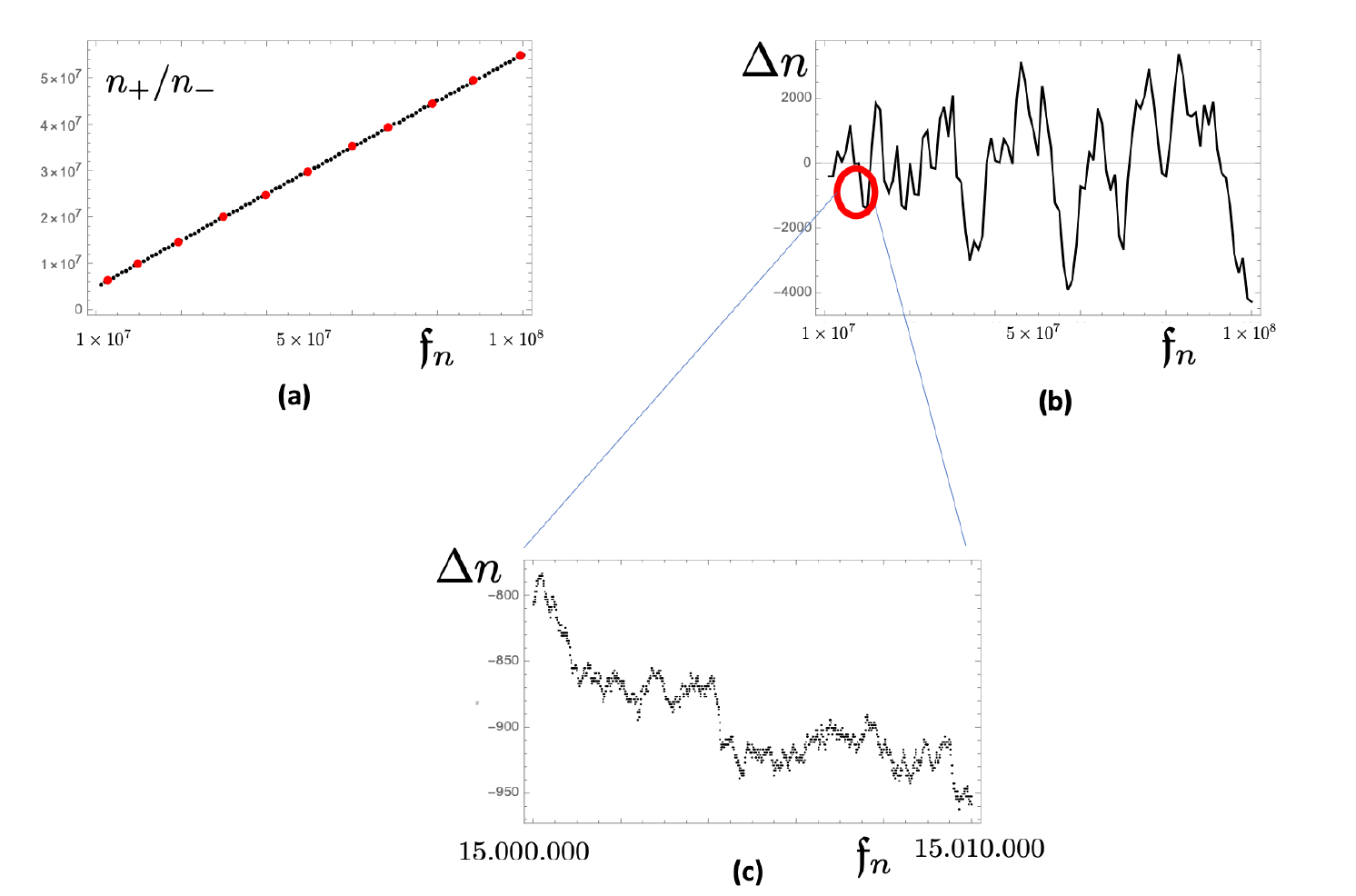}
\caption{(a) Plot of $n_{-} = N_1 + N_3 + N_5 +N_7$ (black dots) and $n_{+} = N_2 + N_4 +N_6 +N_8$ (red dots) in the interval $10^7 \leq \sqf_n \leq 10^8$. 
(b) Plot of $\Delta n = n_{-} - n_{+}$ in the interval $10^7 \leq \sqf_a \leq 10^8$; (c) Zoom on the behaviour of $\Delta n$ in the interval  
 $15.000.000 \leq \sqf_n \leq 15.010.000$. }
\label{spectroscopya}
\end{figure}

Let's now come back to the Mertens function. As seen in the previous Section, it can be written in terms of the $N_k(n)$ as 
\beq
\hat M(n) \,=\, \sum_{k=1}^{q(n)} (-1)^k \, N_k(n) \,\equiv\,n_{+} - n_{-}\,\,\,, 
\label{rewritingMertens2}
\eeq
where we have defined
\beq
\begin{array}{l}
n_{+}\,=\, N_2(n) + N_4(n) + N_6(n) + \cdots \\
n_{-}\,=\, N_1(n) + N_3(n) + N_5(n) + \cdots 
\end{array}
\eeq
Varying the upper limit $n$ of the Mertens function, each $N_k(n)$ is a fluctuating quantity. However, their coarse graining behavior corresponds to very 
{\em smoothly} increasing functions of $n$ (see Figure \ref{evolution}.a), with a mean value asymptotically determined by the Poisson distribution (\ref{Poissondistribution}) 
\beq
N_k(n) \rightarrow  n \, {\mathcal P}(k,n) \,\,\,,
\label{meanN_k}
\eeq
where 
\beq
n \,=\, \sum_{k=1}^{q(n)} N_k(n) \,=\, n_+ + n_-\,.
\label{total}
\eeq
This implies that also $n_{+}$ and $n_{-}$ have a smooth coarse grained behavior, with their values which are always almost coincident (see Figure \ref{spectroscopya}). Hence, the non-zero values of the Mertens function come only from the {\em fluctuations} of $\Delta n = n_{+} - n_{-}$, which are typically of order $\sqrt{n}$. 

The argument we just presented can be easily made more precise employing the Poisson distribution (\ref{Poissondistribution}). Indeed, assuming that in the limit $n\rightarrow \infty$ the $N_k(n)$ satisfy the Poisson distribution (\ref{Poissondistribution}), let's ask the following question: what are the mean value and the variance of the quantity $\hat M(n)$, as expressed in eq.\,(\ref{rewritingMertens2})? For the answer we can argue as follows:

\begin{figure}[t]
\centering\includegraphics[width=0.7\textwidth]{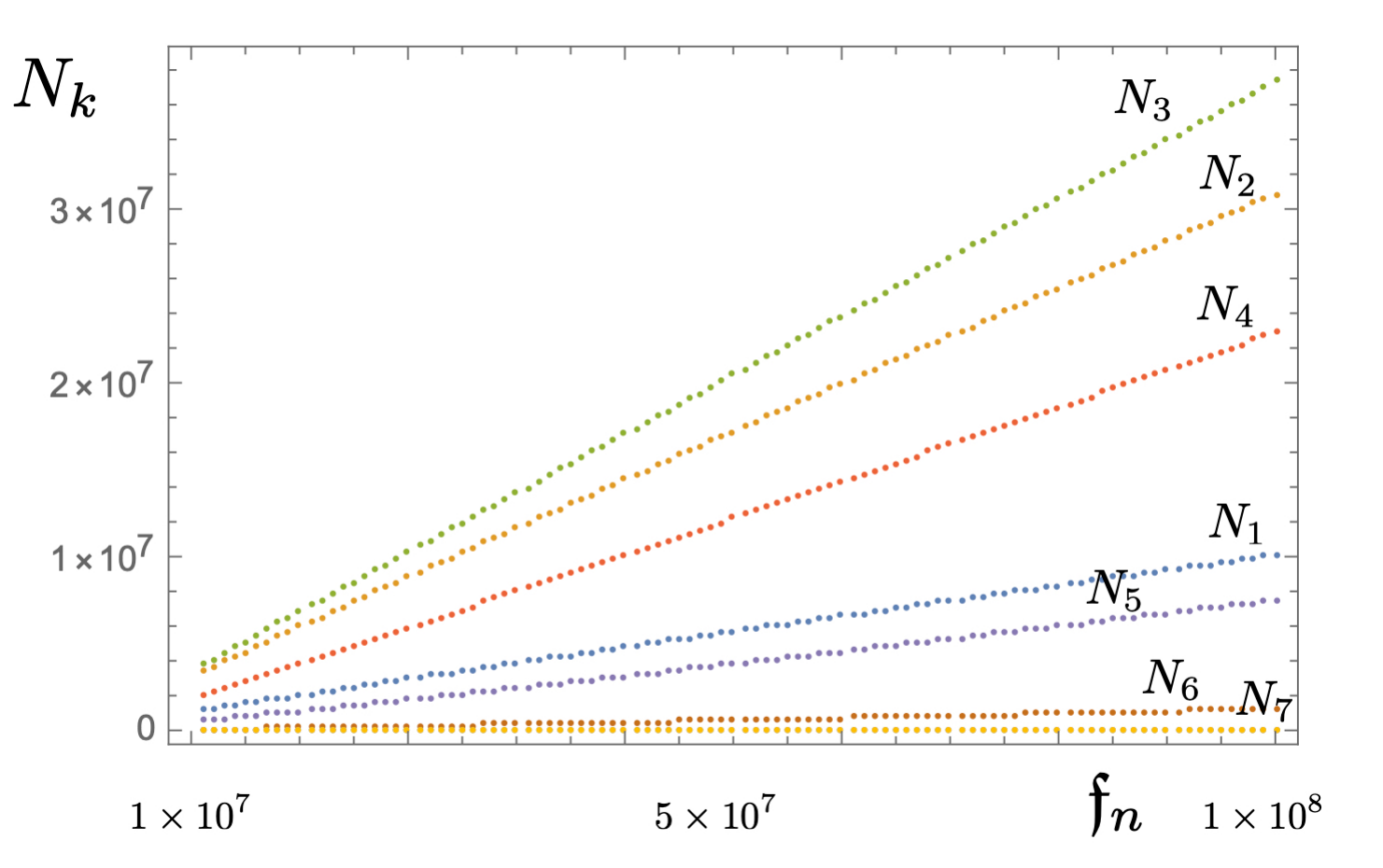}
\caption{The behavior of $N_k$ ($k=1, 2, \ldots 7$) for the square-free number in the interval $10^7 \leq \sqf_n \leq 10^8$. $N_8$ is not reported since it is  too small 
in the entire interval.}  
\label{evolution}
\end{figure}

\begin{itemize}
\item Using the Poisson distribution (\ref{Poissondistribution}), we can define a {\em binomial} variable $a_{\pm}$ with probabilities $P_{\pm}$ associated 
to the event to have a square-free number which belongs to {\em any} of the even or the odd quantities $N_k(n)$. These probabilities are given by 
\begin{eqnarray}
P_+ &\,=\,& \sum_{k=1}^{\infty} \frac{\lambda^{2k-1}}{(2 k-1)!} \,e^{-\lambda} \,=\,\frac{1}{2} \left(1 - e^{-2 \lambda}\right) \,\,\,;\\
P_- &\,=\,& \sum_{k=0}^{\infty} \frac{\lambda^{2k}}{(2 k)!} \,e^{-\lambda} \,\,\,\,\,\,\,\,=\,\,\frac{1}{2} \left(1 + e^{-2 \lambda}\right) \nonumber\,\,\,.
\end{eqnarray}
\item Making use of eq.\,(\ref{total}),  we can express $\hat M(n)$ equivalently as 
\beq
\hat M(n) \,=\, n - 2 n_- \,=\, 2 n_+ -  n \,\,\,,
\label{equivvv}
\eeq  
and therefore, for its mean value, we have  
\beq
\overline M\,=\,\langle \hat M(n) \rangle \,=\,\lim_{n\rightarrow \infty} \frac{1}{n} \left(2 \langle n_+\rangle  - n\rangle\right)\,=\,2\,P_+ -1 \,=\, -e^{-2 \lambda} \rightarrow 0 
\eeq
since $\lambda \rightarrow \infty$ when $n\rightarrow \infty$. 
\item Let's now compute the variance of $\hat M(n)$ using the standard result for a binomial variable 
\begin{eqnarray}
\lim_{n\rightarrow \infty} \langle (\hat M(n) - \overline M )^2 \rangle &\,=\,& \lim_{n\rightarrow \infty} 4 \left(n_+  - \frac{n}{2}\right)^2 \,=\, \nonumber \\
& = &4\, n\, P_+ \, P_-  \,=\, n\, \left(1 - e^{-4\lambda}\right) \,\rightarrow \, n  
\end{eqnarray}
\item By the same token, in the limit $n \rightarrow \infty$ we can compute all higher moments $m_k =(\hat M(n) - \overline M)^k$. 
The odd moments vanish, $m_{2 n +1} =0$ while for the even ones we have  
\begin{eqnarray}
 m_2 & = & n \nonumber \\
 m_4 & = & (3 n -2) n \\
 m_6 & = & (15 n^2 - 30 n + 16) n\nonumber  \\
 m_8 & = & (105 n^3 - 420 n^2 + 588 n - 272) n \nonumber \\
 \ldots & & \ldots\nonumber 
 \end{eqnarray}
 The leading order in $n$, for $n \rightarrow \infty$, gives $m_{2 k} \sim n^k (k-1)!!$, which are precisely the moments of a normal distribution ${\mathcal N}_{0,n} $. 
\end{itemize}

\vspace{3mm}
\noindent
Hence, using the Poisson distribution of the variables $N_k$ we can conclude that, for $n \rightarrow \infty$
\begin{center}
\fbox{
\parbox{15.7cm}{\beq
\langle (\hat M(n))^{2 k+1}\rangle \,=\, 0 
\,\,\,\,\,\,\,\,\,\,\,\,\,\,\,\,
;
\,\,\,\,\,\,\,\,\,\,\,\,\,\,\,\,
\langle (\hat M(n) )^{2 k} \rangle\,=\, n^k \, (k-1)!! 
\eeq
}}
\end{center}
All these considerations lead to the law of iterated logarithms for the variable $\hat M(n)$, considered as a  ``random variable" subject to a normal law distribution
\beq
\lim_{n\rightarrow \infty} \,{\rm sup} \frac{\pm \hat M(n)}{\sqrt{2 n \log\log n}} \,= \, 1 
\hspace{3mm} 
, 
\hspace{3mm}
{\rm a. s. } 
\label{sqrtlogloglog}
\eeq
where ``a.s.'' stands  for ``almost surely'', in  a probabilistic sense \cite{Feller}. 

\section{Summary} In this Part B of the paper, adopting a probabilistic point of view for estimating the restricted Mertens function $\hat M(n)$ in terms of the fluctuation in the numbers $n_+$ and $n_-$ of the square-free numbers with signature $\pm 1$ in their M\"{o}bius coefficients, we have been able to argue that in the limit $n \rightarrow \infty$ these fluctuations can be described by a normal distribution (of zero mean and variance $\sigma^2 \rightarrow n$). If this is indeed the case, we can then apply the law of iterated logarithms for constraining the growth of 
$\hat M(n)$, as in eq.\,(\ref{sqrtlogloglog}). In the course of our argument we have made use of some approximations for evaluating certain number theory functions but it has been rewarding to see 
that these approximated quantities turn out to be very close to their exact values obtained by enumeration. Therefore it becomes quite interesting to check in more detail the emergence of a normal law distribution for the fluctuation of the restricted Mertens function $\hat M(n)$. This is precisely the scope of the Part C of the present paper.

\newpage

\addcontentsline{toc}{section}{Part C}
\begin{center}
{\Large {\bf PART C}}
\end{center}

Part C may be considered the ``experimental section" of the paper, carried on for checking explicitly the hypothesis of randomness of the M\"{o}bius sequences and for the pleasure of finding things out. 
Hereafter the randomness refers to variables which can take only two values (chosen to be either $\pm 1$ or $(0,1)$) with equal probability.
The point of view adopted in the following is interesting {\em per se}: namely, in this part of the paper we pretend to have no a-priori knowledge about the origin of our original sequence $\{\mathcal {S}_n\}$ and our goal will be to infer its true nature by applying a series of statistical tests which become more and more refined. In particular, our aim will be to see whether such a sequence can be considered for all purposes as a truly random series or, vice versa, if there are some hints that signal its non-random nature through the identification of some periodic or other deterministic features. For doing so, we will compare arbitrarily large subsequences\footnote{These subsequences will be often denoted with the notation $\{\epsilon\}$.} of our sequence $\{\mathcal {S}_n\}$ with sequences made of 
ideal random variables with two equiprobable values and will check the similarities of their behavior. Since randomness is a probabilistic property, we are going to characterize such a behavior in terms of probability. When applied to a truly random sequence, the theory of probability tells us the likely outcomes of the statistical tests and this is the basis of our comparisons. In other words, for a given quantity $G$, we compute the probability of its values assuming the true random nature of the binary variables it depends upon and we compare the statistical frequencies of the values of $G$ evaluated on very large subsequences of $\{\mathcal {S}_n\}$. In the next sections we will apply several statistical tests to subsequences extracted from the sequence $\{\mathcal{S}_n\}$ and each statistical test will check the presence or the absence of certain patterns: its output could indicate whether the sequence can be considered as random or not. 



Let's anticipate that the successful and impressive performance of {\bf ALL} the statistical tests we have applied to the sequence $\{\mathcal{S}_n\}$ leads indeed to the conclusion that 
this sequence behaves indeed as a perfect random number generator, with a confidence level \footnote{With enough computer memory, one can increase systematically the confidence level of our tests to $99.9\%$  by 
simply enlarging the size of the sampling of our analysis. As a matter of fact we have done so  for some of our tests,  but not for all.}
of $99\%$. 
We have carried out a huge number of tests in an interval of values of the sequence $\{\mathcal{S}_n\}$ as large as
$$
1 \leq n \leq 10^{16}\,\,\,,
$$
although we have mostly focused our attention on sequences with $n > 10^6$ for avoiding possible atypical behavior of the restricted Mertens function for small values of the index $n$.   
As explained in more detail below, the study was mostly done by analyzing blocks made of $160.000.000$ values and, for all the blocks we have analyzed, we have always 
observed similar results, {\em independently} of the initial value of the blocks. In other words, the random nature of the restricted M\"{o}bius coefficients $\hat\mu(n)$  appears to be a translation invariant feature of their sequence. To the best of our knowledge this is the first time that the random nature of the restricted M\"{o}bius coefficients has been tested so thoroughly and extensively.  The statistical tests we used are mainly those suggested by Donald Knuth in Volume 2 of his {\em The Art of Computer Programming} \cite{Knuth}, those listed in the report of the National Institute of Standards and Technology (NIST) \cite{Nist} or those coming from the battery tests Diehard developed by George Marsaglia \cite{Marsaglia}, together with others which come from other sources. Before presenting our main findings, it is important to set some notation and some basic facts of probability. To this aim, we also warmly invite the reader to consult the NIST report \cite{Nist} for details on the general 
discussion of these protocols since many considerations which follow, in particular those relative to the P-values and level of significance (see Sect. \ref{PPPPVALUES} below), will become {very} clear upon reading the NIST report.  

\section{Statistical Test Suite}\label{StatisticalSuite}

In the following, according to the nature of the test, we will employ either the sequence $\{\mathcal{S}_n\}$, having values $(0,1)$, or the sequence of the restricted M\"{o}bius coefficients $\{\hat\mu_n\}$ which we recall are related to $\{\mathcal {S}_n\}$ by the relation 
\beq
\hat\mu_n \,=\, 2 \, \mathcal{S}_n -1 \,\,\,. 
\label{XXX}
\eeq 
The nature of all our tests will be binary. Namely, their protocol consists in checking a {\em null hypothesis} $H_0$ vs an {\em alternative hypothesis}
which in our case consist of 

\vspace{1mm}

\begin{center}
\fbox{
\parbox{14.7cm}{
 {\bf Null Hypothesis}: the sequence $\{{\mathcal S}_n\}$ 
 is random. \\
\\
 {\bf Alternative Hypothesis}: the sequence $\{{\mathcal S}_n\}$  
 is not random.
}}
\end{center}

\vspace{1mm}
\noindent
For each test, we need to reach a conclusion  deciding whether to accept or reject the null hypothesis. This will be done according to the following probabilistic considerations. 

\subsection{$\chi^2$  distribution}

For each test, we have a random reference scenario which predicts the statistics of all possible values, with the relative probabilities. Suppose that $n$ observations in a random sample from a population are classified into $r$ mutually exclusive classes with respective observed numbers $Y_s$ (for $s = 1,2,\ldots,r$).  According to the null hypothesis, let $p_s$ be the probability that each observation falls into the output $s$ and let $Y_s$ be instead the number of experimental observations that actually do fall into the output $s$. So. the expected numbers are $n p_s$ for all $s$ and we can form the statistical quantity 
\beq
\chi^2 \,=\, \sum_{s=1}^r \frac{(Y_s - n p_s)^2}{n p_s} \,\,\,,
\label{chisquare}
\eeq 
called $\chi^2$ with $k\equiv (r-1)$ degrees of freedom: the degrees of freedom are one less than the total number of possibilities because the total counts have to sum to $n$. Such a quantity clearly measures how the observed numbers of events differ from the theoretical ones. Such a quantity will be defined for each test considered in the following. The important question concerns the distribution of $\chi^2$ values and, more importantly, what constitutes a reasonable value of $\chi^2$. 

The theoretical answers to these questions relies on the $\chi^2$-distribution \cite{Knuth,Nist} which states that the probability density function for the $\chi^2$ distribution with $k$ degrees of freedom is given by 
\beq
P_k(x) \,=\,\frac{1}{2^{k/2} \Gamma(k/2)} \, x^{k/2 -1} \, e^{-x/2}
\label{chisquaredis}
\hspace{5mm}
,
\hspace{5mm} x\geq 0\,\,\,.
\eeq
The easy way to understand this distribution is to think of  a random sampling of quantities taken from a normal distribution:  then, the $\chi^2$ distribution is nothing else but the distribution of the sum of the squares of these random samples. The degrees of freedom $k$ in a $\chi^2$ distribution is also its mean and half of its variance since we have 
\beq
\langle \chi^2 \rangle \,=\, \int_0^\infty \chi^2 \, P_k(\chi^2)\, d\chi^2 \,=\, k 
\hspace{5mm}
;
\hspace{5mm}
\langle (\chi^2 - k)^2 \rangle \,=\,\int_0^\infty \, (\chi^2 -k )^2\, P_k(\chi^2) \, d\chi^2 \,=\, 2 k \,\,\,. 
\eeq
These expressions tell us what is the most reasonable value of $\chi^2$, alias $\chi^2 \sim k$, and its normal fluctuation range, $\chi^2 \sim k\pm \sqrt{2 k}$. Notice that the probability distribution (\ref{chisquaredis}) depends only on $k$ and not on the values $n$ and the probabilities $p_s$. However, it is useful to remind that the chi-square distribution is an approximation that is valid only for large enough values of $n$: the rule of thumb is that $n$ has to be large enough so that each of the expected values $n p_s$ is around $5$ or more. 

\subsection{Normal distribution} 
In some of our tests we will also employ the familiar normal distribution ${\mathcal N}_{\mu,\sigma^2}(s)$
\beq
{\mathcal N}_{\mu,\sigma^2} (s) \,=
\,
\frac{1}{\sqrt{ 2 \pi \sigma^2}} \, \exp\left[{- \frac{(s-\mu)^2}{2 \sigma^2}}\right] 
 \,\,\,,
\eeq
with mean $\mu$ and variance $\sigma^2$. This will happen when our test variable  takes  the form $(s - \mu)/\sigma$, where $s$ is the sample test statistic value while $\mu$ and $\sigma^2$ are the expected theoretical values of the mean and the variance of our test statistics. Notice that, for $r\rightarrow \infty$, the previous $\chi^2$-distribution tends to a normal distribution. 

\subsection{P-value and level of significance}\label{PPPPVALUES}

For each test, according to its nature, a relevant randomness statistic must be chosen (e.g. the $\chi^2$ or the normal distributions) and used to determine the acceptance or the rejection of the null hypothesis. Indeed, if the sequence under test is non-random, the calculated test statistic will fall into extreme regions of the reference distribution. 
This is equivalent to say that, from the reference distribution, we have chosen a {critical value} $\hat{\bf t}$,  
typically in the tail of the distribution, and compare it to the measured test statistical value: if the test statistic value exceeds the critical value, the null hypothesis for randomness is rejected, otherwise, the null hypothesis will be accepted. 
In order to quantify this concept, it is necessary to introduce the $P$-value which, under the null hypothesis of randomness, is the probability that the chosen test statistics will assume values that are equal or worse than the observed test statistic value when considering the null hypothesis. In other words,  the $P$-value is the probability of obtaining test results at least as extreme as the results actually observed, under the assumption that the null hypothesis is correct. The $P$-value is also called the {\em tail probability}. The $P$-values are always in the range $(0,1)$, with the rule of thumb that  

\vspace{3mm}

\noindent\fbox{\begin{minipage}[t][3\height][c]{\dimexpr\textwidth-2\fboxsep-2\fboxrule\relax}
\centering
{\bf  
The higher the $P$-value,  the better the statistical confirmation of the null hypothesis.  }

\end{minipage}}

\vspace{3mm}
\noindent
So, if the $P$-value for a test is determined to be equal to $1$, then the sequence appears to have perfect randomness. Viceversa, a $P$-value equal to $0$ indicates that the sequence is completely non-random. 

{In order to handle the intermediate cases, a {\bf significance level} $\alpha$ must be chosen for the tests, with the rule that if $P$-value $\geq \alpha$, then the null hypothesis of randomness is accepted, i.e. the sequence appears to be random. Vice versa, if $P$-value $\leq \alpha$, then the null hypothesis is rejected, i.e. the sequence appears to be non-random. Typically $\alpha$ is chosen in the range $(0.001,0.01)$. In the following we will mostly choose $\alpha$ to be $\alpha=0.01$. In order to clarify better this notion 
of significance level, let's quote direclty the statement written in the NIST report \cite{Nist}
\begin{itemize}
\item 
A value of $\alpha$ equal to $0.01$ indicates that one would expect $1$ sequence in $100$ sequences to be rejected. A $P$-value $\geq 0.01$ would mean that the sequence would be considered to be random with a confidence of $99\%$. Vice versa, a $P$-value $ < 0.01$ would mean that the conclusion is  that the sequence is non-random with a confidence of $99\%$. 
\end{itemize}
}

Concerning the statistical variables, in the rest of the paper we consider the  $P$-value coming either from the $\chi^2$ or from the gaussian distributions, according to the nature of the 
variable under scrutiny.  

\begin{figure}[t]
\centering\includegraphics[width=0.5\textwidth]{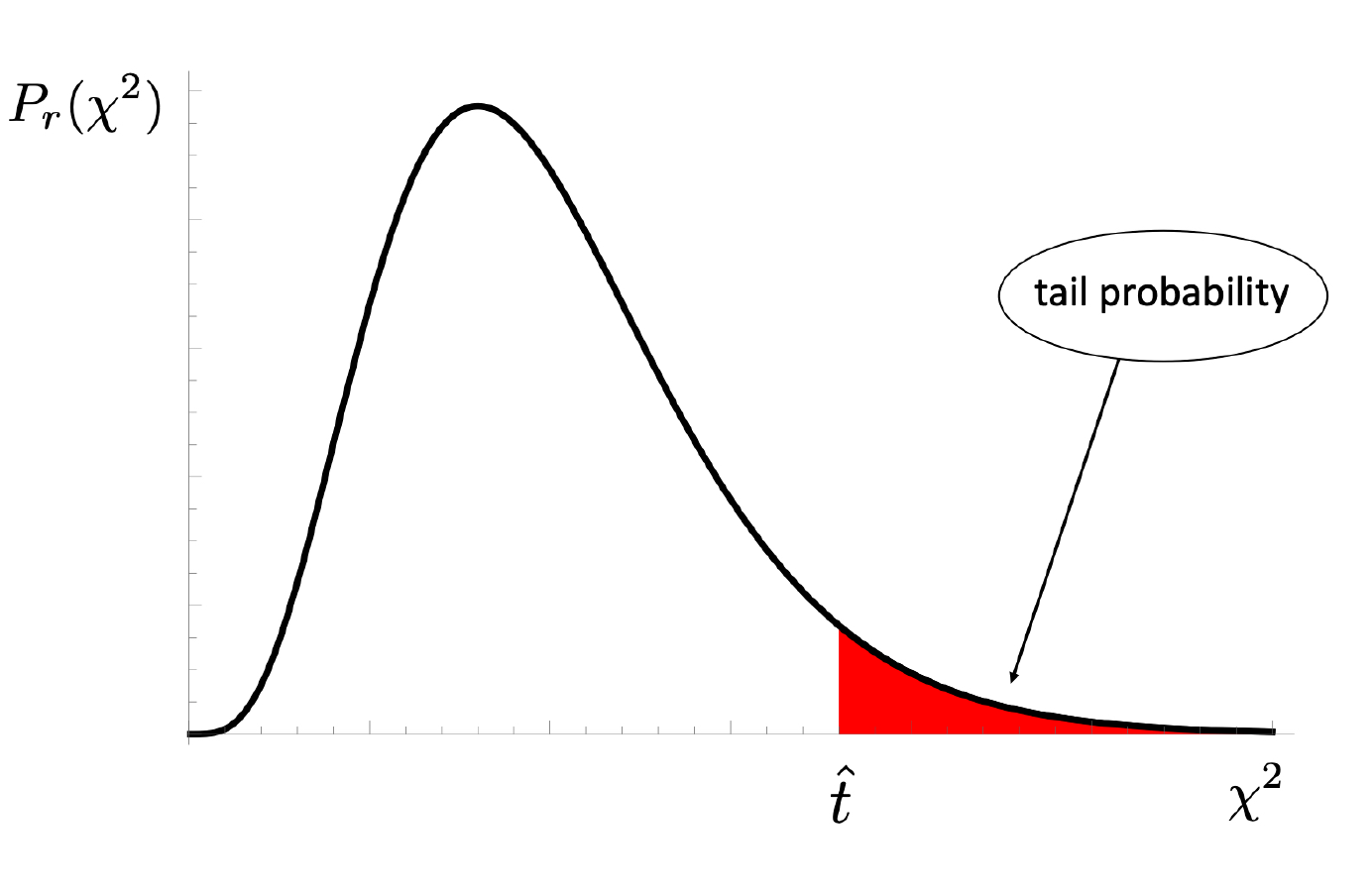}
\caption{The integral on the tail of the chi-square distribution, for values larger than a critical value $\hat t$, determines its relative $P$-value.}  
\label{tailintegral}
\end{figure}

\begin{figure}[b]
\centering\includegraphics[width=0.6\textwidth]{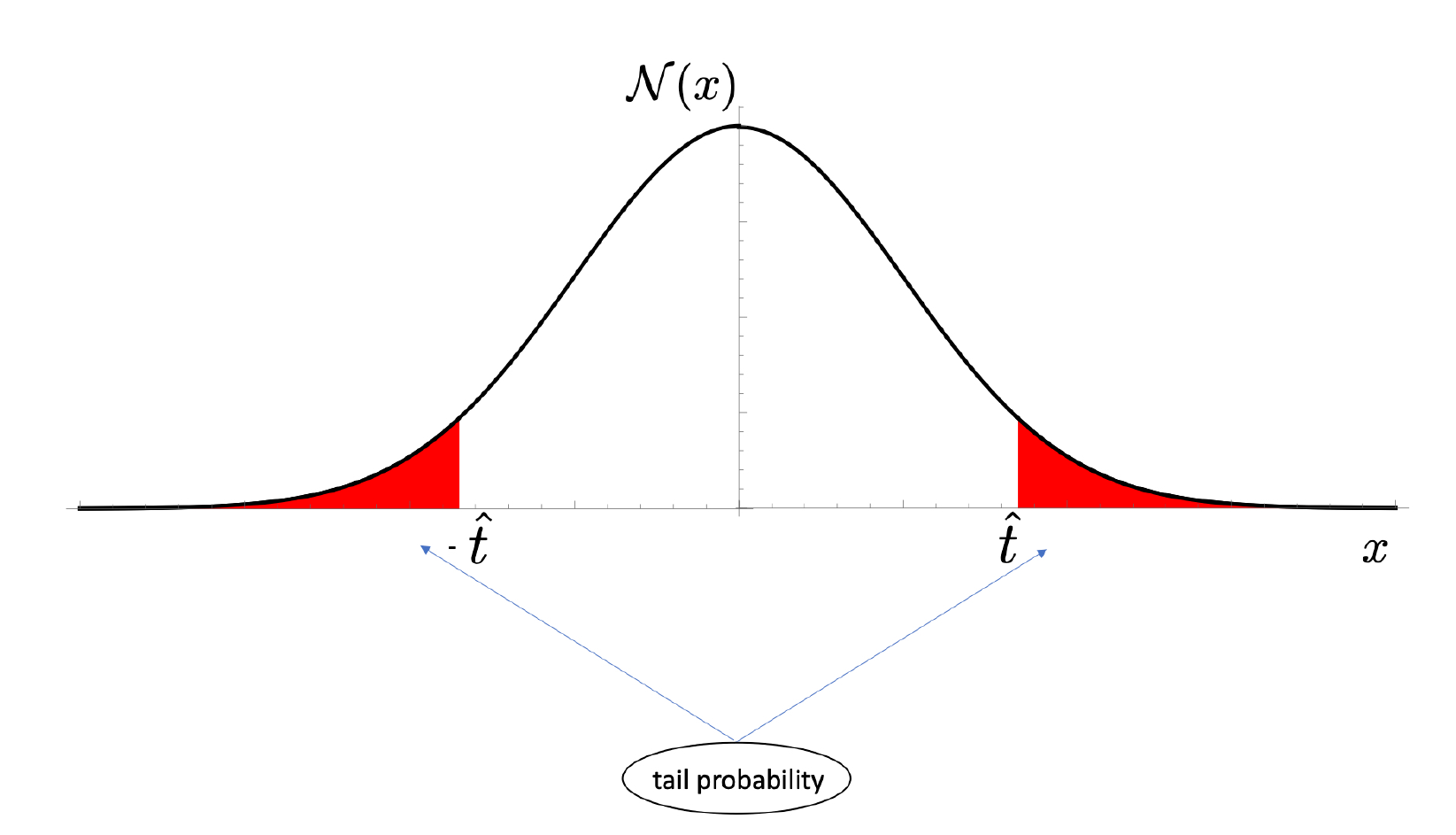}
\caption{The integral on both tails of the normal distribution, for values (in absolute term) larger than a critical value $\hat t$, determines its relative $P$-value.}  
\label{tailintegralgaussian}
\end{figure}

\begin{itemize}
\item For the $\chi^2$ distribution with $k$ degrees of freedom, which is defined only for positive values $x> 0$, the $P$-value is given by the integral on the tail of the distribution (see Figure \ref{tailintegral}) and expressed in terms of the function 
\beq
Q(a,z) \,=\,\frac{1}{\Gamma(a)} \int_z^\infty e^{-t} \, t^{a-1} \, dt \,\equiv \, \frac{\Gamma(a,z)}{\Gamma(a)} \,\,\,,
\label{qax}
\eeq
where $\Gamma(a,z)$ is the Incomplete Gamma Function, with $a = k/2$.   The $P$-value is equal to $Q(a,z)$ with parameters $a,z$ that depend on the application. 
\item For the normal distribution, the extreme values of the distribution are placed {\em both} at the left and the right parts of the curve (see Figure \ref{tailintegralgaussian}).  Therefore the $P$-value is relative to the two red areas and expressed in terms of the complimentary erf function: 
\beq 
{\rm erfc}(z) \,\equiv 1-{\rm erf} (z) \equiv \,\frac{2}{\sqrt{\pi}} \, \int_z^\infty e^{-u^2} \, du \,=\,\frac{\Gamma\left(\frac{1}{2},z^2\right)}{\Gamma\left(\frac{1}{2}\right)}\,\, ,
\label{erfc}
\eeq
\end{itemize}
\beq
\label{erfc2}
P(z) = {\rm erfc } (z/\sqrt{2}) .
\eeq

\subsection{Other statistical tests for continuous distributions}
It is useful to stress that some of our statistical tests presented in the next sections will involve random quantities that range over infinitely many values. In this case we will need to compare the ``experimental" distribution with the reference continuous distribution derived under the hypothesis of randomness for the variables involved. 
In addition to the $\chi^2$ test, the degree of validity of the reference continuous distribution will be certified by a statistical test such as the Kolmogorov-Smirnov test \cite{KolSmitest}, 
 which is a nonparametric test of the equality of continuous, one-dimensional probability distributions that can be used to compare a sample with a reference probability distribution. 
In addition to the Kolmogorov-Smirnov test, for testing the gaussianity of certain data we will also use other statistical tests, such as: 
 Anderson-Darling test \cite{ADtest}; Baringhaus-Henze test \cite{BHtest}; Cramer-von Mises test \cite{CMtest}; Jarque-Bera ALM test \cite{JBtest}; Mardia Combined test, Mardia Kurtosis  test, Mardia Skewness test \cite{Mardia1test}; Shapiro-Wilk test \cite{ShWitest}. We refer the reader to the original literature for the details of all these tests. For our purposes it is sufficient to say that all these tests are nowadays standardly implemented, for instance, in Mathematica \cite{Mathematica} and how good a given theoretical distribution is able to describe the numerical samples is, as usual, captured by a probability value (P-value) and, as a rule of thumb, a small P-value suggests that it is unlikely that the data follows the theoretical distribution.

\subsection{Single Brownian Trajectory Problem.} 
As in our previous studies on the Generalized Riemann Hypothesis \cite{ML,LM}, also in this case for the statistical analysis of the Mertens function we have to face the problem that there is  {\em one and only one} sequence $\{{\mathcal S}_n\}$. This poses a natural question: how can we then define a  statistical ensemble $\CE$ relative to the possible values ${\mathcal S}_n$ of the sequence and assign the relative probabilities to the associated Mertens function $\hat M(n)$? 
As is well known, this is a common problem in many time series, in particular for all those time series which refer to situations for which it is impossible to \textquotedblleft turn back time\textquotedblright. Indeed, in these cases, at any point of the sequence it is impossible to have access to all possible outputs (at the $k$-step there is as a matter of fact one and only one output, which in our case is just $\mathcal{S}_k$) and therefore it seems equally impossible to define the relative probabilities\footnote{It is the same problem to be faced in order to give meaning to the statement "the probability of rain tomorrow is $60\%$. We cannot recreate today's weather conditions hundreds of times and check how often it rains.}. 
In the literature, this is known as the {\em Single Brownian Trajectory Problem} (see, for istance \cite{brow1,brow2,brow3,brow4} and references therein). The solution of the problem passes through the {\em block variables}, as originally proposed in our previous papers \cite{ML,LM} and briefly discussed hereafter. 

\subsection{ Intervals and block variables}\label{intervalsss}
\begin{figure}[b]
\centering\includegraphics[width=0.6\textwidth]{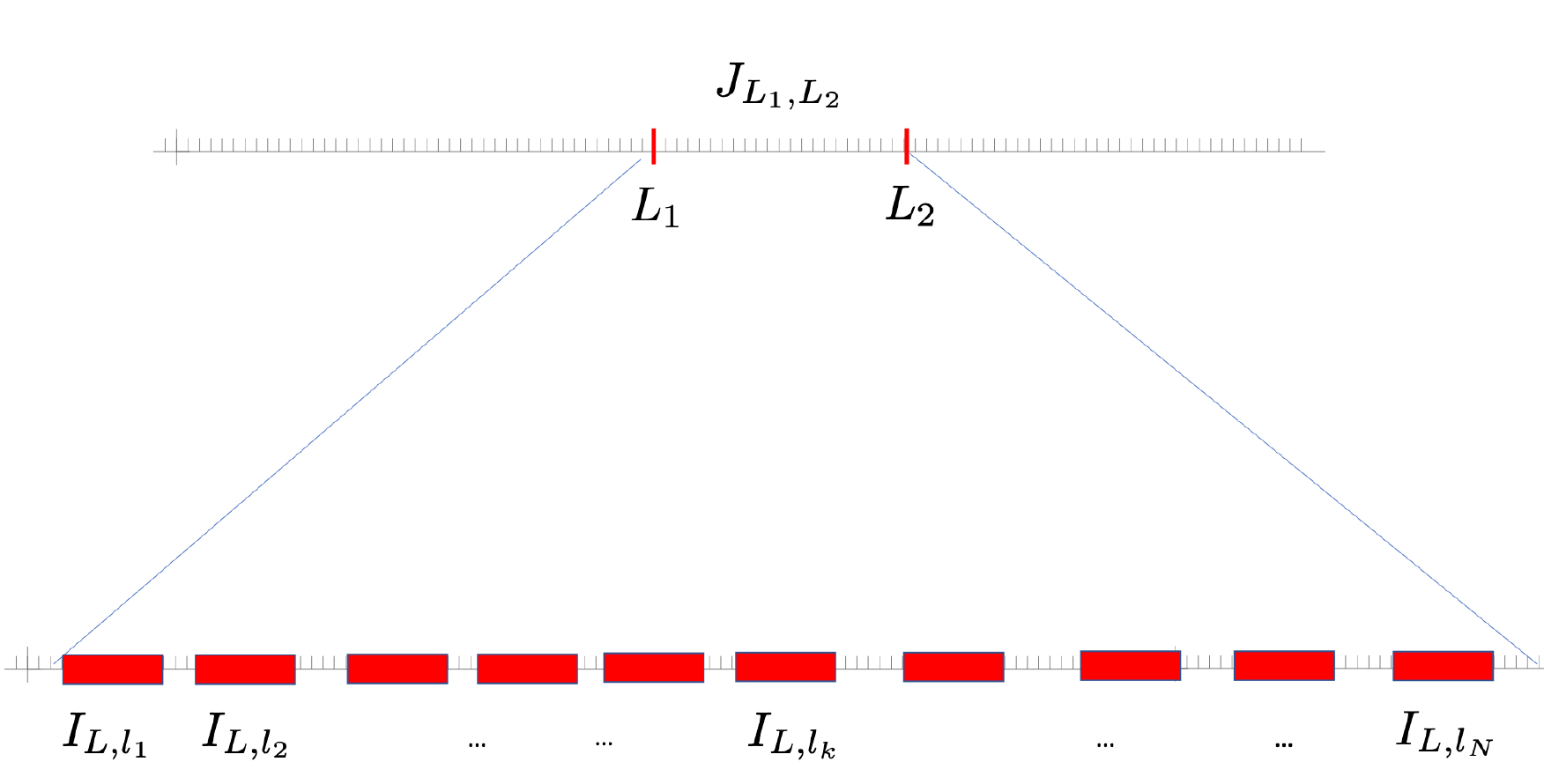}
\caption{The set of intervals $I_L (\ell_k)= I_{L, \ell_k} $ within the large interval $J_{L_1,L_2}$.}  
\label{subintervals}
\end{figure}

 In order to deal with the Single Brownian Trajectory Problem, imagine we have an arbitrarily long time series as our $\{\mathcal{S}_n\}$: to sample it, we take \textquotedblleft stroboscopic\textquotedblright \,snapshots of it. Let's, first of all, introduce some definitions: let $J_{L_1,L_2}$ (with $L_1 < L_2$) denote a set of consecutive integers given by 
 \beq
 J_{L_1,L_2} \,=\,\{L_1, L_1 +1, L_1 + 2, \ldots, L_2 -1, L_2 \}\,\,\,.
 \eeq
It is also useful to define ordered intervals of integers of length $L$ starting at a point $\s$  
\beq
I_L(\s) =\{\s, \s+1, \s+2, \ldots, \s+L -1 \} \,\,\,. 
\,\,\,
\label{intervals}
\eeq 
To sample the sequence $\{\mathcal{S}_n\}$, first of all we can use $J_{L_1,L_2}$ (with very large values of both $L_1$ and $L_2$) to select huge portions of values along the 
sequence $\{\mathcal {S}_n\}$; secondly, within these portions, we can use well separated intervals $I_L(l)$ to extract sub-sequences\footnote{Later we drop the lower indices and we denote these generic subsequences simply as $\{\epsilon\}$.} $\{\epsilon_{L,l}\}$ of the $\mathcal {S}_n$'s given by 
\beq
\{\epsilon_{L,l}\} \,=\, \{ \mathcal{S}_l, {\mathcal S}_{l+1}, \ldots, {\mathcal S}_{ \s+L -1}\}
\,\,\,.
\eeq 
The analogous subsequences for the restricted M\"{o}bius coefficients will be denote by $\{\hat \epsilon_{L,l}\}$, i.e. 
\beq
\{\hat \epsilon_{L,l}\} \,=\, \{\hat\mu_l, \hat\mu_{l+1}, \ldots, \hat\mu_{\s+L -1}\}\,\,\,.
\eeq
The subsequences $\{\epsilon\}$ or $\{\hat\epsilon\}$ are those sequence of numbers which will be submitted later to our statistical tests. The situation is illustrated in Figure \ref{subintervals}.

Concerning the restricted Mertens function itself, in correspondence with these intervals, we define the {\it block variables} $B_L (\s)$ as 
\beq
B_L(\s) \,=\,\sum_{k\in I_L(\s)}  \hat\mu(k) \,=\, \sum_{k=\s}^{\s + L-1}  \hat\mu(k)  \,\,\,. 
\label{groupvariables}
\eeq
These block variables are of course parts of the total restricted Mertens function which, using the same notation, can be written in fact as 
\beq
\hat M(n) \,=\,B_n(1)\,=\,\sum_{k=1}^n \hat\mu(k)\,\,\,. 
\eeq
Using the intervals $J_{L_1,L_2}$ together with $I_L (\ell)$,  we can also define many and well separated block variables $B_L(\s)$ of the {\em same length} $L$  and we can use them as members of the ensemble to which belongs the original sequence $B_n(1)$! This is equivalent to the stroboscopic snapshot procedure which is behind the solution of the Single Brownian Trajectory Problem (see the forthcoming subsection).  The validity of this self-averaging procedure relies on three assumptions on the corresponding time series which are indeed satisfied: its ergodicity, stationarity and scalability.
\begin{enumerate}
\item {\bf Ergodicity}. In the case of our subsequences $\{\hat\epsilon_{L,l}\}$, their ergodicity means the presence of {\em all} possible outputs $\pm 1$ of the restricted M\"{o}bius 
coefficients.
\item {\bf Stationarity}. The stationarity of the subsequences $\{\hat\epsilon_{L,l}\}$ means that, fixing a very large value of $L$, the output of our tests is independent of the starting point $l$. 
\item {\bf Scalability}. The scalability of the subsequences $\{\hat\epsilon_{L,l}\}$ means that we can always arbitrarily enlarge our sampling moving toward large values of $l$ in our subsequences, 
i.e. any test applied to the subsequences $\{\hat\epsilon_{L,l_1}\}$ can be equally applied to any other of such randomly chosen subsequences $\{\hat\epsilon_{L,l_2}\}$
with $l_2 > l_1$. In other words, there is no limit to the value of the starting point $l$ of the subsequences.   
 
\end{enumerate}

\subsection{Statistical Ensemble $\CE$ for  the series ${\bf \hat M(n)}$}\label{eenseble}
The block variables $B_L(\s)$ are the equivalent of the  
\textquotedblleft stroboscopic\textquotedblright \,images of length $L$ of a single Brownian trajectory (see Figure \ref{partitiontimeseriesGN}) and they allow us to control the irregular behavior of the original series $B_n(1)$ by proliferating it into a collection of sums of the same length $L$. It is this collection of sums that forms the {\em set of events}, i.e. the ensemble $\CE$ relative to the sums of $L$ consecutive terms $\hat\mu(k)$.  This ensemble is defined as follows \cite{ML,LM}: 

\begin{figure}[t]
\centering
\includegraphics[width=0.41\textwidth]{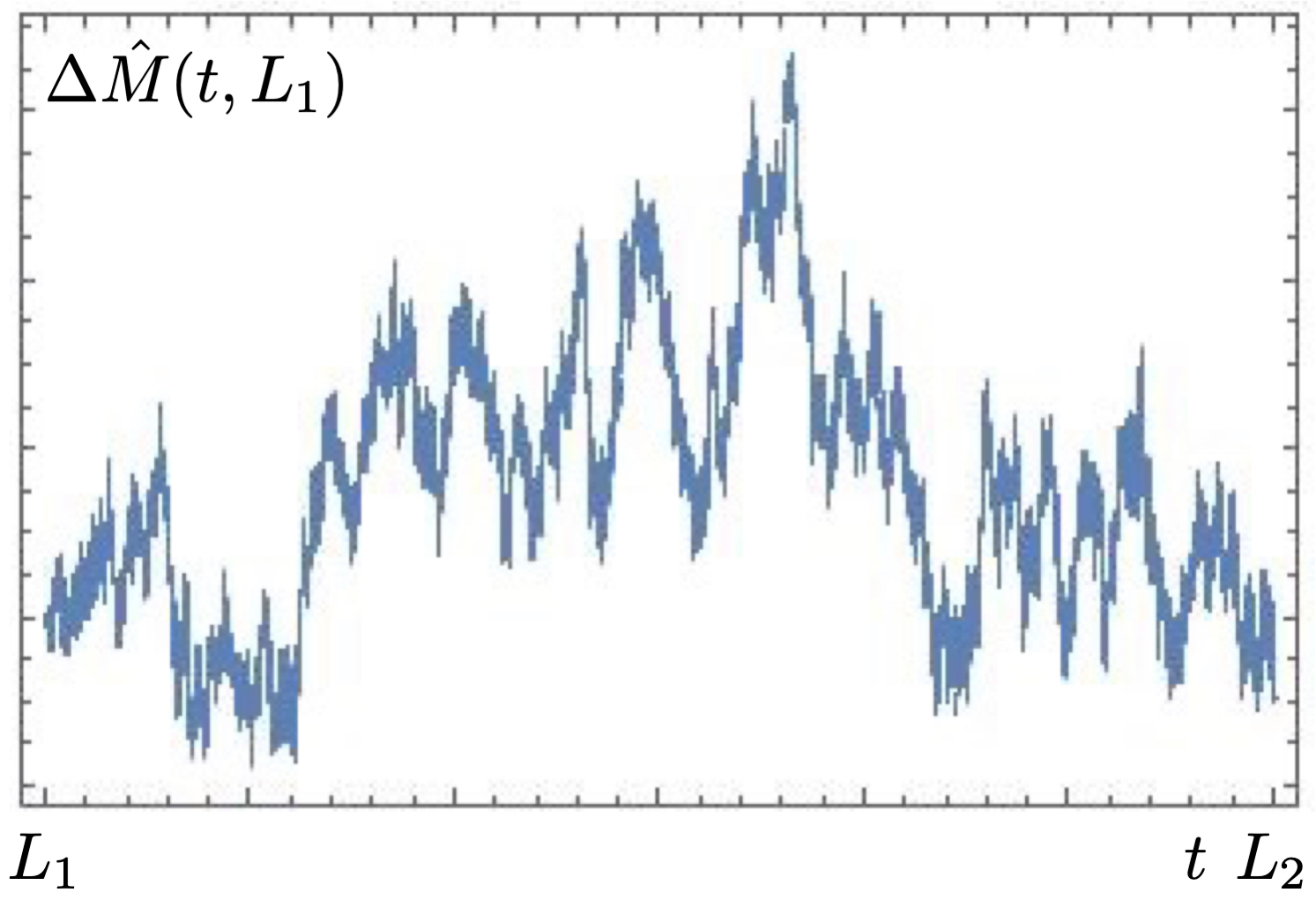}
\,\,\,
\includegraphics[width=0.40\textwidth]{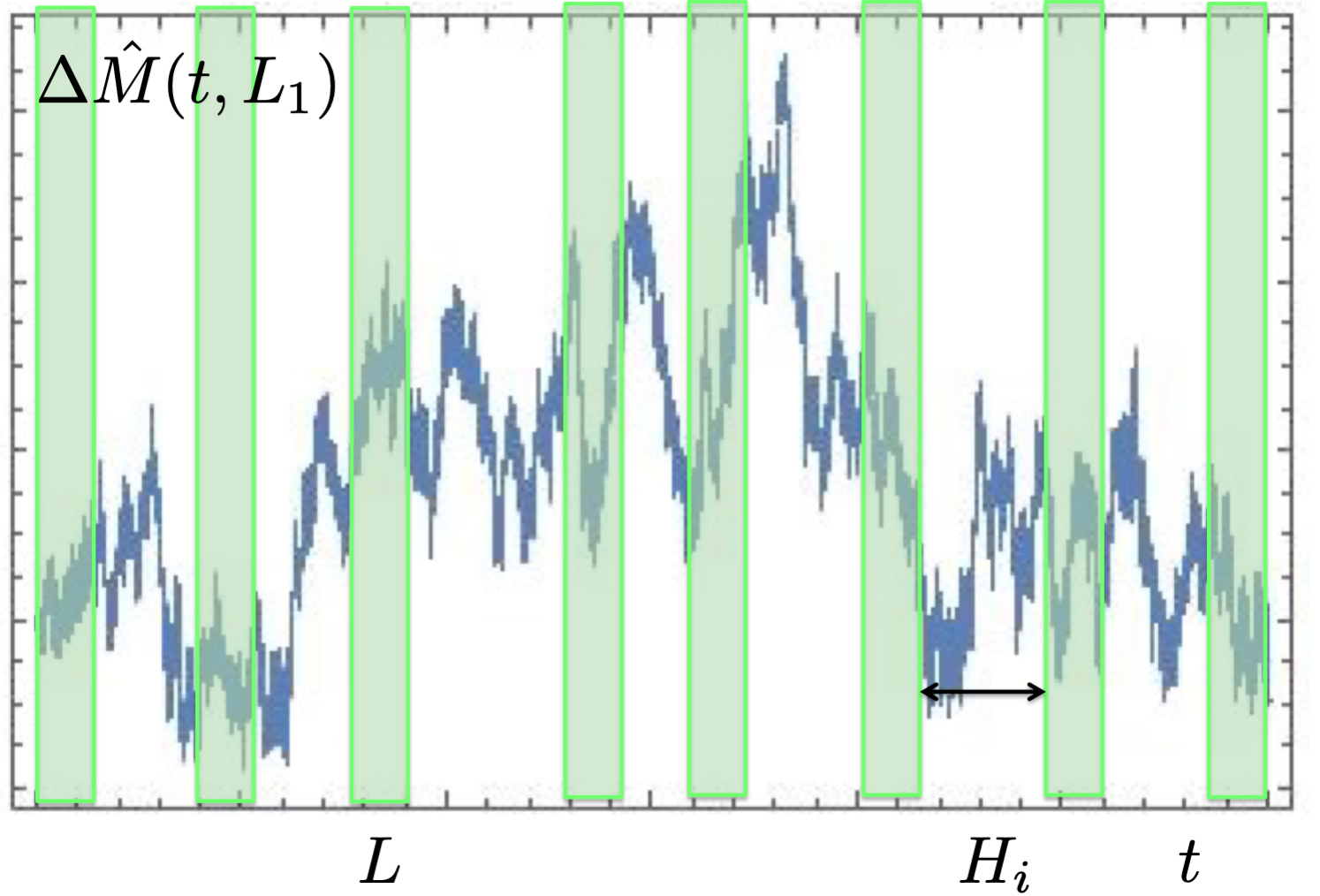}
\caption{Left hand side: series $\Delta\hat M(t,L_1) \equiv
\sum_{k=L_1}^{t} \hat\mu(k)$ vs $t$ in the interval $(L_1,L_2)$. 
Right hand side: sampling of the time series done in terms of block variables $B_L (\s)$ of length $L \gg 1$ relative to the green intervals separated by distances $H_i$. Under the hypothesis of stationarity of the sequence $A_{L_1,L_2}$, the values of these blocks define a probabilistic ensemble $\CE$ for the quantity $B_L$ relative to the sum of the first $L$ values $\hat\mu(n)$. } 
\label{partitiontimeseriesGN}
\end{figure}
\begin{enumerate}
\item Consider two very large integers $L_1$ and $L_2$  (which eventually we will send to infinity), with $L_1\gg 1$, $L_2\gg1$  but also $D \equiv (L_2 - L_1) \gg 1$.
\item For any fixed integer $N$, with $1 \ll N \ll D$, consider the union of  sets 
\beq 
{\cal G}_N \,=\, \bigcup_{i=1,\ldots N} I_L(i) 
\,\,\,\,\,\,\,\,
,
\,\,\,\,\,\,\,\,
L_1 \leq i < L_2 
\,\,\,\,\,\,,\,\,\,\,\, I_L(i) \cap I_L(j) = 0 
\,\,\,,\,\,\,  i \neq j 
\eeq
made of $N$ {\em non-overlapping} and also well separated intervals of length $L$ whose origin is between the two large numbers $L_1$ and $L_2$ (see Figure \ref{partitiontimeseriesGN}). The integer $N$ is the cardinality of the set ${\cal G}_N$.    These conditions ensure that the block variables $B_L(i)$ computed on such disjoint intervals are very weakly correlated and therefore we can assume that we are dealing statistically with $N$ separated copies of the original series $B_n(1)$. 
\item At any given $L_1$ and $L_2$, the cardinality ${\rm card} ({\cal G}_N) = N$ of these sets cannot be larger of course than $D/L$. There is however a large freedom in generating  them: 
\begin{enumerate}[label={\bf \alph*}]
\item We can take, for instance, $N$ intervals $I_L(\s)$ separated by a fixed distance $H$, with the condition that $N(L + H) = D$; 
\item Alternatively,  we can take,  $N$ intervals $I_L(\s)$ separated by random distances $H_i$ such that $N L + \sum_{i=1}^N H_i = D$. 
\end{enumerate}
\item The ensemble $\CE$ is then defined as  the set of the $N$ block variables $B_L(\s)$ relative to the intervals $I_L(\s) \in {\cal G}_N$:
\beq
\label{CEdef}
\CE = \{ B_L (\s) \},  ~~~~{\rm with~~}  I_L (\s) \in {\cal G}_N
\eeq
\end{enumerate}

\noindent 
In summary, choosing two very large and well separated  integers $L_1$ and $L_2$, we can generate a large number of sets of intervals ${\cal G}_N$ and use the corresponding block variables of length $L$ to sample the typical values taken by a series consisting of a sum of $L$ consecutive terms $\hat\mu(k)$. In view of the ergodicity, stationarity and scalability of the subsequences $\{\hat\epsilon_{L,l}\}$  
this is equivalent to determining  the statistical properties of the original series $B_n(1)$.


\section{Frequency Tests}
Starting with this section, we are going to present all statistical tests which we applied to check explicitly the randomness of our sequence $\{{\mathcal S}_n\}$. In this section we present the outputs of the simplest of these statistical checks which consist of frequency tests: we test, in order, the frequency of monobits, pairs, triplets, quadruplets and quintuplets. It is worth 
emphasizing that, in all cases, the null hypothesis of the random nature of the sequence $\{{\mathcal S}_n\}$ has been confirmed with high confidence levels: for instance, in the case of the frequency test of the quintuplets, the pattern of 5 consecutive numbers such as $(1, 1, 1, 0, 1)$ along the sequence $\{{\mathcal S}_n\}$ turns out to be as probable as $(1,0,0,1,1)$ or any of the other $32$ possibilities. It is worth saying that the successful outputs of our analysis emerged either applying the statistical tests to random subsets of the sequence $\{{\mathcal S}_n\}$, where we have made use of the intervals $I_L(l)$ previously discussed, or to progressively increasing lengths of the sequence $\{{\mathcal S}_n\}$ itself.   

\begin{table}[t]
\begin{center}
\begin{tabular}{|c| c | c | |c || c || c | c | |c | |c||}\hline
\multicolumn{5} {|c|} {${\mathcal S}_n$ sequence} & \multicolumn{4}{|c|}{ random sequence}\\\hline 
 block &\, \# 1  &  \# 0 & v &  P-value &  \, \# 1  &  \# 0 & v &  P-value \\\hline
1 &9999306 & 10000694 & 0.98 &  0.33 & 10001800 & 9998200 & 0.80 & 0.42\\
2 & 9998263 & 10001737 & 0.77 &  0.43 &10001502 & 9998498 & 0.67 &  0.50 \\
3 & 10000928 & 9999072 & 0.41 &  0.67 & 9998258 & 10001742 &  0.78 & 0.43\\
4 &10001456 & 9998544 & 0.65 &  0.51 & 10000964 & 9999036 & 0.43 & 0.67 \\
5 &9998082 & 10001918 & 0.85 &  0.39 & 9998122 & 10001878 &  0.83 & 0.40\\
\hline
\end{tabular}
\end{center}
\caption{
Frequencies of $1$'s and $0$'s for the sequence $\{{\mathcal S}_n\}$ (left hand side) and for a random sequence (right hand side) for five randomly chosen blocks of length $L = 2\times 10^8$, 
in the interval $(L_1, L_2) = (10^{14}, 10^{16})$. 
$v$ and $P$-value are respectively the value of the test statistic variable and the $P$-value of test.} 
\label{monobit1}
\end{table}
\begin{table}[b]
\begin{center}
\begin{tabular}{| c | c | c | |c | |c | | c | c | |c | |c||}\hline
\multicolumn{5} {|c|} {${\mathcal S}_n$ sequence} & \multicolumn{4}{|c|}{ random sequence}\\\hline 
 \,L&  \# 1  &  \# 0 & v &  P-value &  \, \# 1  &  \# 0 &v &  P-value \\\hline
$10^3$ & 492 & 508 & 0.50 &  0.61 & 494 &506 &0.37 & 0.70\\
$10^4$ & 4985 & 5015 & 0.30 &  0.76 &4924 & 5076 & 1.52 &  0.12 \\
$10^5$ & 50006 & 49994 & 0.03 &  0.96 & 50019 & 49981 &  0.12 & 0.90\\
$10^6$ &500093  & 499907 & 0.18&  0.85 & 499916 & 500084& 0.17 & 0.86 \\
$10^7$ & 5000402 & 4999598 & 0.25 &  0.79 & 5000491 & 4999509 &  0.31 & 0.75\\
$10^8$ & 50000577 & 49999423 & 0.11& 0.90 & 50005027 & 49994973 & 1.00 &0.31 \\
\hline
\end{tabular}
\end{center}
\caption{
Frequencies of $1$'s and $0$'s for the sequence $\{{\mathcal S}_n\}$ (left hand side), starting from $L_1=10^{14}$ and for a random sequence (right hand side) for blocks of increasing length $L$. 
  $v$ and $P$ are respectively the value of the test statistic variable and the $P$-value of test.  
} 
\label{monobit2}
\end{table}

\subsection{Monobit test} In this test we have checked the frequencies of the $0$'s and $1$'s along the sequence $\{{\mathcal S}_n\}$ and compared this output with the corresponding frequencies of true random sequences. We have performed the analysis, choosing the variables  $\hat\mu_k = 2 {\mathcal S}_k -1$. Defining 
\beq
V_n \,=\, \hat\mu_1 + \hat\mu_2 + \cdots \hat\mu_n\,\,\,,
\eeq
we take as the test statistic variable the absolute value of $V_n$ divided the square number of measurements\footnote{An alternative choice is to use the $\chi^2$ statistical test. We have also employed this test statistical variable, arriving to the same conclusions presented in the text.}
\beq
v\,=\, \frac{|V_n |}{\sqrt{n}} \,\,\,.
\eeq
In this test, the $P$-value is given by eq.\,(\ref{erfc}) 
\beq
{\rm P-value} \,=\,{\rm erfc}\left(\frac{v}{\sqrt{2}}\right)\,\,\,.
\eeq

We have performed two types of analysis: 
\begin{enumerate}
\item In the first type of analysis we have chosen 5 random non-overlapping intervals $I_l(l)$ with $L =2 \times 10^8$ and $l > 10^{14}$. 
The results of this analysis, together with the value of $v$  and the $P$-value, is reported in Table \ref{monobit1}. In the same table, for a direct comparison we have also reported the corresponding data of a random sequence. One can see that there were always very high $P$-values, i.e. our sequence $\{{\mathcal S}_n\}$ passes successfully this statistical test.  
\item We have made the same statistical analysis but including {\em progressively} more and more terms in the sequence. The results are shown in Table \ref{monobit2}. 
Also in this case there were always very high $P$-values, i.e. our sequence $\{{\mathcal S}_n\}$ passes successfully this statistical test.  
\end{enumerate}

\begin{table}[t]
\begin{center}
\begin{tabular}{||c| c | c | c | c | |c | |c|| } \hline
 \,block&  $f_1$  &  $f_2$ & $f_3$ & $f_4$ & $\chi^2$ &  P-value \\ \hline
1 & 0.12504 & 0.12501 & 0.12497 & 0.12497 & 0.51 & 0.91\\
2 & 0.12496 & 0.12496 & 0.12500 & 0.12508 & 1.55 & 0.67 \\
3 & 0.12503 & 0.12499 & 0.12500 & 0.12498 & 0.19 & 0.97 \\
4 &0.12505 & 0.12492 & 0.12506 & 0.12497 & 2.08 & 0.55 \\
5 & 0.12498 & 0.12489 & 0.12502 & 0.12511 & 3.97 & 0.26 \\
\hline
\end{tabular}
\end{center}
\caption{
Frequencies of the pairs $(0,0)$, $(0,1)$, $(1,0)$, $(1,1)$ for the sequence ${\mathcal S}_n$ for randomly chosen blocks of length $L = 2\times 10^7$ 
 in the interval $(L_1, L_2) = (10^{14}, 10^{16})$. 
$\chi^2$ and $P$ are respectively the value of the test statistic variable and the $P$-value of test.} 
\label{pairsfre1}
\end{table}
\begin{table}[b]
\begin{center}
\begin{tabular}{|| c | c | c | c | c | |c | |c|| } \hline
 \, L&  $f_1$  &  $f_2$ & $f_3$ & $f_4$ & $\chi^2$ &  P-value \\ \hline
 $ 10^3 $ &0.1264 & 0.12480 & 0.1254 & 0.12330 & 1.33 & 0.72\\
 $ 10^4 $  & 0.12496 & 0.12496 & 0.12500 & 0.12508 & 0.40 & 0.93\\
 $ 10^5 $  & 0.12553 & 0.12447 & 0.12546 & 0.12453 & 0.79 & 0.85 \\
 $ 10^6$   & 0.12480 & 0.12528 & 0.12545 & 0.12537 & 3.72 & 0.29 \\
 $10^7$    & 0.12512 & 0.12496 & 0.12491 & 0.12501 &  1.85 & 0.60 \\
$ 10^8$    & 0.12499 & 0.12502 & 0.12496 & 0.12502 &1.84 & 0.60\\ 
     \hline
\end{tabular}
\end{center}
\caption{
Frequencies of the pairs $(0,0)$, $(0,1)$, $(1,0)$, $(1,1)$ for the sequence $\{{\mathcal S}_n\}$ for blocks of increasing length $L$, starting from $L_1=10^{14}$, together with 
the values of $\chi^2$ and the $P$-values of the test. }
\label{pairsfre2}
\end{table}

\subsection{Pair frequencies}

In this test we have analysed the frequency distribution of the pair $({\mathcal S}_{2k}, {\mathcal S}_{2k+1})$ along our sequence $\{{\mathcal S}_n\}$. There could be 
$4$ possible outputs, given by 
\beq
a_1 = (0,0) \,\,\,,\,\,\, 
a_2 = (0,1)  \,\,\,,\,\,\, 
a_3 = (1,0)  \,\,\,,\,\,\, 
a_4 = (1,1) \,\,\,.
\label{pairoutput}
\eeq
If the sequence $\{{\mathcal S}_n\}$ is random, the frequencies $f_k$  of each of these outputs must be $1/4$. In this case, our test statistical variable is the chi-square 
\beq
\chi^2 \,=\, \sum_{k=1}^4 4 (f_k - 1/4)^2\,\,\,,
\label{chipairs}
\eeq
where $f_k$ is the measured frequency of each output given in (\ref{pairoutput}) and the $P$-value is given by the formula (\ref{qax}), i.e. 
\beq
P-{\rm value}\,=\,Q\left(\frac{3}{2},\frac{\chi^2}{2}\right) \,\,\,.
\eeq
As before, we show the results of our analysis for five random blocks of length $L = 2 \times 10^7$ (see Table \ref{pairsfre1}), in the interval $(L_1, L_2) = (10^{14}, 10^{16})$, 
and for blocks of increasing lengths, starting from $L_1 = 10^{14}$ (see Table \ref{pairsfre2}). Also in this case, the higher values of the $P$-values shows that our sequence 
passes successfully also this test.

\begin{table}[t]
\begin{center}
\begin{tabular}{||c| c | c | c | c | c | c | c | c | |c | |c|| } \hline
 \,  block & $f_1$  &  $f_2$ & $f_3$ & $f_4$ & $f_5$ &  $f_6$ & $f_7$ & $f_8$ & $\chi^2$ &  P-value\\ \hline
1 & 4.165 & 4.168  & 4.166 & 4.173 & 4.165 & 4.161 & 4.166 & 4.167 &  3.89 & 0.79 \\
2 & 4.165 & 4.168. & 4.166 &  4.163 & 4.166 & 4.169 & 4.167 & 4.167 &  1.47 & 0.98 \\
3 & 4.160 & 4.165 &  4.167 &4.174 &4.166 & 4.168 & 4.166 & 4.166 & 4.87 & 0.67 \\
4& 4.170 & 4.163 & 4.177 & 4.161 & 4.171 & 4.160 & 4.161 & 4.165 & 13.7 & 0.05 \\
5& 4.162 & 4.170 & 4.170 & 4.168 & 4.162 & 4.166 & 4.169 & 4.163 & 3.68 & 0.81\\\hline
 \end{tabular}
\end{center}
\caption{
Frequencies (in units of $10^{-2}$) of the triplets $a_1 = (0,0,0)$,  $a_2 = (0,0,1)$, 
$a_3 = (0,1,0)$,   
$a_4 = (1,0,0) $,
$a_5 = (1,1,1)  $,
$a_6 = (1,1,0)  $,
$a_7 = (1,0,1)  $,
$a_8 = (0,1,1) $ for the sequence $\{{\mathcal S}_n\}$ for randomly chosen blocks of length $L = 2\times 10^7$, 
 in the interval $(L_1, L_2) = (10^{14}, 10^{16})$.
$\chi^2$ and $P$-value are respectively the value of the test statistic variable and the $P$-value of test. }
\label{tripletsfre1}
\end{table}

\subsection{Triplet frequencies} 
In this test we have analysed the frequency distribution of the triplets $({\mathcal S}_{3k}, {\mathcal S}_{3k+1}, S_{3k +2})$ along our sequence $\{{\mathcal S}_n\}$. There could be 
$8=2^3$ possible outputs, given by 
\begin{eqnarray}
&& a_1 = (0,0,0) \,\,\,,\,\,\, 
a_2 = (0,0,1)  \,\,\,,\,\,\, 
a_3 = (0,1,0)  \,\,\,,\,\,\, 
a_4 = (1,0,0) \,\,\, \label{tripletutput}
\\
&& 
a_5 = (1,1,1)  \,\,\,,\,\,\, 
a_6 = (1,1,0)  \,\,\,,\,\,\, 
a_7 = (1,0,1)  \,\,\,,\,\,\, 
a_8 = (0,1,1) \nonumber 
\end{eqnarray}
If the sequence ${\mathcal S}_n$ is random, the frequencies of each of these outputs must be $1/8$. In this case, our test statistical variable is the chi-square 
\beq
\chi^2 \,=\, \sum_{k=1}^8 8 (f_k - 1/8)^2\,\,\,,
\label{chipairsa}
\eeq
where $f_k$ is the measured frequency of each output given in (\ref{tripletutput}) and the $P$-value is given by the formula (\ref{qax}), i.e. 
\beq
P-{\rm value} \,=\,Q\left(\frac{7}{2},\frac{\chi^2}{2}\right) \,\,\,.
\eeq
As before, we show the results of our analysis for five random blocks of length $L = 2 \times 10^7$ (see Table \ref{tripletsfre1}) 
and, also for the triplets, the values of the $P$-values show that our sequence passes successfully this test.

\begin{figure}[b]
\centering\includegraphics[width=1.0\textwidth]{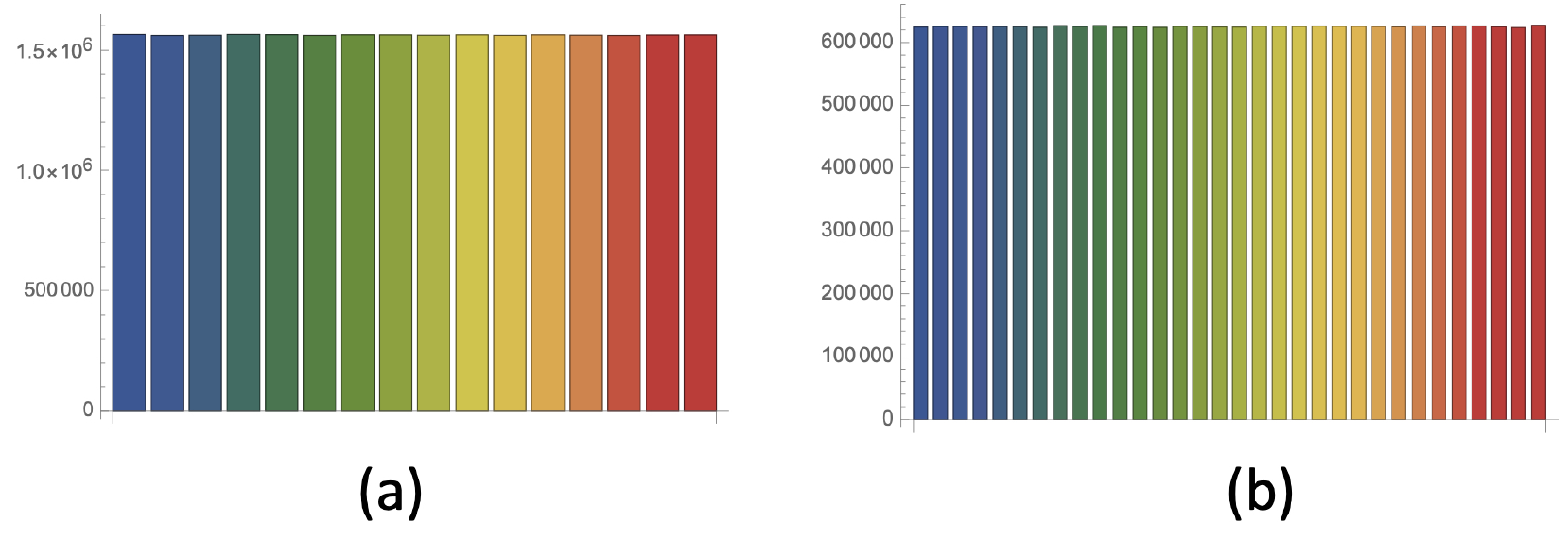}
\caption{Histograms of the outcome frequencies of (a) quadruplets and (b) quintuplets for a block variable of length $N = 10^8$,
 in the interval $(L_1, L_2) = (10^{14}, 10^{16})$. }  
\label{chartquadr}
\end{figure}

\subsection{Quadruplet and quintuplet frequencies} 

For the statistical analysis of quadruplets and quintuplets it becomes rather cumbersome to show all the frequencies of the relative outcomes (there are $16$ outputs for the quadruplets and $32$ for the quintuplets). Hence, we simply decide to show in Figure \ref{chartquadr} the histograms of the quadruplet and quintuplet frequencies for a block variable of length $L = 10^8$  extracted around $L_1=10^{14}$ and, from these histograms, one can see a perfect balance among all frequencies of the events. We also report in Table \ref{poker1}   the $\chi^2$ values (given by eq.\,(\ref{chisquare}), with $r=16$ and $p_s = 1/16$) and the $P$-values for the frequency distribution of the quadruplets for blocks of increasing length $n$, with $L_1=10^{14}$. The same for the quintuplets, whose chi-square is computed according to eq.\,(\ref{chisquare}) with $r=32$ and $p_s = 1/32$, and with all relative data reported in Table \ref{poker1}. Also in these case, the higher values of the $P$-values do not indicate a deviation of randomness of the sequence $\{{\mathcal S}_n\}$.

\section{Proportion of sequences passing a test}\label{passingtest}

We can have an additional control on the significance of a statistical test we have done by computing the proportion of the results that pass the test\footnote{We strongly advice the reader to read Section 4.2 of the NIST report \cite{Nist} for further details. on this topic.}.
We illustrate here the result for one particular test but we have checked that the same conclusions also apply to all other tests discussed in this part of the paper. Consider, for instance,  the test about the equidistribution of triplets along the sequence 
$\{{\mathcal S}_n\}$.  The protocol of the procedure is as follows

\begin{table}[t]
\begin{center}
\begin{tabular}{lll}
\begin{tabular}{||c |c  |c|| } \hline
 \,L &    $\chi^2$ &  P-value \\ \hline
$10^3$   & 12.6 & 0.63 \\
$10^4$   & 7.13. & 0.95\\
$10^5$   &6.13 &  0.97\\
$10^6$   & 20.54& 0.15 \\
$10^7$   & 8.89 & 0.88\\
$ 10^8$  & 15.6& 0.40 \\
\hline
 \end{tabular}
& \hspace{3cm} & 
\begin{tabular}{||c |c  |c|| } \hline
 \,L &    $\chi^2$ &  P-value \\ \hline
$10^3$   & 37.2 & 0.20 \\
$10^4$   & 38.1. & 0.16\\
$10^5$   &28.4 &  0.60\\
$10^6$   & 38.8& 0.15 \\
$10^7$   & 21.6 & 0.89\\
$ 10^8$  & 24.1& 0.80\\
\hline
 \end{tabular}
\end{tabular}
\end{center}
\caption{Table for the quadruplets (left hand side) and for the quintuplets (right hand side). 
 Block variables of increasing length $L$ extracted from the sequence $\{{\mathcal S}_n\}$ starting from $L_1 = 10^{14}$. 
$\chi^2$ and $P$ are respectively the value of the test statistic variable and the $P$-value of test.  }
\label{poker1}
\end{table}

\begin{enumerate}
\item For our sequence $\{{\mathcal S}_n\}$ we consider $N$ non-overlapping intervals $I_L(l)$  of length $L$ with increasing starting points $l$ (see eq.(\ref{intervals})). 
Associated to each of these intervals, there are $N$ subsequences $\{\epsilon^{(1)}\}, \{\epsilon^{(2)}\},  
\ldots, \{\epsilon^{(N)}\}$ of $\{\mathcal{S}_n\}$, each of length $L$. We compute the frequencies
of the various triplets in each of these subsequences and we compute the relative $P$-values. Considered altogether, these steps form what we call a {\em sweep}. 
\item We compute the proportion of subsequences that pass the test. In our case, having fixed $\alpha = 0.01$, we count the number of subsequences for which the $P$-values $ \geq 0.01$.  
\item We repeat the sweep of point 1 for another $T$ times, choosing, each time, $N$ intervals of length $L$ never overlapping with any previous one. 
\item There is now a statistical formula \cite{Nist} which determines the range of acceptable proportions in terms of the confidence interval: this is given by 
\beq
{\mathcal I}_{\pm} \,=\, 
(1 - \alpha) \pm 3 \sqrt{\frac{\alpha (1-\alpha)}{N}}
\label{proportionsequnce}
\eeq
where $N$ is the sample size, sufficiently high in order to have a reasonable statistics (say $N > 10$). Note that as the sample size $N\to \infty$, 
$\, \,  {\mathcal I}_\pm = 1-\alpha$ as expected. 
\item If the proportion falls outside this interval, then there is evidence that the data is non-random.  
\end{enumerate}

As our example we have chosen $N=100$, $L = 10^5$ and $T =14$. So, in our case, the confidence interval is given by  ${\mathcal I}_+ = 1.02$ and ${\mathcal I}_- = 0.96015$. In Figure \ref{proportions}.a we show the result of this type of check for the triplets and all proportions are perfectly within the confidence interval. Similar results also hold for the tests concerning monobits, pairs, quadruplets and quintuplets.

\subsection{ Uniform distribution of P-values} 
\begin{figure}[b]
\centering\includegraphics[width=1.0\textwidth]{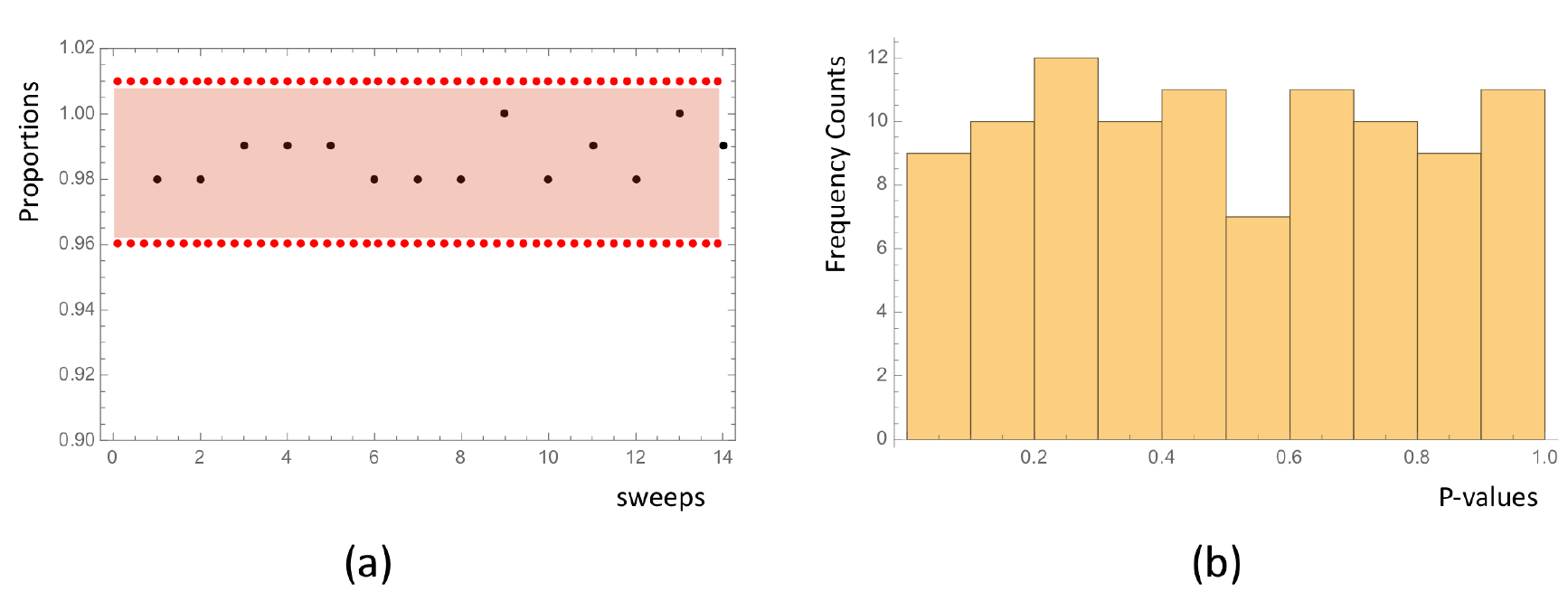}
\caption{(a) Proportions of sequence that pass the statistical test $P$-value $\geq 0.01$. The red region is between the lower and higher values of the confidence interval identified 
by ($\mathcal{I}_-,\mathcal{I}_+$), see eq.\,(\ref{proportionsequnce}).  (b) Histogram of $P$-values for the triplet frequency test for a particular sweep.}
\label{proportions}
\end{figure}

One can also check the uniform distribution of the $P$-values. This is easily done using a histogram of the various values obtained for the 
$P$-values with the interval $(0,1)$ divided into 10 sub-intervals, with the relative counting of how many $P$-values lie within each sub-interval. The quantitative degree of uniformity can be controlled in terms of a $\chi^2$ test, in this case given by 
\beq
\chi^2 \,=\, \sum_{i=1}^{10} \frac{\left(F_i - S_S/10\right)^2}{S_S/10} \,\,\,,
\label{aboveqq}
\eeq
where $S_S$ is the sample size and $F_i$ are the frequencies of the $P$-values relative to each interval. 
We can now calculate a $P$-value for this distribution of the P-values! Let's denote it by $\overline P$: it is associated to the $\chi^2$ of (\ref{aboveqq}) and it  
is given by
\beq
\overline P-{\rm value} \,=\,Q\left(\frac{9}{2},\frac{\chi^2}{2}\right) \,\,\,. 
\eeq
According to \cite{Nist}, if $\overline P \geq 0.0001 $, then the sequences can be considered to be uniformly distributed.  
 In Figure \ref{proportions}.b we report a typical histogram of the $P$-values for 
one of our sweeps of the triplet frequencies, with the corresponding $P$ value equal to $\overline P = 0.0351$ and therefore satisfies the bound for a uniform distribution of the $P$-values.  


\section{Correlation function and discrete Fourier transform}

Let's consider a generic time series whose generic output is here denoted $a_t$. Useful information on the nature of this time series can be obtained by computing the correlation between the two outputs $a_t$ and $a_{t+s}$, 
where the interval $s$ is called the lag.  For a sequence of length $T$, the correlation function of lag $s$ is defined as 
\beq
\gamma(s) \,=\, \left(\sum_{i=1}^{T-s} (a_i - \overline m)(a_{i+s} - \overline m)\right)/\sum_{i=1}^T (a_i - \overline m)^2\,\,\,,
\eeq
where $\overline m $ is the average of the sequence $a_n$.  If we now take as time series our sequence $\{\mathcal {S}_n\}$, we can refine our definition of the correlation function and take as definition 
the quantity
\beq
F(s) \,=\, \langle \gamma(s) \rangle \,\,\,,
\eeq
where the average is taken wrt our stroboscopic sampling of our sequence $\{\mathcal {S}_n\}$. This means that we consider $N$ non overlapping intervals $I_L(l)$ randomly chosen, with $L \gg T$, 
and, in each of these intervals,  we consider a sub-interval of length $T$ randomly chosen to extract a subsequence $\{\epsilon\}$ and to compute $\gamma(s)$. Finally we average on all the $N$ $\gamma(s)$'s computed in this way. In our analysis we have typically chosen $N=100$ and $T  =  25000$. 

\begin{figure}[t]
\centering\includegraphics[width=1.0\textwidth]{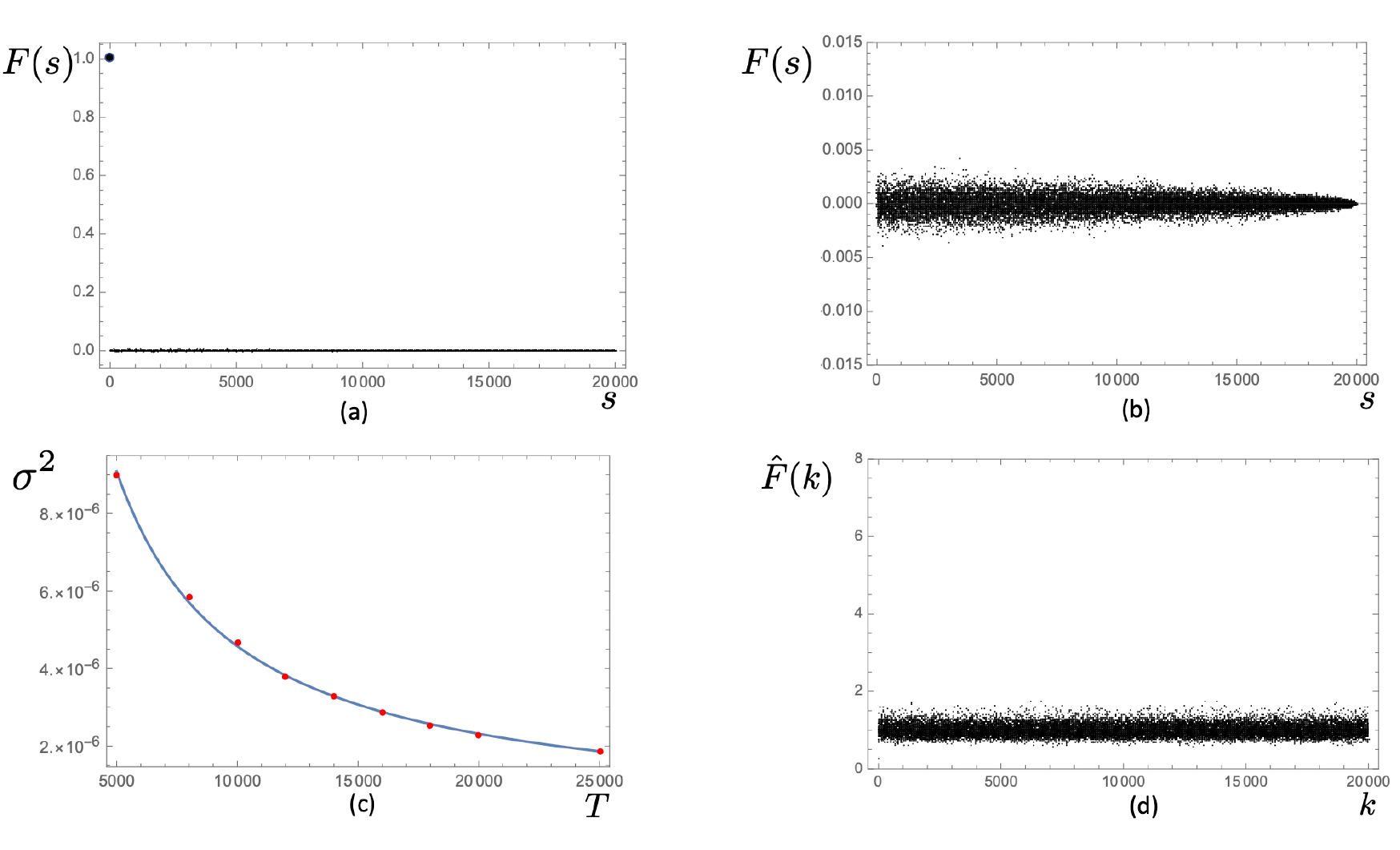}
\caption{(a) Correlation function $F(s)$ up to lag $s=20000$ for the sequence $\{\mathcal {S}_n\}$. Except at $s=0$, the correlation function is almost zero. (b) Zoom of the correlation function, which shows that its values are always of the order or smaller than $1/\sqrt{T}$, where here $T=20.000$. (c) Variance $\sigma^2$ of the values of the correlation function versus the length $L$ of the sequence, together with the curve of the best fit $\sim 1/T$. (d) Squared magnitudes of the Fourier transform $\hat F(k)$ of the correlation function, almost constant in the entire interval of the frequencies.}
\label{correlationfunctionfigure}
\end{figure}

For an infinitely long ideal true random sequence (i.e. $T \rightarrow \infty$) it holds  
\beq
F(s) \,=\,\left\{ 
\begin{array}{lll}
1 & , & {\rm if}\, s=0 \\
0 & , & {\rm if } \, s\neq 0 
\end{array}
\right.
\eeq 
However, our definition involves a {\em finite} sequence of length $T$: in this case, for uncorrelated quantities of length $T$, due to finite-size effects we expect that the variance $\sigma^2$ of the sequence of the $\gamma(s)$ is instead of the order $1/T$. This is indeed what is shown in Figure \ref{correlationfunctionfigure}.a and \ref{correlationfunctionfigure}.b, while the finite size scaling 
of the variance $\sigma^2$ with respect to $T$ is show in Figure  \ref{correlationfunctionfigure}.c.  

The Discrete Fourier transform of the correlation function is defined as 
\beq
\hat F(k) \,=\,\sum_{k=0}^{T-1} F(s) e^{i s k/L}\,\,\,,
\eeq
and its eventual peaks provide an indication of the periodicities present in the sequence. If $\{S_n \}$ was a true random sequence, ${\hat F}(k)$ would be a constant function of value $1$. As shown in Figure \ref{correlationfunctionfigure}.d, for our sequence $\{\mathcal {S}_n\}$ the squared magnitudes of the discrete Fourier transform is a function almost constant around the expected value $1$, with a spread due to the finite length $T$ of the sequences analyzed.  

\subsection{ Discrete Fourier Transform Test} 
\begin{figure}[t]
\centering\includegraphics[width=0.7\textwidth]{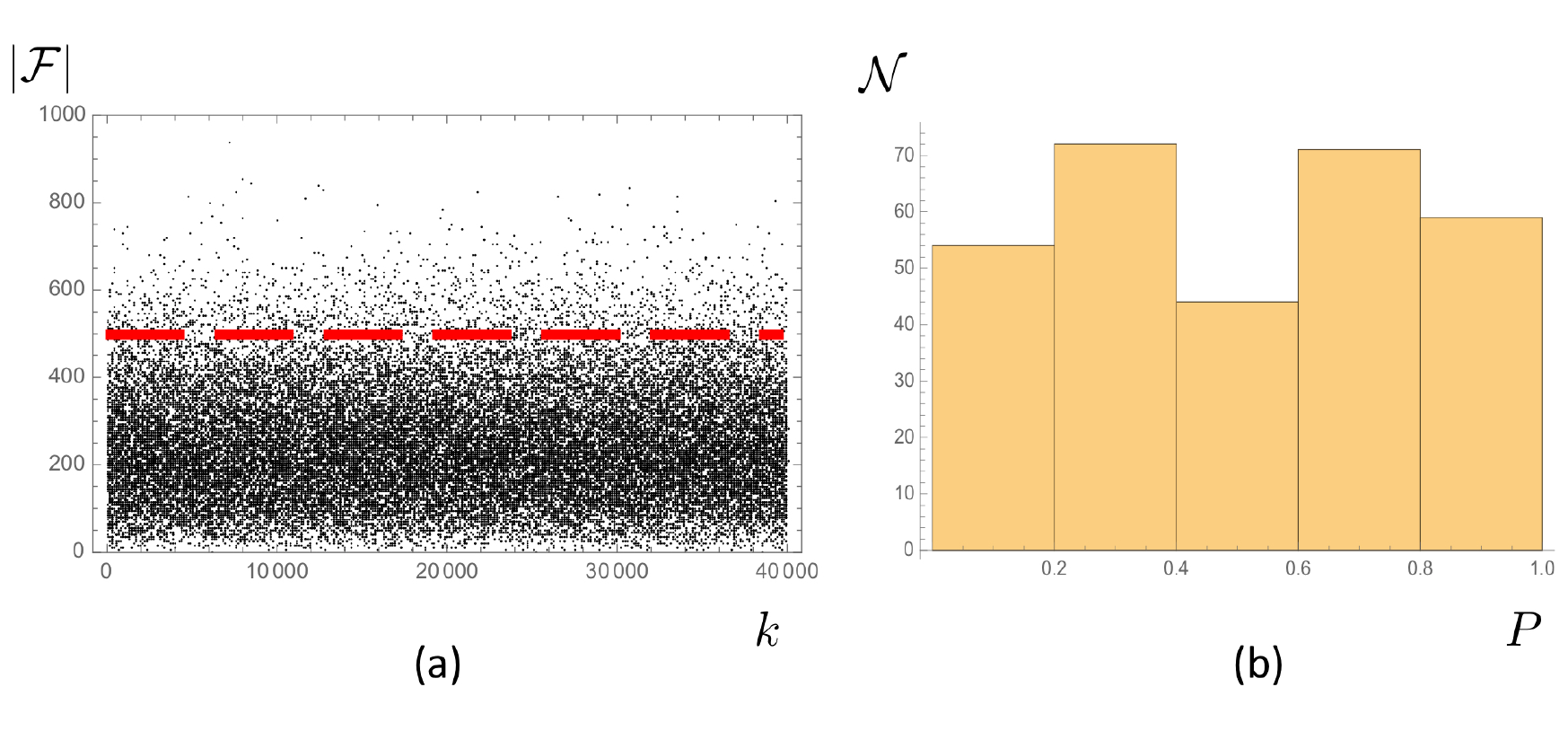}
\caption{(a) Squared magnitudes of the Discrete Fourier transform of $L=8 \times 10^4$ consecutive values of $\{\hat\mu_n\}$. The red dashed line is the peak height threshold value. (b) Distribution of the $P$-values for this test. All the obtained $P$-values pass the statistical test $P$-value $ > 0.01$. }
\label{DiscreteFourierAnalysis}
\end{figure}

The purpose of this test is to study the peak heights in the Discrete Fourier Transform of our sequence in order to detect, if any, its periodic features. The presence of these periodic features would indicate in fact a deviation from the randomness. We focus our attention to detect whether the number of peaks exceeding  $95\%$ of the threshold value ${\mathcal T}$ (given below) is significantly different than $5\%$. The way the test is implemented is as follows
\begin{enumerate}
\item Within a large interval $J_{L_1,L_2}$ around $L_1 \sim 10^{12}$ and $L_2 \sim 10^{16}$, we select randomly 
$N$ intervals $I_L(l)$ associated to subsequences $\{\hat\epsilon\}$ of length $ L$ of our sequence $\{\hat\mu_n\}$. In more detail, we have taken 
$N = 300$, $L=8 \times 10^4$, with $l \geq 2\times 10^{12}$. 
\item We compute the Discrete Fourier Transform $\mathcal {F}$ of the subsequence made of $ L$ consecutive restricted M\"{o}bius coefficients.
\item We calculate the modulus $|{\mathcal F}|$ of the first $L/2$ elements of the Discrete Fourier Transform.  
\item We compute the threshold value 
\beq
{\mathcal T}\, =\, \sqrt{\left(\log\frac{1}{0.05}\right) \,  L} \,\,\,,
\eeq
(for instance, in our case ${\mathcal T}= \, 489.549..$). 
This is the $95\%$ peak height  threshold value. If the sequence is random, $95\%$ of the values of $|{\mathcal F}|$ should not exceed ${\mathcal T}$. 
\item We define the quantity $N_0 \equiv 0.95 \,  L/2$, which is the expected ($95\%$) theoretical number of peaks that have to be less than the threshold ${\mathcal T}$ if the sequence is random. 
\item We count $N_e$, i.e. the actual observed number of peaks in $|{\mathcal F}|$ that are less than ${\mathcal T}$ and compute the statistical test 
\beq
v \,=\, \frac{(N_e - N_0)}{\sqrt{ L (0.95) (0.05) /4}}\,\,\,.
\eeq
\item Then we compute the $P$-value of this test, given by 
\beq
P-{\rm value} \,=\, {\rm erfc}\left(\frac{|v|}{\sqrt{2}}\right)
\,\,\,.
\eeq
\item If the computed $P$-value is $ < 0.01$, then the sequence is not random. 
\end{enumerate}
A typical output of this test is shown in Figure \ref{DiscreteFourierAnalysis}.a, where we also show (the red dashed line) the threshold value ${\mathcal T}$.  We perform our tests on sequences made of 
$ L= 8 \times 10^4$ consecutive numbers chosen randomly among some stroboscopic large intervals of the overall sequence $\{\hat\mu_n\}$ (with initial point of our analysis $L_1 =10^{14}$, 
see Section \ref{intervalsss}), always obtaining $P$-values larger than the decisional value $0.01$. Hence the sequence $\{{\mathcal S}_n\}$ passes successfully also this test. A histogram of the $P$-values we obtain is shown in Figure \ref{DiscreteFourierAnalysis}.b.  

\subsection{Correlation with random vectors} 
\begin{figure}[b]
\centering\includegraphics[width=0.6\textwidth]{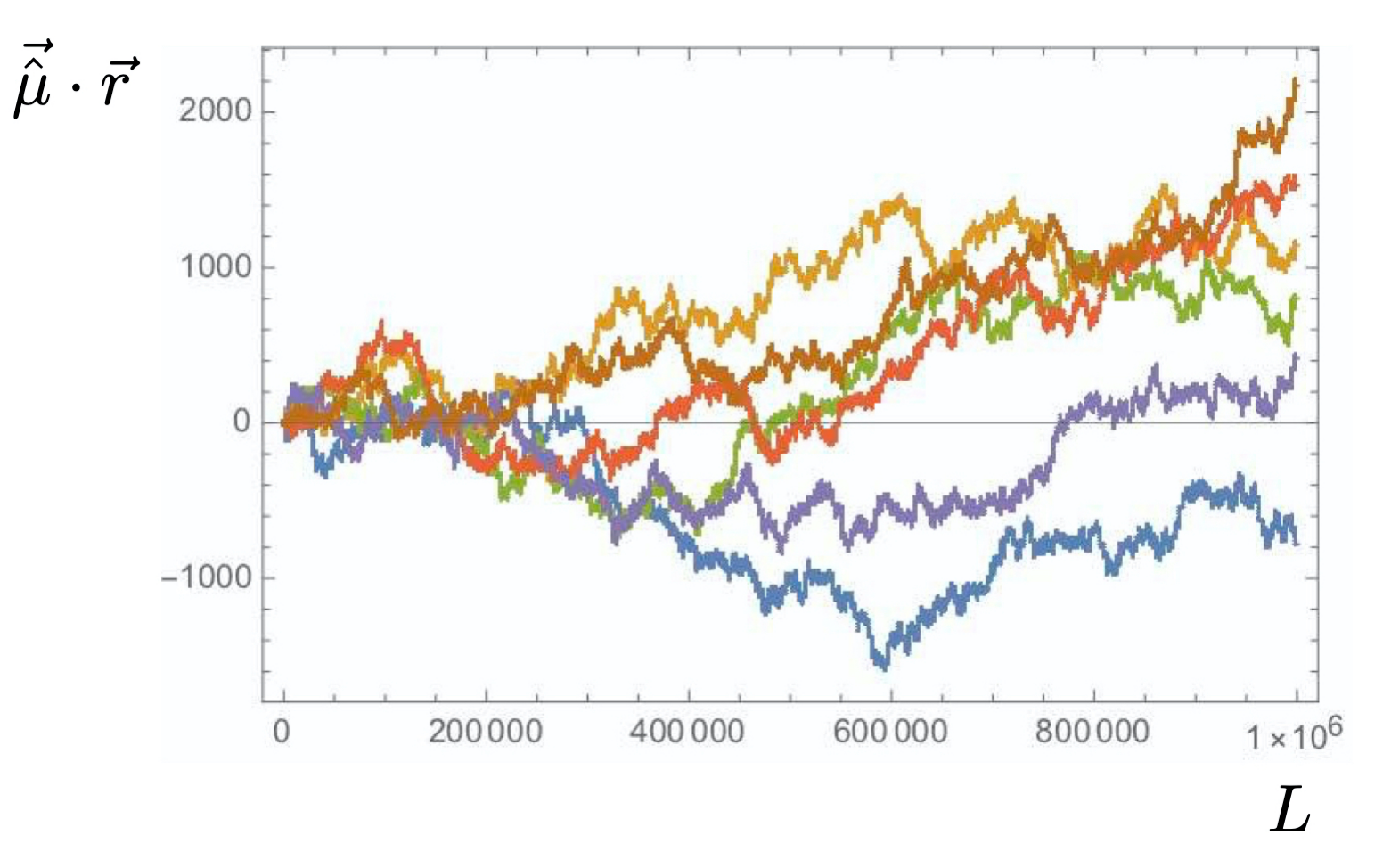}
\caption{Correlation between a sequence of $L$ consecutive values of the restricted M\"{o}bius coefficients along the square-free numbers and six random sequences of numbers $\{\pm 1\}$. }
\label{randomcorrelation}
\end{figure}

As stated in \cite{Iwaniec,Green-Tao}, the randomness of the M\"{o}bius coefficients would imply that they do not correlate with any reasonable sequence $r(n)$. In our case we can reformulate this statement as follows: (a) first of all we focus our attention to the restricted values of the M\"{o}bius function to square-free numbers only, i.e. the sequence $\{\hat \mu(n)\}$; (b) secondly, we extract a random consecutive sequence of $L$ values from this sequence; (c) thirdly, we compute the correlation of this sequence with another sequence $r(n)$ made of $L$ values. This ends up in the formula  
\beq
\hat\mu   \cdot r \,=\,\sum_{k=1}^L \hat\mu(k) \,r(k) \,\,\,
\sim o(L)
\label{littleo}
\eeq
where $o(L)$ means that, for $L \rightarrow \infty$, $\mu \cdot r /L \rightarrow 0$. 
In Figure \ref{randomcorrelation} we present the correlation of $\{\hat \mu_n\}$ versus $L$ with six random sequences $r(n)$ made of $\{\pm 1\}$: all these correlations seem not to grow faster than $\sqrt{L}$, so that the asymptotic behaviour dictated by (\ref{littleo}) seems to be indeed ``experimentally" well verified.

\section{Oscillation and long run of ones tests}

In this chapter we address two more refined tests concerning the randomness of our sequence $\{{\mathcal S}_n\}$: the first concerns the variation of the $1$'s and $0$'s present in it, the second the longest strings of these symbols. 

\subsection{ Oscillation test} Let's consider the oscillations of our sequence $\{{\mathcal S}_n\}$ and their statistics. An oscillation is considered to be a change from a 1 to 0 and vice versa. So, for instance, in this short piece of a sequence 
$$
0 \, 0 \, 0 \, 1 \, 1  \,  0  \, 1  \,  0  \, 1  \,  0 \,  0 \,  0  \, 1 \, 1 \, 0\, 
$$
the oscillations are in coincidence with the vertical bar hereafter 
$$
0 \, 0 \, 0 \, || \, 1 \, 1  \, ||\, 0  \, ||\, 1  \, || \, 0  \, || \, 1  \, ||\,  0 \,  0 \,  0  \, ||\, 1 \, 1 \, ||\, 0\, 
$$
We are in presence of a fast oscillation when there are a lot of changes. To have a quick birds-eye view of the configurations, it is convenient to adopt a binary graphical notation, in which a black box is associated to $1$ while a white box is associated to $0$. For instance, the alternating sequence 
$$ 
0 \, 1\, 0 \, 1 \, 0  \,  1  \, 0  \,  1  \, 0  \,  1 \,  0 \,  1  \, 0 \, 1 \, 0\, 
$$ 
(which has the highest possible oscillation), has the graphical representation 
\vspace{1mm}

\begin{center}

\begin{Chessboard}
\Row{B,A,B,A,B,A,B,A,B,A,B,A,B,A,B}
\end{Chessboard}
\end{center}
\vspace{1mm}

\noindent
while a sequence as 
$$
0 \, 0 \, 0 \, 0 \, 0  \,  0  \, 0  \,  1  \, 1  \,  1 \,  1 \,  1  \, 1 \, 1 \, 1\, 
$$ 
which has just one oscillation, graphically appears as 
\begin{center}

\begin{Chessboard}
\Row{B,B,B,B,B,B,B,A,A,A,A,A,A,A,A}
\end{Chessboard}
\end{center}

\vspace{1mm}
\noindent
In statistical mechanics, if $0$ and $1$ denote the two ground states of an Ising like system, each oscillation corresponds to a kink state which interpolates between the ground states. 

\subsection{ Bernoulli variables} 
In order to set up a statistical test for the oscillations, let's fix some notation and some basic relations: hereafter $\tilde{\mathcal S}_k$ (with $k=1,2,\ldots,L$) denotes the $k$-th element 
of a sub-sequence $\{\epsilon\}$ of 
$\{\mathcal {S}_n\}$, while  $n_0$ and $n_1$ respectively denote the number of $0$'s and $1$'s in this sub-sequence. Therefore we have
\beq
n_0 + n_1 \,=\,  L \,\,\,.
\eeq
In the following we denote by $\rho$ the fraction of $1$'s present in the sequence 
\beq
\rho \,=\,\frac{1}L \, 
\sum_{j=1}^L \tilde{\mathcal {S}}_j\,\,\,,
\eeq
Let's also define the variable $r(k)$ as 
\beq
r(k) \,=\,\left\{\begin{array}{lll}
0 & & \,\,{\rm if} \, \,\,\tilde{{\mathcal S}}_k = \tilde{{\mathcal S}}_{k-1} \\
1 & & \,\,{\rm otherwise}
\end{array}
\right.
\eeq
Then we can set up our statistical test $V$ in terms of the oscillations present in the sequence, defined as 
\beq
V \,=\, 1+ \sum_{k=2}^{L} r(k) \,\,\,.
\label{numoscillation}
\eeq
$V$ simply measures the number of kinks (oscillations) present in the system: a large value of $V$ indicates that there are fast oscillations in the sequence while a small value that there are oscillations  that are too slow.  Both cases are not typical for a  true random sequence. 

The quantity $r_k = r(k)$ is a Bernoulli variable with parameter $p = 2 n_0 n_1/(L (L-1))\sim 2 \rho (1-\rho)$, since  we have 
\begin{equation}
\begin{array}{cllll}
\tilde S_k \, \tilde S_{k+1} & & r_k  & & {\rm probability}\\
1 \, 1 & & 0 & & \rho^2 \\
0 \, 1 & & 1 & & \rho (1-\rho) \\
1 \, 0 & & 1 & & \rho(1-\rho) \\
0 \, 0 & & 0 & & (1-\rho)^2 
\end{array}
\end{equation}
Hence, $r_k$ takes value $1$ only in the two cases in the middle and their probability is then $2 \rho(1-\rho)$.  Therefore 
\beq
\langle r_k \rangle \,=\,\langle r_k^2 \rangle \,=\,2 \rho (1-\rho) 
\,\,\,.
\eeq
Since $V$ is a linear combination of the variables $r_k$ we have 
\beq
\overline V \,=\,\langle V \rangle \,=\, 1 + \sum_{k=2}^{L} \langle r_k \rangle \,=\, 1 + 2 (L-1) \rho (1-\rho)\sim 2\,L \,\rho \,(1- \rho) \,\,\,.
\eeq
For its variance we have 
\beq
\langle (V - \overline V)^2\rangle \,= \, \langle V^2 \rangle - \overline V^2 \,\,\,,
\eeq
where 
\begin{eqnarray}
\langle V^2 \rangle & = & \left\langle \left(1 + \sum_{k=2}^L r_k\right) \left(1 + \sum_{j=2}^L r_j \right) \right\rangle\,=\, 
1 + 2 \sum_{k=2}^L \langle r_k \rangle + \sum_{j,k=2}^L \langle r_k r_j \rangle \\
&=& 1 + 4 (L-1) \rho (1- \rho) + (L-1)^2 \langle r^2_k \rangle + \sum_{2 \leq j \neq k \leq L} \langle r_j r_k \rangle \nonumber\,\,\,.
\end{eqnarray}
For the $2 (L-2)$ cases in which $j = k \pm 1$ we have 
\beq
\langle r_j r_k \rangle \,=\,  \frac{n_0 n_1 (n_1 -1) + n_0 n_1 (n_0-1)}{L (L-1) (L-2)} \,=\,\frac{n_0 n_1}{L (L-1)} \,\sim \rho (1-\rho) \,\,\,,
\eeq
while for the remaining $(L-1) (L-2) - 2 (L-2) = (L-2) (L-3)$ cases where $j\neq k$ we have 
\beq
\langle r_j r_k 
\rangle \,=\, \frac{4 n_0 n_1 (n_0-1)(n_1-1)}{L (L-1) (L-2) (L-3)} \,\,\,,
\eeq
\begin{figure}[t]
\centering\includegraphics[width=1.0\textwidth]{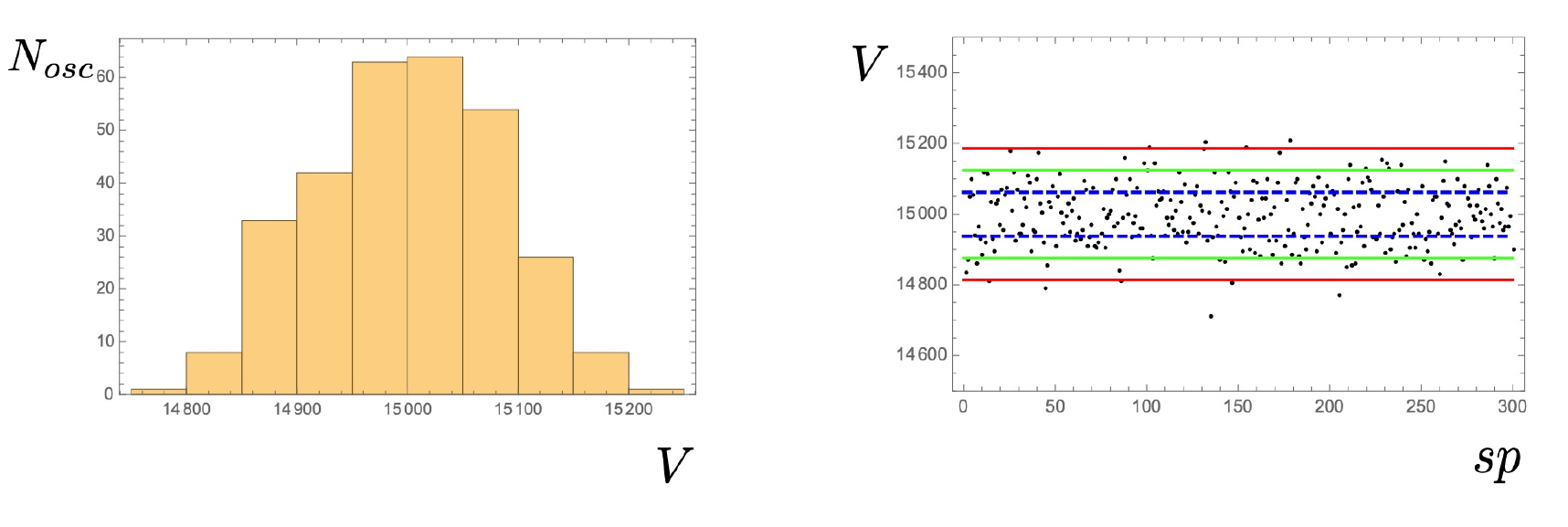}
\caption{(a) Histogram of the number of oscillations $N_{osc}$ for block sequences of length $L=30000$ of ${\mathcal S}_n$ taken randomly from $\tilde N=300$ different non-overlapping stroboscopic intervals $I_{L}(l)$ of length $L = 50000$, where $l \in ( 10^{14},  10^{15})$. (b) Values of the oscillations for each of the $300$ sweep, together with the intervals relative to the $1$, $2$ and $3$ standard deviations (blue, green and red). }
\label{oscillationnn1}
\end{figure}
\begin{figure}[b]
\centering\includegraphics[width=1.0\textwidth]{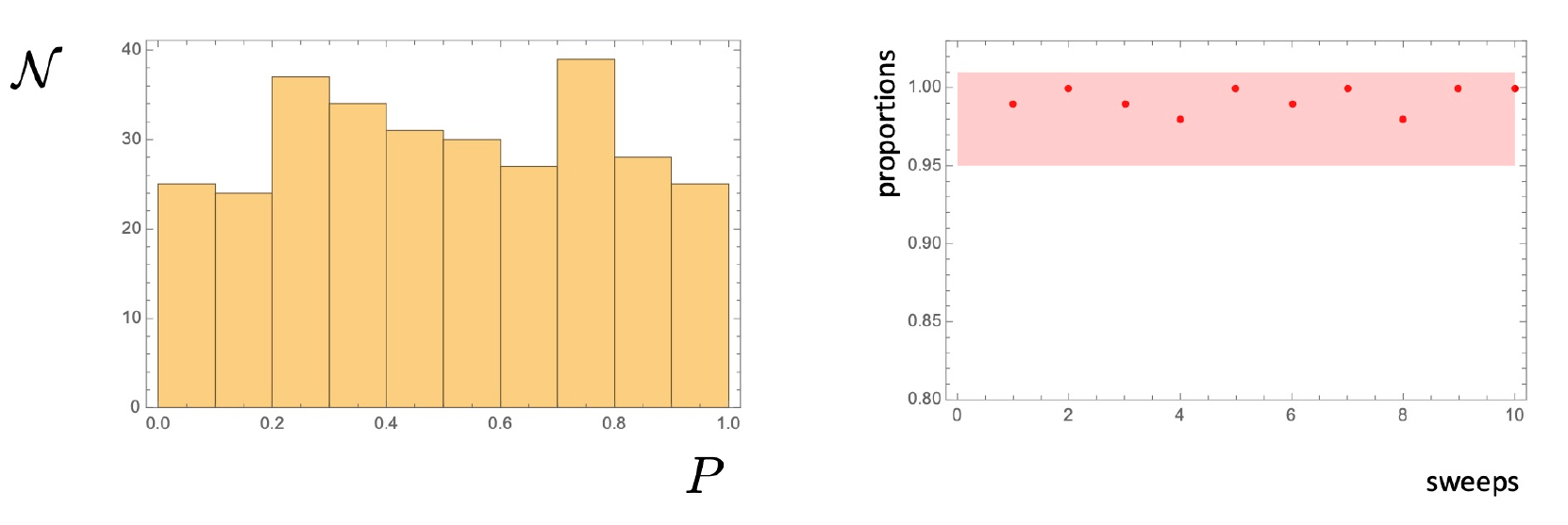}
\caption{(a) Histogram of the $P$-values for the oscillations shown in Figure \ref{oscillationnn1}. (b) Proportions of sequence that pass the statistical test $P$-value $\geq 0.01$. In red the lower and higher values of the confidence interval.}
\label{oscillationnn2}
\end{figure}
Hence, putting all the pieces together, we have 
\beq
\langle (V - \overline V)^2\rangle \,=\, \frac{2 n_0 n_1 (2 n_0 n_1 - n_0 - n_1)}{L^2 (L-1)} \sim 4 L\, \rho^2 (1-\rho)^2 \,\,\,.
\eeq
In light of the formulas given above, in the limit where $L \gg 1$ it is natural to define the normalized variable  
\beq
v \,\equiv\, \frac{V - 2 L \rho (1-\rho)}{2 \rho (1-\rho) \sqrt{L}}
\eeq
and compute the P-value according to 
\beq
{\rm P-value} \,=\, {\rm erfc}\left(\frac{v}{\sqrt{2}}\right) \,\,\,.
\eeq
If the computed $P$-value is $< 0.01$ the sequence has to be considered non-random, otherwise the sequence is promoted to be random.

We have applied this test to sequences of $L=30000$ consecutive values of $\{\mathcal{S}_n\}$ taken randomly from $ N=300$ different non-overlapping stroboscopic intervals $I_{L}(l)$ of length $L = 50000$, where $l \in (L_1,L_2) = (10^{14}, 10^{15})$. The result of this analysis is reported in Figure \ref{oscillationnn1}, where  in Figure  \ref{oscillationnn1}.a we show the histogram of the numbers of values of $V$ obtained in our test while in Figure  \ref{oscillationnn1}.b we show the value obtained in each sweep, together with the intervals relative to $1$, $2$ and $3$ standard deviations. Although most of the values are within $2$ standard deviations, there are few of them which fall beyond $3$ standard deviations, with a $P$-value
which was below $0.01$. Are these deviations significant or are just statistical fluctuations? To answer this question, we have performed the analysis of the proportion of sequences passing the test: using blocks of $15$ sweeps each, we compute the proportion of sequences which pass the test, with the range of acceptable proportion determined by eq.\,(\ref{proportionsequnce}), with $N=30$. Hence, in this case  $({\mathcal I}_-,{\mathcal I}+) = (0.93, 1.04)$. The relative histogram of all obtained $P$-values is reported in Figure \ref{oscillationnn2}.a, while the proportion of sequences passing the test is in Figure \ref{oscillationnn2}.b, in which we see that the few cases which have $P$-values lower than the threshold $0.01$ can be considered statistical fluctuations.

\subsection{Longest Run of $1$'s}

\begin{figure}[t]
\centering\includegraphics[width=1.0\textwidth]{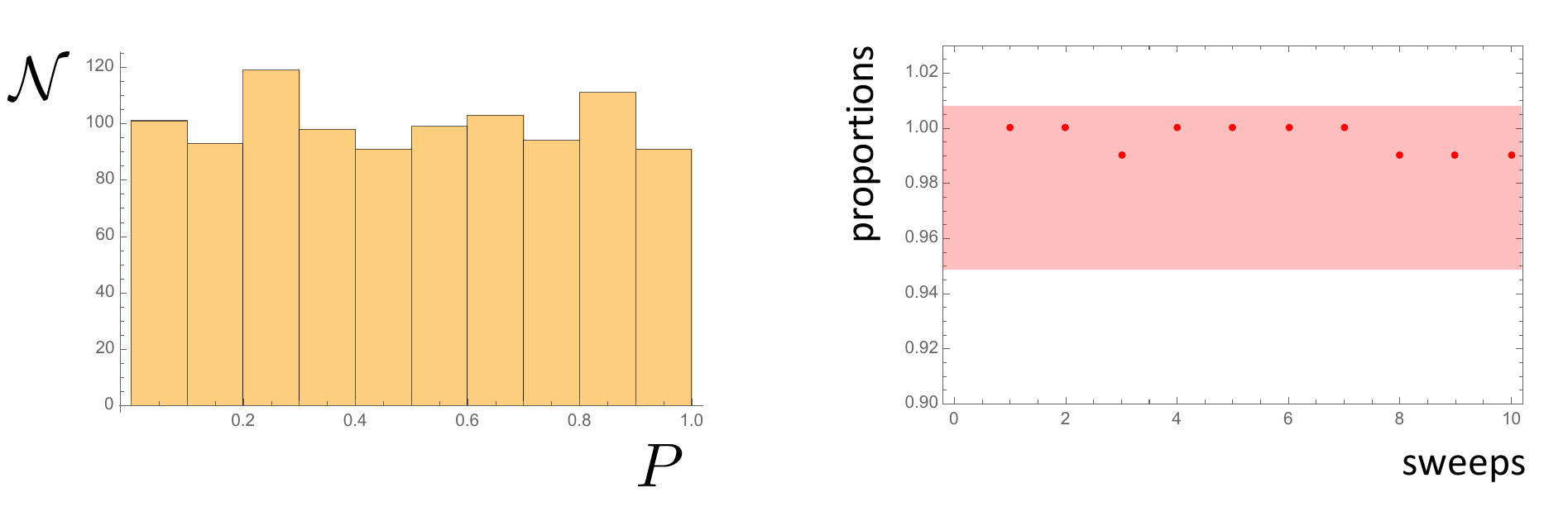}
\caption{(a) Histogram of the $P$-values for the longest runs of $1$'s, with $n=6272$, $M=128$ and $K=5$. (b) Proportions of sequence that pass the statistical test $P$-value $\geq 0.01$. In red the lower and higher values of the confidence interval.}
\label{longestrunss}
\end{figure}

The purpose of this test is to check whether the length of the longest run of $1$'s is compatible with the length expected in a true random sequence. Of course irregular patterns in the longest runs of $1$'s also implies irregular patterns in the longest runs of $0$'s. To perform such a test, first of all we sample stroboscopically our sequence 
$\{{\mathcal S}_n\}$ in terms of our intervals $I_L(l)$. Then, in each interval we extract a subsequence\footnote{The lengths of the variables $n=6272$, $M=128$, $K=5$ and $N=49$ are fixed in these amounts by the request to optimise the statistical test, see \cite{Nist}.} $\{\epsilon\}$ 
 of $n=6272$ consecutive values from $\{{\mathcal S}_k\}$; we divide such a subsequence in $M$-bit blocks (which we take here to be $M=128$) and we measure how all the observed longest run length within the $M$-bit blocks matches the expected longest length within $M$-bit blocks. We take the following $K=5$ categories of events $\nu_k$, for which  is known their relative tabulated probabilities $\pi_k$ (see Table \ref{longprob}) \cite{Nist}.  In this case the corresponding $\chi^2$ variable is given by 
\beq
\chi^2 = \sum_{i=0}^K \frac{(\nu_i - N \pi_i)^2}{N \pi_i} \,\,\,,
\eeq
where $N = 6272/128 = 49$, with the $P$-value given by 
\beq
P-{\rm value} \,=\,Q\left(\frac{K}{2},\frac{\chi^2}{2}\right)\,\,\,.
\eeq
The relative data, reported in Figure \ref{longestrunss}, shows that the sequence $\{ {\mathcal S}_n\}$  also passes successfully  this test. 

\begin{table}[h!]
\begin{center}
\begin{tabular}{||c | |c|| } \hline
 \,$\nu_k$  & $ \pi_k$ \\ \hline
$\nu_0 \,\,\rightarrow $\,\,\,\,\,\, $\# 1's \leq 4\,\,$  & 0.1174 \\
$\nu_1 \,\,\rightarrow$\,\,\,\,\,\, $\# 1's = 5\,\,$& 0.2430\\
$\nu_2 \,\,\rightarrow$\,\,\,\,\,\, $\# 1's = 6\,\,$ &  0.2493\\
$\nu_3 \,\,\rightarrow$\,\,\,\,\,\, $\# 1's = 7\,\,$& 0.1752 \\
$\nu_4 \,\,\rightarrow$\,\,\,\,\,\, $\# 1's = 8\,\,$   & 0.1027\\
$\nu_5 \,\,\rightarrow$\,\,\,\,\,\, $\# 1's \geq 9\,\,$ & 0.1124\\
\hline
 \end{tabular}
\end{center}
\caption{
Categories of events $\nu_k$, related to the long runs of $1$'s in a block whose length is $128$, with the relative tabulated probabilities $\pi_k$. } 
\label{longprob}
\end{table}

\section{Non-overlapping matching test}
We can screen our sequence for the number of occurrences of a given pre-defined aperiodic target $m$-bit string $B$ and compare the actual frequency of this string with its theoretical probability. In more detail, the purpose of this test is to detect whether there are {\em too many} or {\em less} sequences than expected. In other words, this test rejects sequences which exhibits irregular occurrences of a given a-periodic pattern. 
The protocol goes as follows \cite{Barbour}: 
\begin{itemize}
\item we extract a subsequence $\{\epsilon\}$ from $\{{\mathcal S}_n\}$. Such a subsequence is made of the union of $N$ neighborhood blocks of length $L$. 
\item we denote by $W_j$ ($j=1,2,\ldots, N$) the number of times that the target string $B$ occurs within the block $j$. The search for matches proceeds by creating an 
$m$-window on the sequence, comparing the bits within that window against the target $\epsilon$. 
\item if there is {\em no} match, the window slides over one bit. If there {\em is} instead a match, the window slides over $m$ bits. 
\item at the end of the process, for each block we see the final values of the variables $W_j$. 
\end{itemize}
Under the hypothesis of the randomness of the sequence, the expected mean $\mu$  and variance $\sigma^2$ of the variable $W$ are \cite{Nist}
\beq
\mu \,=\, (L - m +1)/2^m 
\hspace{3mm}
,
\hspace{3mm}
\sigma^2 \,=\, L \, \left(\frac{1}{2^m} - \frac{2 m -1}{2^{2m}} \right) \,\,\,.
\eeq
Hence, the $\chi^2$ variable is defined as 
\beq
\chi^2 \,=\, \sum_{j=1}^N \frac{(W_j - \mu)^2}{\sigma^2} 
\,\,\,, 
\eeq
and the corresponding $P$-value is therefore given by 
\beq
{\rm P-value} \,=\, Q\left(\frac{N}{2}, \frac{\chi^2}{2} \right)
\,\,\,.
\eeq
Before presenting our data, an example will help in clarifying the nature of this test. 

\vspace{3mm}
\noindent
{\bf Example}. Let's choose $N=2$ and $L=10$, with $m=3$ and target $m$-string $B = 0 0 1$. Let's imagine that our subsequence is 
\beq
\epsilon = 1 0 0 1 0 0 0 1 1 0 0 1 0 0 0 1 1 0 0 1
\eeq
made of $N=2$ blocks of length $L=10$ 
 \beq
\epsilon_1 = 1 0 0 1 0 0 0 1 1 0 
\hspace{5mm}
,
\hspace{5mm}
\epsilon_2 = 0 1 0 0 0 1 1 0 0 1
\eeq
As we stated, as  target string $B$ we choose $B = 0 0 1$. Hence we have the situation shown in the Table \ref{matchingt1}. Therefore, for this example, we have 
\begin{eqnarray}
\mu & \,= \, & (10 - 3 +1)/2^3 \,=\, 1 \\
\sigma^2 & \,=\, & 10 \,\left(\frac{1}{2^3} - \frac{2 \times 3 -1}{2^{6}}\right) \,=\, 0.46875
\end{eqnarray}
The $\chi^2$ variable is then 
\beq
\chi^2 \,=\, \frac{(2-1)^2 + (2-1)^2}{0.46875} \,=\, 4.26667\,\,.
\eeq
and the corresponding $P$-value is 
\beq
{\rm P-value}\, =\, Q\left(\frac{2}{2}, \frac{4.26667}{2} \right)\,=\, 0.118442\,\,\,.
\eeq

\begin{figure}[b]
\centering\includegraphics[width=1.0\textwidth]{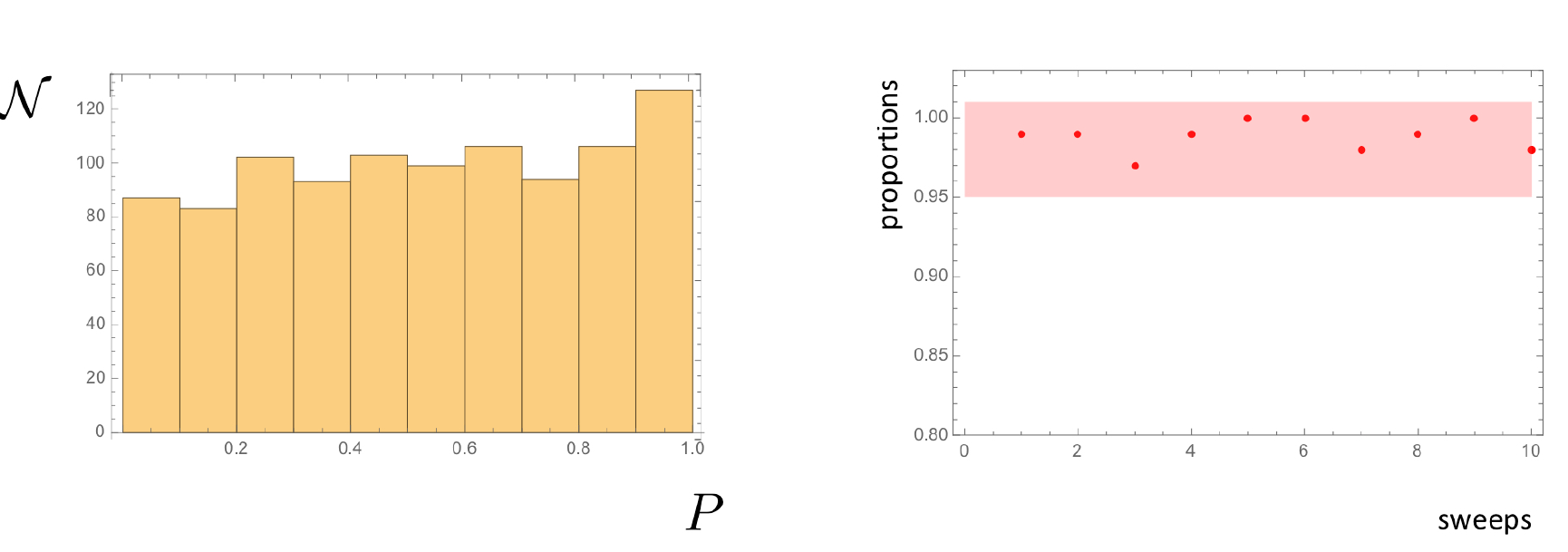}
\caption{(a) Histogram of the $P$-values for the non-overlapping matching test, with $N=80$ blocks with length $L=80$ in $1000$  sequences extracted randomly in the interval of values of the 
sequence $\{{\mathcal S}_n\}$ given by $(L_1,L_2) = (10^{14}, 10^{16})$. The length of the target string $B$ is $m=8$.  (b) Proportions of sequence that pass the statistical test $P$-value $\geq 0.01$. In red the lower and higher values of the confidence interval $(0.95,1.02)$.
}
\label{matchfigure}
\end{figure}

\begin{table}[t]
\begin{center}
\begin{tabular}{| c|  c | c| c | c| }\hline
\multicolumn{1}{|c}{ } &
\multicolumn{2} {|c|} {Block 1} & \multicolumn{2}{|c|}{ Block 2 }\\\cline{2-5} 
 Bit Positions  &\, Bits  &  $W_1$ & Bits &  $W_2$   \\\hline\hline
1-3 &100 & 0 & 010 &  0 \\\hline
2-4 & 001 (hit) &  1 & 010 &  0  \\\hline
3-5 & Not examined &  & 000 &  0 \\\hline
4-6 &Not examined &  & 001 (hit) &  1 \\\hline
5-7 &000 & 1 &  Not examined & \\\hline
6-8 &001 (hit) & 2 & Not examined  &  \\\hline
7-9 &Not examined  & & 100 &  1 \\\hline
8-10 &Not examined &  & 001 (hit)  &  2 \\
\hline
\hline
Total & & 2 & & 2 \\
\hline
\end{tabular}
\end{center}
\caption{
Non-overlapping Matching Test Example for a sequence made of $N=2$ blocks of length $M=10$ and target bit $B =001$ made of $m=3$ bits.} 
\label{matchingt1}
\end{table}

\vspace{3mm}
\noindent
{\bf Results for the sequence $\{{\mathcal S}_n\}$}. In implementing this test on our sequence $\{{\mathcal S}_n\}$ we have chosen an aperiodic string of $m=8$ bits, in particular 
we have employed the target string 
\beq
B\,=\, \{0, 0, 1, 0, 1, 1, 0, 1\}\,\,\,.
\eeq
as well as other target strings such as 
\beq 
\{0,,0,1,1,1,1,0,1\}
\hspace{3mm}
,
 \hspace{3mm}
\{0,0,0,0,1,0,0,1\}
\hspace{3mm}
{\rm or} 
 \hspace{3mm}
 \{0,1,0,1,1,0,1,1\}
\eeq
We have always chosen $N=80$ blocks of length $L=80$ in sequences extracted randomly in the interval of values of the sequence $\{{\mathcal S}_n\}$ given by 
$(L_1,L_2) = (10^{14}, 10^{16})$. The number of subsequences analysed were $1000$. The data relative to the P-values and the proportions of them that pass the statistical test are reported in Figure \ref{matchfigure}: from these plots one can see that our sequence $\{{\mathcal S}_n\}$ passes successfully also this non-overlapping matching test.

\section{Rank of matrices}
\begin{figure}[b]
\centering\includegraphics[width=0.8\textwidth]{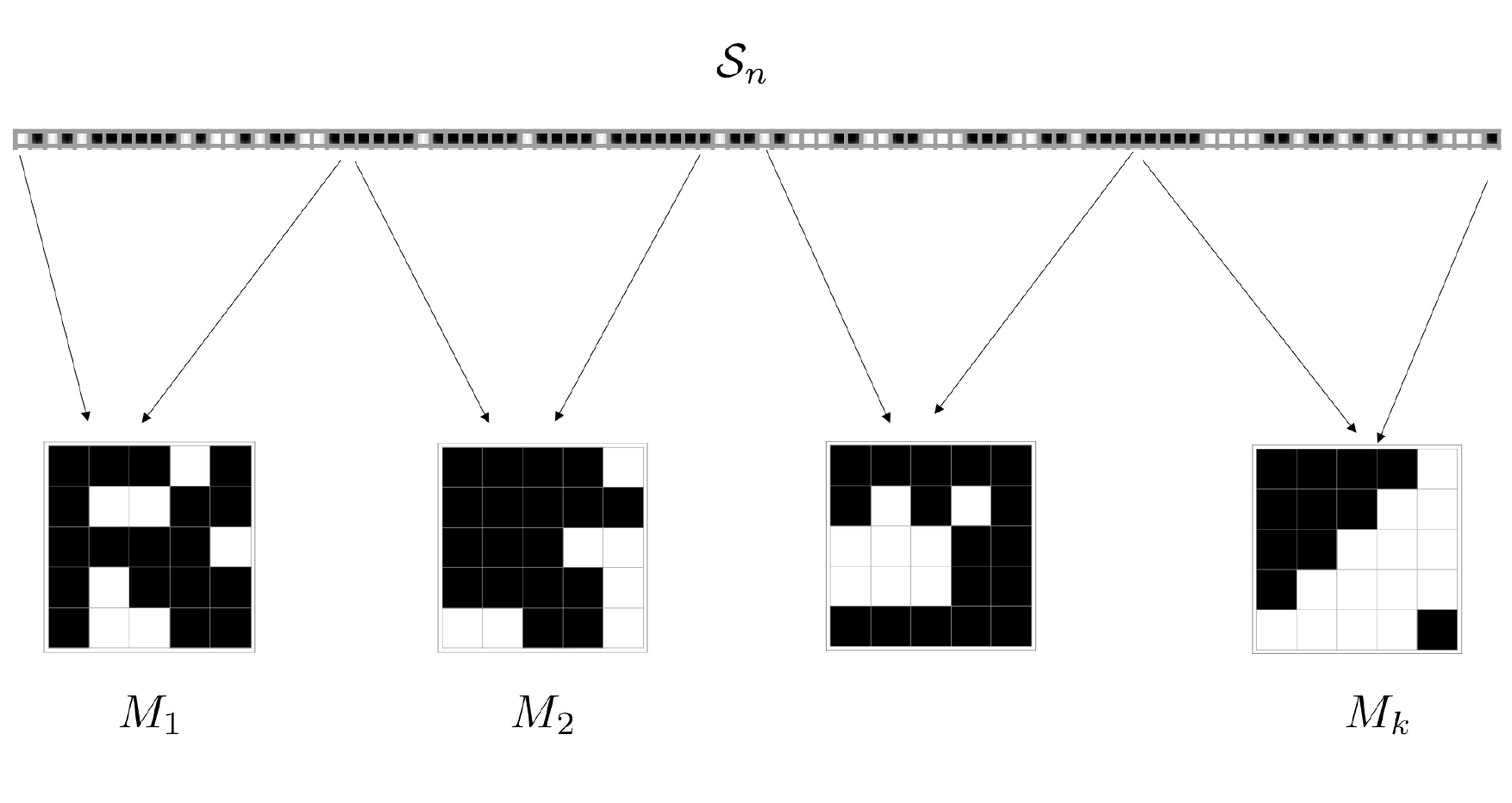}
\caption{
The way of extracting matrices of dimensions $H \times H$ from the original sequence ${\mathcal S}_n$.}
\label{matrixextr}
\end{figure}

Among the different tests which we can make to verify the randomness of our sequence $\{\mathcal{S}_n\}$, one of them concerns the distribution of the rank of binary matrices that can be built up in terms of the elements $\mathcal{S}_n$. The test goes as follows \cite{Marsaglia}:  given an integer length $H$ (in our case $ 10 \leq H \leq 32$) and a segment of length $L\gg H^2$ of our sequence $\{\mathcal{S}_n\}$, we have $N = [L/H^2]$ disjoint blocks of $\{\mathcal{S}_n\}$: each of them can be rearranged in such a way to form a $H \times H$ binary matrix.  Each row of these matrices is filled with successive $H$ terms of the original sequence (see Figure \ref{matrixextr}). Once we have organized our original segment of $\{\mathcal{S}_n\}$ of length $L$ in this way, we can enquire about the probability distribution of the ranks of these matrices on the binary field $\mathbb{F} =\{0,1\}$: this means that, denoting by $v_i$ either the $i$-th binary row or column of these matrices, we can make linear combinations of these vectors but only using as coefficients the two elements of $\mathbb{F} $
\beq
a_1 \, v_1 + a_2 \, v_2 + \cdots a_n \, v_n 
\hspace{5mm}
,
\hspace{5mm}
a_i \in \mathbb{F}
\eeq

\begin{figure}[t]
\centering\includegraphics[width=1.0\textwidth]{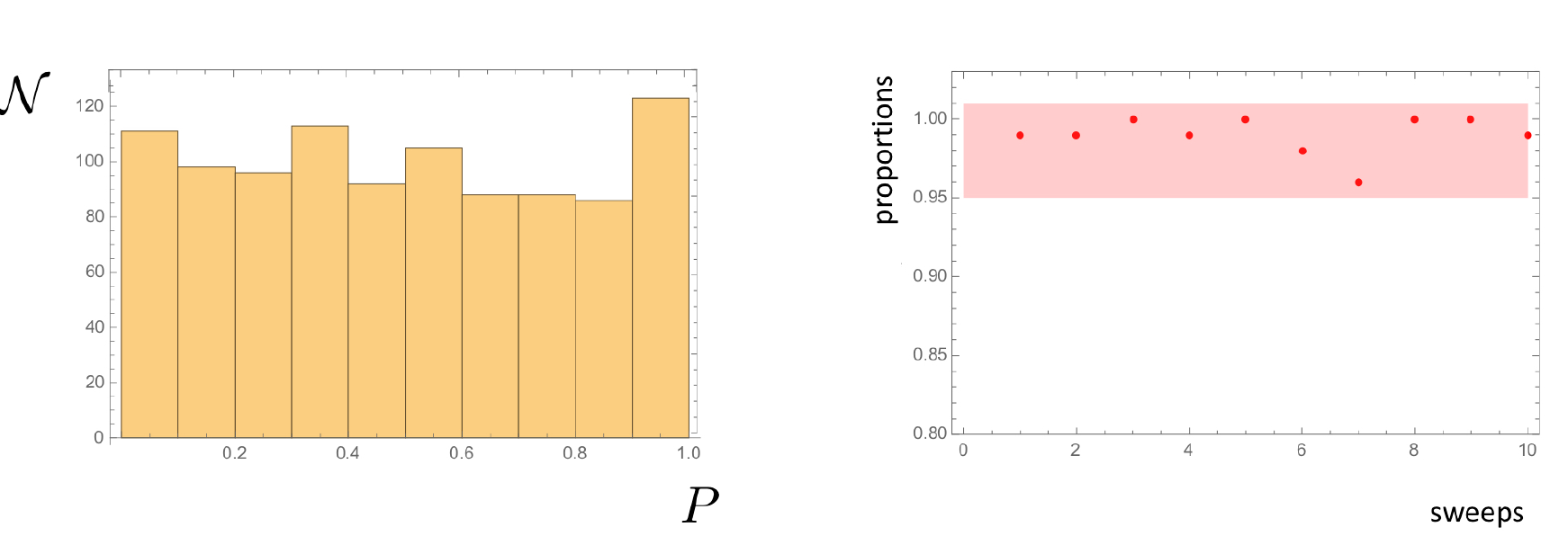}
\caption{(a) Histogram of the $P$-values for the rank of the binary matrices ($10\times 10$) for $1000$ sequences. (b) Proportions of sequence that pass the statistical test $P$-value $\geq 0.01$. In red the lower and higher values of the confidence interval.}
\label{rank10fig}
\end{figure}

\subsection{Probability of the ranks}  For a purely random binary matrix, after making properly linear combination of the rows, its rank is essentially the number of nonzero rows and this can be determined by a simple counting. Hence, the probability that the rank of a random $H \times H$ binary matrix takes the value $r =1,2 , . . . , H$ is given by \cite{MarsagliaTsay}
\beq
P(r) \,=\, 2^{r (2 H -r) -H^2} \,\prod_{i=0}^{r-1} \frac{(1 - 2 ^{i-H})^2}{1 - 2^{i-r}} \,\,\,.
\label{marsagliaform}
\eeq 
Surprisingly enough, the maximum of this probability distribution is not at $r = H$ but at $r = H-1$: this is easily seen taken $H \rightarrow \infty$ limit, where we have 
\begin{eqnarray}
P(H) &\,=\,& \prod_{j=1}^\infty \left(1 - \frac{1}{2^j}\right) \sim 0.2888...\nonumber \\
P(H-1) &\,\sim \,& 2 \, P(H) \,\sim \, 0.5776... \\
P(H-2) &\,\sim\,&\frac{4}{9} \, P(H) \,\sim\, 0.1284...\nonumber
\end{eqnarray}
and the rest of all other probabilities having very small values ($\leq 0.0005$) when $H \geq 10$. Below we denote by $\tilde P$ the probability relative to all other ranks but $r=H$ and $r=H-1$, 
i.e $\tilde P = 1 - P(H) - P(H-1)$.

\subsection{Protocol}  As usual, we have stroboscopically extracted from our sequence $\{\mathcal{S}_n\}$ various subsequences $\{\epsilon\}$ of length $L$ in terms of the intervals $I_L(l)$. 
Here we present the result relative to $1000$  intervals of length ${ L}=10^6$ around $l=10^{15}$ both for matrices with $H=10$ and $H=32$. When $H=10$, we have a sample of 
$N=10^4$ matrices while, when $H=32$, we have a sample of $N=10^3$ matrices. For each interval, we have determined the rank $r$ of each matrix and the frequencies 
$F_H$, $F_{H-1}$ and $N-F_H - F_{H-1}$ of the values $r=H$, $r=H-1$ and the ranks not exceeding $H-2$
\begin{eqnarray*}
&&F_H \,=\, \#\{r_k = H\}  \\
&& F_{H-1} \,=\,\#\{r_k = H-1\} 
\end{eqnarray*}
The chi-square value is easily defined 
\beq
\chi^2 \,=\, \frac{(F_H - P(H) N)^2}{N\,P(H) } +  \frac{(F_{H-1} - P(H-1) N)^2}{N\,P(H-1) } + 
\frac{(N-F_H - F_{H-1} - \tilde P)^2}{N\,\tilde P )}
\eeq
In Figure \ref{rank10fig} we report the distribution of the P-values and the proportions of the subsequences which are within the interval of confidence for the case $H=10$ (for $H=32$ 
the results are analogous). In both cases analysed ($H=10$ and $H=32$), our sequence $\{{\mathcal S}_n\}$ has successfully passed this test.

\section{Maurer Universal statistical test}
A non-random sequence has a certain level of predictability. This permits to compress the sequence without losing information: this idea is at the basis of all file compression, irrespectively of whether they are text files, video or audio files. On the other hand, a purely random sequence presents a  ``strong resistance" to its compression, because there are 
no periodic or predictive patterns that can be used to this aim. Therefore,  the stronger is the resistance of a sequence to its compression,  the higher is its level of randomness.

\begin{figure}[t]
\centering\includegraphics[width=0.8\textwidth]{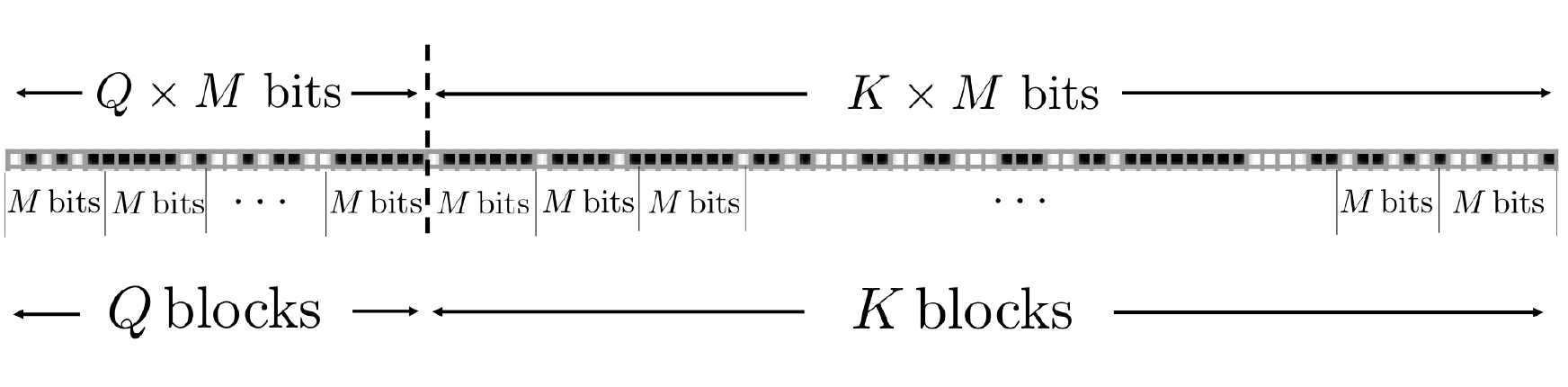}
\caption{
The partition of the $\epsilon$ sequence into $Q$ and $K$ blocks made of $M$-bits non-overlapping sub-blocks.}
\label{partition}
\end{figure}

The purpose of the Maurer Universal statistical test is to detect whether or not a sequence can be significantly compressed keeping its level of information. The protocol goes as 
follows \cite{Maurer}:
\begin{itemize}
\item we denote by  $f_L$ the sum of the ${\rm log}_2$ distances between 
matching $M$-bit templates, i.e. the sum of the number of digits in the distance between 
$M$-bit templates. 
\item we consider a sub-sequence $\epsilon$ of length $L = M (Q+K)$ extracted from our sequence $\{{\mathcal S}_n\}$ and divided into two segments: an 
initialization segment made of $Q$ non-overlapping blocks of length $M$, and a test segment of $K$ non-overlapping blocks also of length $M$, as shown in Figure \ref{partition}. 

As an example, let's assume that 
$
\epsilon = 1 0 0 1 0 0 0 1 1 0 0 1 0 0 0 1 1 0 0 1
$, so that $L=20$. Assume that $M=2$ and $Q=4$, so that $K=6$. The initialization segment is $1 0 0 1 0 0 0 1 $, while the test segment is $1 0 0 1 0 0 0 1 1 0 0 1$. The $M$-bit blocks are shown in the following Table 
\begin{table}[h!]
\begin{center}
\begin{tabular}{| c |  c | c |}\hline
 Type  &\, Block  &  Contents \\\hline\hline
 &  1 & 10 \\\cline{2-3}
Initialization & 2&  01 \\\cline{2-3}
Segment & 3 & 00 \\\cline{2-3}
& 4 & 01 \\\hline
& 5 & 10 \\\cline{2-3}
Test & 6 & 01\\\cline{2-3}
Segment & 7 &  00 \\\cline{2-3}
 & 8 & 00\\\cline{2-3}
& 9 & 11\\\cline{2-3}
& 10 & 01\\
 \hline
\end{tabular}
\end{center}
\end{table}
\item We can use the initialization segment to create a table for each possible $M$-bit value, namely the $M$-bit value is used as an index for the table. The block number of the last occurrence of each $M$-bit block is noted in the table, with the rule that the index $i$ goes from $1$ to $Q$, where $T_j = i$, where $j$ is the decimal representation of the content of the $i$-th $M$-bit block. 

Following the example above, the following Table is created in terms of the $4$ initialization blocks 
\begin{table}[h!]
\begin{center}
\begin{tabular}{| c | c | c |c | c | }\hline
\multicolumn{1}{|c}{ } &
\multicolumn{4} {|c|} {Possible M-bit Value}\\ \cline{ 2-5}
 & 00 & 01 & 10 &  11 \\\cline{2-5}
 & (saved in $T_0$ & (saved in $T_1$)  & (saved in $T_2$) &  (saved in $T_3$)  \\\hline
Initialization & 3 & 4 & 1 & 0 
\\\hline 
\end{tabular}
\end{center}
\end{table}
\item We then examine each of the $K$ blocks in the test segment and determine the number of blocks since the last occurrence of the same $M$-bit block (i.e. $i- T_j$). Replace the value in the table with the location of the current block (i.e. $T_j = i$). Add the calculated distance between re-occurrences of the same $M$-bit block to an accumulating ${\rm log}_2$ sum of all the differences detected in the $K$ blocks, namely ${\rm sum} \rightarrow {\rm sum} + \log_2(i-T_j)$.

Following the example given above, the table and the cumulative sum are constructed as follows: 
\begin{itemize}
\item For block 5 (i.e. the first test block): 5 is placed in the ``10" row of the table (i.e. $T_2$, and 
${\rm sum} = \log_2(5-1) = 2$
\item For block 6: 6 is placed in the ``01" row of the table (i.e. $T_1$) and 
${\rm sum} =  2 + \log_2(6 - 4) = 3$ 
\item For block 7: 7 is placed in the ``00" row of the table (i.e. $T_0$) and 
 ${\rm sum} =  3 + \log_2(7 - 0) = 5.80735$. 
 \item For block 8: 8 is placed in the ``00" row of the table and 
  ${\rm sum} =  5.80735 + \log_2(8-7) = 5.80735$. 
 \item For block 9: 9 is placed in the ``11" row of the table (i.e. $T_3$) and 
   ${\rm sum} =  5.80735 + \log_2(9-0) = 8.97728$. 
 \item For block 10: 10 is placed in the ``01" row of the table and 
 ${\rm sum} =  8.97728 + \log_2(10 -6) = 12.97728$  
 \end{itemize} . 
 The states of the Table are given by
 \begin{table}[h!]
\begin{center}
\begin{tabular}{| c | c | c | c | c |}\hline
\multicolumn{1}{| c |}{ Iteration} &
\multicolumn{4} {| c |} {Possible M-bit Value}\\ \cline{ 2-5}
Blocks & \,\,\,\,00\, & \,\,\,\,01\, & \,\,\,\,10\, & 11 \\\hline
4 &  \,\,\,3 &  \,\,\,4 &  \,\,\,1 & 0 \\\hline
5 & \,\,\,3 & \,\,\,4 & \,\,\,5 & 0 \\\hline
6 & \,\,\,3 & \,\,\,6 & \,\,\,5 & 0 \\\hline
7 & \,\,\,7 & \,\,\,6 &\,\,\, 5 & 0 \\\hline
8 & \,\,\,8 & \,\,\,6 & \,\,\,5 & 0\\\hline
9 & \,\,\,8 & \,\,\,6 & \,\,\,5 & 9\\\hline
10 & \,\,\,8 & \,\,\,10 & \,\,\,5 & 9\\\hline
\end{tabular}
\end{center}
\end{table}
 \item We compute the test statistic variable 
 \beq
 f_L \,=\,\frac{1}{K} \sum_{i=Q+1}^{Q+K} \log_2(i-T_j)\,\,\,,
 \label{snn}
 \eeq
 where $T_j$ is the table entry corresponding to the decimal representation of the content of the $i$-th M-bit block.  
  In our example, $f_L = \frac{12.97728}{6} = 2.16288$. 
 \item We then compute the $P$-value, given in this case by 
 \beq
 {\rm P-value}\, =\, {\rm erfc}\left(\left|\frac{f_L - \mu}{\sqrt{2} \sigma}\right|\right) \,\,\,,
 \eeq
 where $\mu$ and $\sigma^2$ are the theoretical expected values for the mean and the variance of the variable $s_n$. These quantities are tabulated for different values of $M$ and 
 $K$ \cite{Nist}. Since, in our case, we decide to use $M=6$, $Q=640$ and $K=233227$ (these numbers somehow optimize the performance of the test), the theoretical mean and the standard deviation are given by 
 \beq
 \mu \,=\, 5.2177052 
 \hspace{3mm}
 ,
 \hspace{3mm}
 \sigma  = 0.0020213
 \,\,\,.
 \eeq
 
\end{itemize}

\begin{figure}[b]
\centering\includegraphics[width=1.0\textwidth]{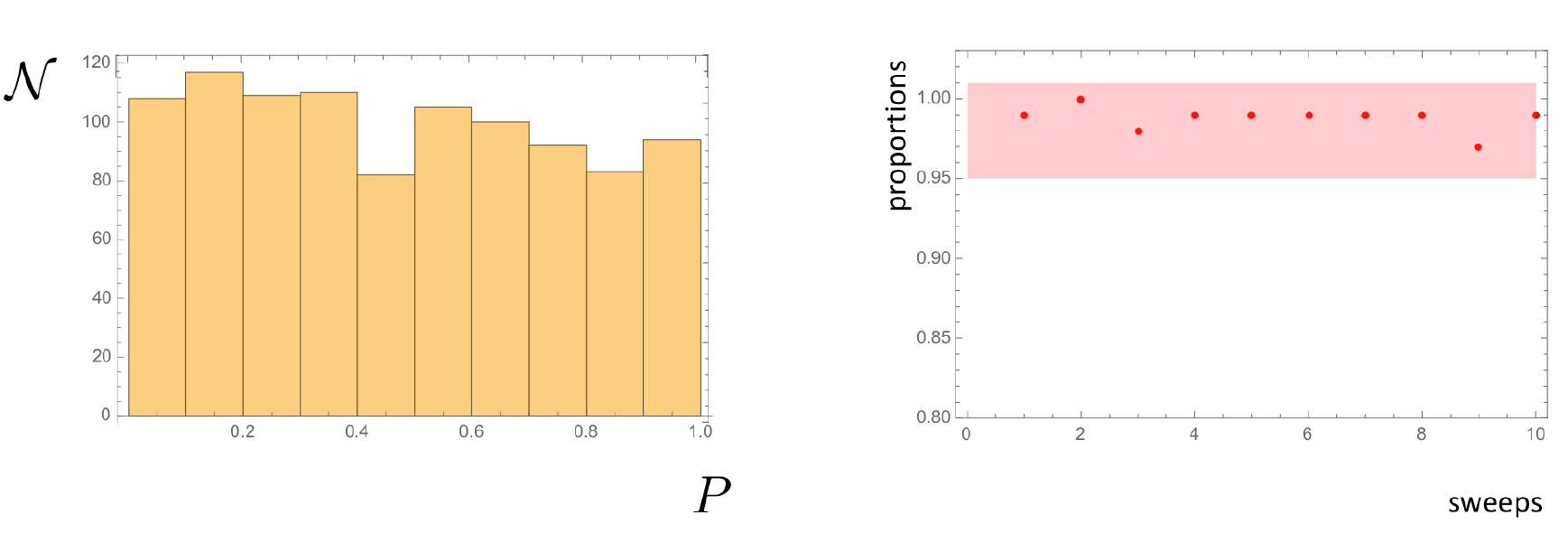}
\caption{(a) Histogram of the $P$-values for the Maurer Universal Statistical Test for $1000$ sequences (b) Proportions of sequence that pass the statistical test $P$-value $\geq 0.01$. In red the lower and higher values of the confidence interval.}
\label{maurerrr}
\end{figure}

We have performed our test on $1000$ sequences $\epsilon$ of total length $L=1.400.000$ extracted from our original sequence $\{{\mathcal S}_n\}$ in the interval 
$(L_1,L_2) = (10^{14}, 10^{16})$. 
The histogram of the $P$-values and the proportions of the sequences that pass the statistical test are shown in Figure \ref{maurerrr}. 
These plots show that our sequence $\{{\mathcal S}_n\}$ successfully passes also this Maurer Universal Statistical Test. This means that our sequence cannot be compressed 
significantly without losing information\footnote{The same conclusion is also reached using standard file compression applications implemented on Mac: applied to large sequences 
extracted by ${\mathcal S}_n$, the compressed files had sizes equal to $99.9 \%$ of the original files. Whoever is interested in performing such a kind of test can request a sample of these sequences by writing to one of us.}.

\section{Entropy test}
How unpredictable is our sequence $\{\mathcal{S}_n\}$? To answer this question, we are going to test the frequency of all possible overlapping $m$-bit pattern across the entire sequence. The aim of this statistical test is to compare the frequency of overlapping blocks of two consecutive/adjacent lengths ($m$ and $m+1$) against the expected result for a purely random sequence. In the following we denote with $\epsilon$ the subsequence of $\{\mathcal{S}_k\}$ to be tested of length $n$. The protocol goes as follows \cite{Pincus}:  
\begin{enumerate}
\item Given the sequence $\epsilon$, we create another one by appending $(m-1)$ bits from the beginning of the sequence to the end of the sequence. 

So, for example, if $\epsilon = 0110010101$ (with $n=10$) and $m=3$, we append the $0$ and $1$ at the beginning of $\epsilon$ to its end, giving rise to the new 
sequence $\epsilon' = 011001010101$. 

\item We make a frequency count of the $n$ overlapping blocks (if a block containing $\epsilon_j$ to $\epsilon_{j+m-1}$ is examined at time $j$, the block containing 
$\epsilon_{j+1}$ to $\epsilon_{j+m}$ is examined at time $(j+1)$). We denote by $C_i^m$ the count of the possible $m$-bit, where $i$ is the $m$-bit value. 

For the example above, the overlapping $m$-bit blocks (here $m=3$) are: $011, 110, 100, 001, 010, 101, 010, 101, 010, 101$. The calculated counts for the $2^m=2^3 = 8$ 
possible $m$-bit strings are 
$$
\# 000 = 0 \,,\, \# 001 = 1 \,,\, \# 010 = 3\,,\, \# 011 = 1\,,\, \# 100 = 1\,,\, \# 101 = 3\,,\, \# 110 = 1 \,,\, \#111=0
$$
\item We compute $C_i^m \,=\, \frac{\# i}{n}$ for each value of $i$. 
For our example, 
$$ 
C_{000}^3 =0 \,,\, C_{001}^3 = 0.1 \,,\, C_{010}^3 = 0.3\,,\, C_{011}^3 = 0.1\,,\, C_{100}^3 = 0.1\,,\, C_{101}^3 = 0.3\,,\, C_{110}^3 = 0.1 \,,\, C_{111}^3=0
$$
\item We then compute 
\beq 
\varphi^{(m)} \,=\, \sum_{i=0}^{2^m -1} \pi_i \log\pi_i
\hspace{3mm}
,
\hspace{3mm}
\pi_i = C_j^3 
\,\,
{\rm and}
\,\, 
j = \log_2 i 
\eeq

\vspace{1mm}
For our example 
\begin{eqnarray*}
\varphi^{(3)} &= & 0 (\log 0) + 0.1 (\log 0.1) + 0.3 (\log 0.3) + 0.1 (\log 0.1) + \\
&+& 0.1 (\log 0.1) + 0.3 (\log 0.3) + 0.1 (\log 0.1) + 0 (\log 0) =
-1.64342
\end{eqnarray*}
\item We then repeat steps 1-4 replacing $m$ with $(m+1)$ 

\vspace{1mm}
Step 1: in our example, $m$ is now $4$ and the sequence to be tested become $\epsilon' = 0110010101011$. 

\vspace{1mm}
Step 2: the overlapping blocks become: 0110, 1100, 1001, 0010, 0101, 1010, 0101, 1010, 1011. The calculated values are 
$$
\# 0110 = 1 \,,\, \# 1100 = 1 \,,\, \# 1001 = 1 \,,\, \# 0010 = 1 \,,\, \#\, 0101 = 2 \,,\, \#1010 = 2 \,,\,\# 1011 =1 
$$
and all the other patterns values are zero. 

\vspace{1mm}
Step 3. Hence we have 
$$
C_{0110}^4 = 0.1 \,,\, C_{1100}^4 = 0.1 \,,\, C_{1001}^4 = 0.1 \,,\, C_{0010}^4 = 0.1 \,,\, C_{0101}^4 = 0.2 \,,\, C_{1010}^4 = 0.2 \,,\,C_{1011}^4 =0.1
$$

\vspace{1mm}
Step 4. Hence we have $\varphi^{(4)} = 5 \times 0.1 (\log 0.1) + 2 \times 0.2 (\log 0.2) =  -1.79507$ 

\item Since for a fixed block length $m$ it is expected that in long random irregular strings, one has 
$$ 
(\varphi^{(m)} - \varphi^{(m+1)}) \sim \log 2 
$$
the limiting distribution of the variable 
\beq 
\chi^2 = n [\log 2 - (\varphi^{(m)} - \varphi^{(m+1)}]\,\,\,,
\eeq
should coincide with that of a $\chi^2$-random variable with $2^m$ degrees of freedom. Therefore the $P$-value is given by 
\beq 
{\rm P-value} \,=\, Q\left(2^{m-1}, \frac{\chi^2}{2}\right) \,\,\,.
\eeq

\vspace{1mm}
In our example, $\chi^2 = 5.41497$ and $Q = (4,2.70748) = 0.712442$
\end{enumerate}

\begin{figure}[t]
\centering\includegraphics[width=1.0\textwidth]{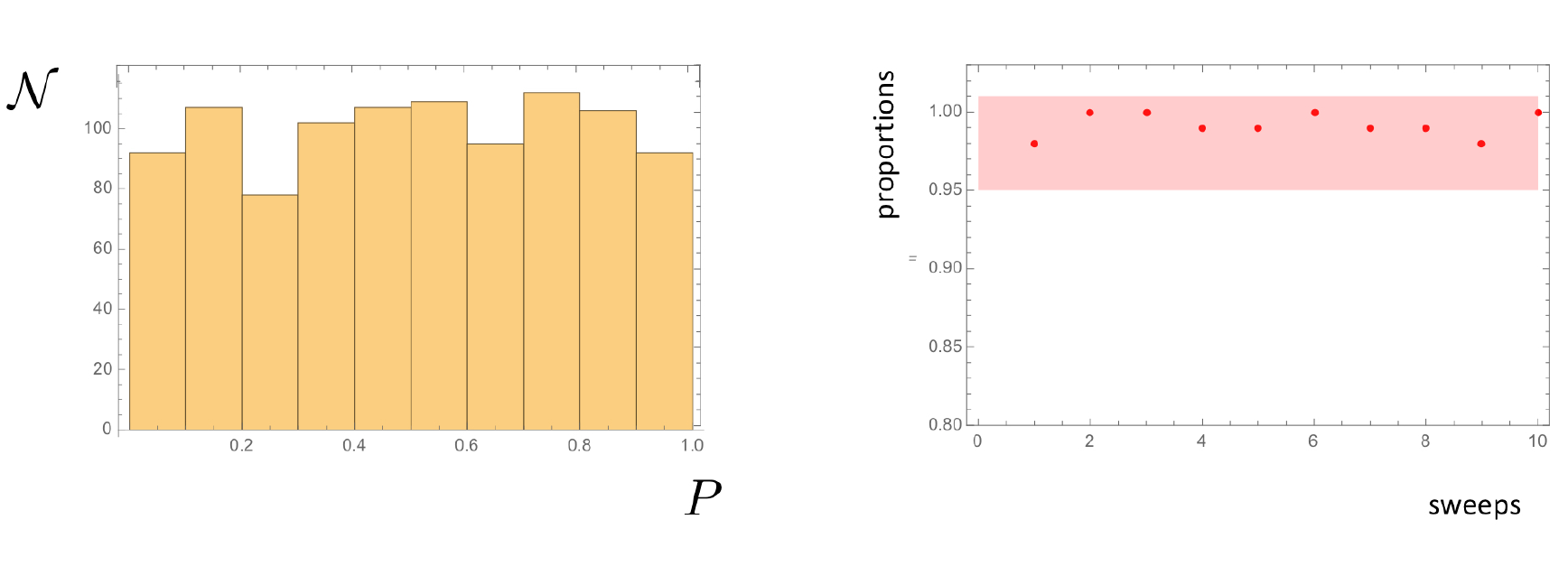}
\caption{(a) Histogram of the $P$-values for the Entropy Test for $1000$ sequences (b) Proportions of sequence that pass the statistical test $P$-value $\geq 0.01$. In red the lower and higher values of the confidence interval.}
\label{entropy}
\end{figure}

We have performed the entropy test on $1000$ subsequences of length $n=1.400.000$ in the interval $(L_1,L_2) = (10^{14},10^{16})$ with a pattern of $m=4$-bits.
The results of the statistical analysis are reported in Figure \ref{entropy} and so, our sequence $\{\mathcal{S}_n\}$ passes successfully also this test. 

\section{Cumulative sum Test }

The aim of this test is to analyze the maximal excursion (from zero) of the random walk defined by the cumulative sum of the subsequences $\{\hat\epsilon\}$ of 
 $\{\hat\mu_n\}$  with values $\{\pm 1\}$. The feature under scrutiny is whether the cumulative sum of the partial sequences occurring in the tested sequence is too large or too small with respect to the expected behaviour of the cumulative sum of a purely random walk. While for a random sequence, the excursion of the random walk should be near to zero, for non-random sequences one expects that the excursion from zero may be significant. The theoretical framework behind this test is the standard random walk (see, for instance   \cite{Yuval,Mazo,Rudnick,Levy,diffusion}). 

\begin{figure}[t]
\centering\includegraphics[width=1.0\textwidth]{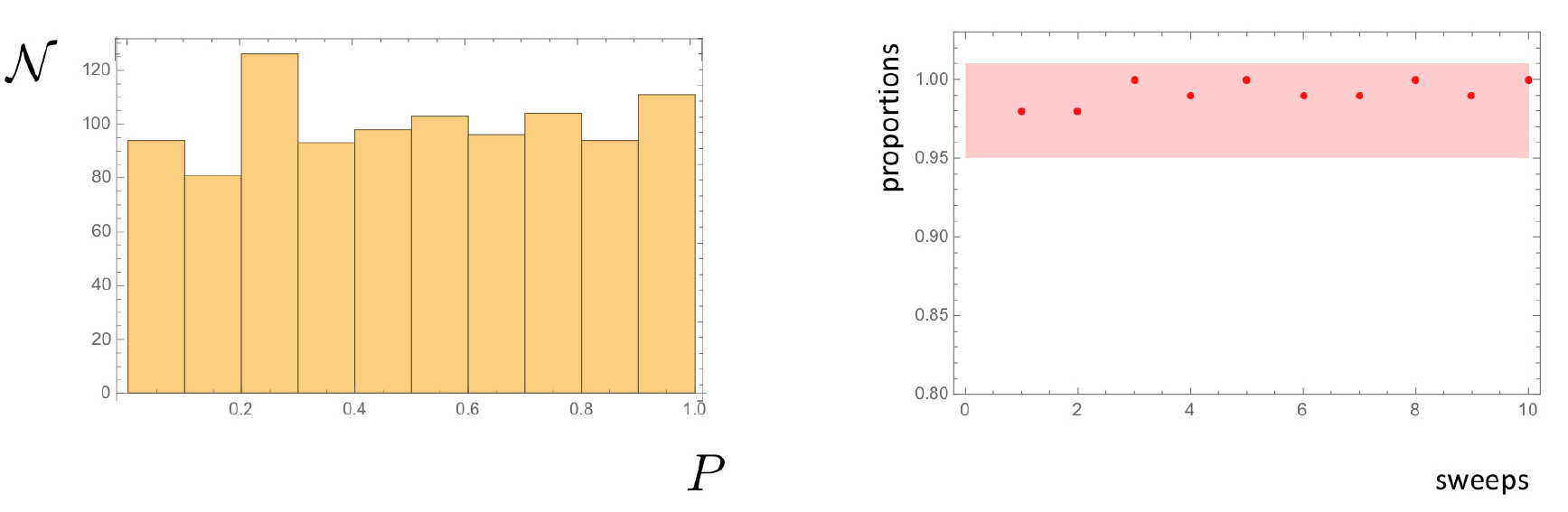}
\caption{(a) Histogram of the $P$-values for the Cumulative Sum Test for $1000$ sequences (b) Proportions of sequence that pass the statistical test $P$-value $\geq 0.01$. In red the lower and higher values of the confidence interval.}
\label{cumsum}
\end{figure}

To implement this test let's consider the block variables, defined earlier in eq.(\ref{groupvariables})
\beq 
B_L(\s) \,=\,\sum_{k\in I_L(\s)}  \hat\mu(k) \,=\, \sum_{k=\s}^{\s + L-1}  \hat\mu(k)  \,\,
\eeq
The statistical variable of the test is 
\beq
t\,=\, {\rm max}_{1\leq k \leq L} |B_L(l) |\,\,\,,
\eeq
i.e. the largest of the absolute values of the partial sums $B_L(l)$. From the theory of the random walk, the limiting distribution of the absolute values of the partial sum
is 
\begin{eqnarray}
&& G(z) \,\equiv\, \lim_{L\rightarrow \infty} P\left(\frac{{\rm max}_{1\leq k \leq L} |B_L(l)|}{\sqrt{L}} \leq z\right) \,=\, \frac{1}{\sqrt{2 \pi}} \,
\int_{-z}^{z} \sum_{k = - \infty}^{\infty} (-1)^k \,\exp\left[- \frac{(u- 2 k z)^2}{2}
\right] \, du \nonumber \\
&& = \sum_{k=-\infty}^{k=\infty} (-1)^k \left\{\Phi[(2 k+1) z] - \Phi[(2 k-1) z]\right\} \nonumber \\
&&= \Phi(z) - \Phi(-z)- 2 \sum_{k=1}^\infty \left\{ 2 \Phi[(4 k -1) z] - \Phi[(4 k+1) z] - \Phi[(4 k -3) z]\right\} \\
&&\sim \Phi(z)- \Phi(-z) - 2 \{ 2 \Phi(3 z) - \Phi(5 z) - \Phi(z)\} \nonumber \\
&& \sim 1 - \frac{4}{\sqrt{2\pi} z} \,\exp\left[-\frac{z^2}{2}\right] \nonumber
\hspace{3mm}
,
\hspace{3mm}
z \rightarrow \infty
\end{eqnarray}
where, above, $\Phi(x)$ is the standard $erf(x)$ function 
\beq
\Phi(x) \,\equiv\, \frac{1}{\sqrt{2\pi}} \int_{-\infty}^x \exp\left[-\frac{u^2}{2}\right] \, du \,\,\,.
\eeq
The corresponding $P$-value is given by 
\beq
{\rm P-value} \,=\, 1 - G\left({\rm max}_{1 \leq k \leq L} |B_L | /\sqrt{L}\right) \,\,\,.
\eeq
We have analyzed $1000$ subsequences $\{\hat\epsilon\}$ of length $L=1.400.000$ in the interval $(L_1,L_2) = (10^{14},10^{16})$ and the histogram of the corresponding $P$-values 
are reported in Figure \ref{cumsum}, together with the proportion of the sequences that pass the statistical test $P$-value $\geq 0.01$. Also in this case, the conclusion is that the 
sequence $\{{\mathcal S}_n\}$ passes successfully this test. 

\section{Random Excursions Test}
In this test we consider the subsequences $\{\hat\epsilon\}$ with values $\{\pm 1\}$ and we test them versus certain predictions which come from purely random walk theory. 
The protocol of this test is as follows \cite{baron}: 
\begin{figure}[b]
\centering\includegraphics[width=0.5\textwidth]{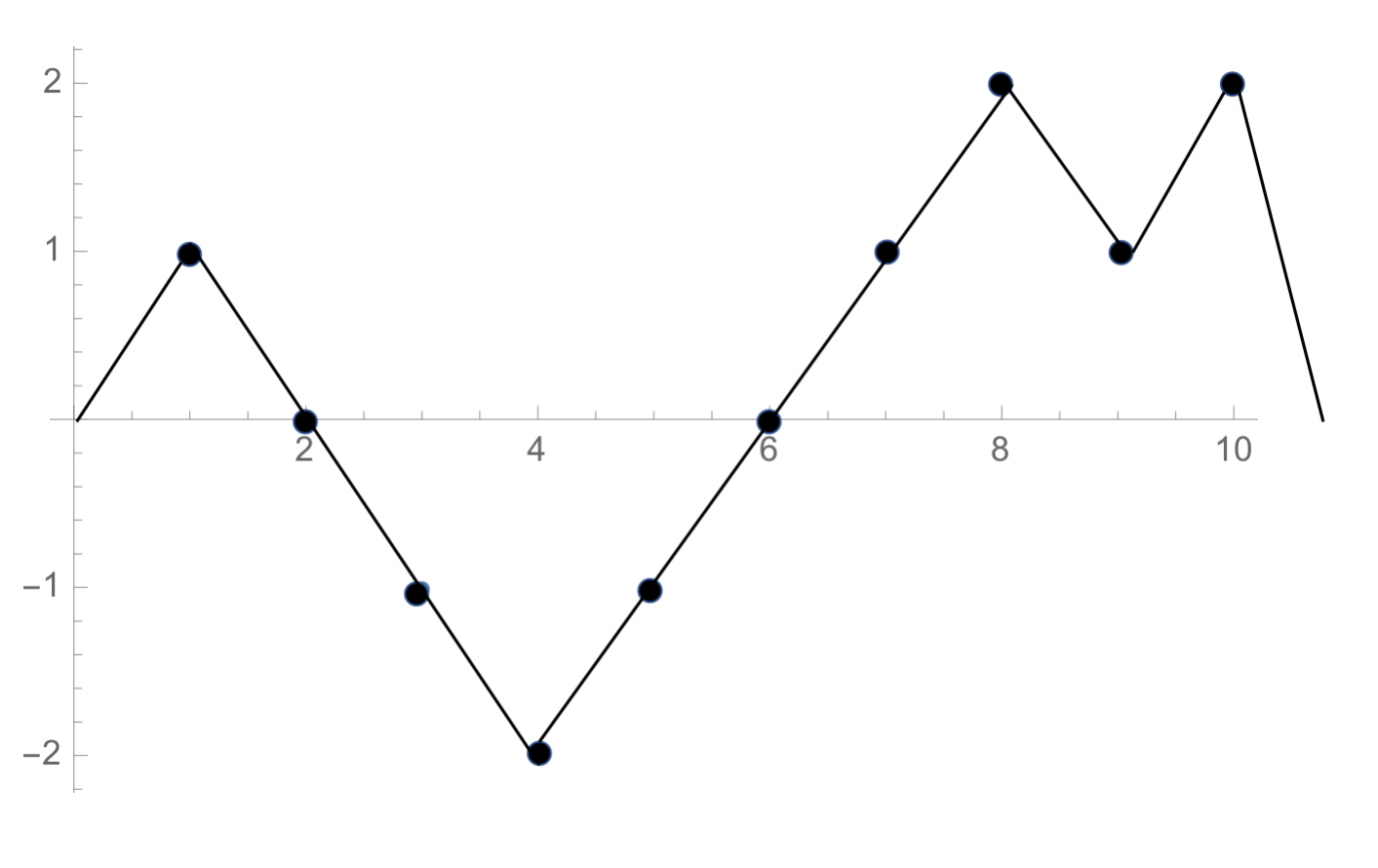}
\caption{(Example of random walk. The $J=3$ cycles are: (0, 1,0), (0,-1,-2,-1,0) and (0,1,2,1,2,0).}
\label{cycle}
\end{figure}
 \begin{table}[t]
\begin{center}
\begin{tabular}{| c | c | c | c  |}\hline
\multicolumn{1}{| c |}{ } &
\multicolumn{3} {| c |} {Cycles}\\ \cline{ 2-4}
State & \,\,\,\,Cycle 1\, & \,\,\,\,Cycle 2 \,  & Cycle 3 \\\cline{2-4}
x & (0,1,0) & (0,-1,-2,-1,0) & (0,1,2,1,2,0) \\\hline
-4 &  \,\,\,0 &  \,\,\,0 &  \,\,\,0 \\\hline
-3  & \,\,\,0 & \,\,\,0 & \,\,\,0 \\\hline
-2 & \,\,\,0 & \,\,\,1 & \,\,\,0 \\\hline
-1 & \,\,\,0 & \,\,\,2 &\,\,\, 0  \\\hline
1 & \,\,\,1 & \,\,\,0 & \,\,\,2 \\\hline
2 & \,\,\,0 & \,\,\,0 & \,\,\,2 \\\hline
3 & \,\,\,0 & \,\,\,0 & \,\,\,0 \\\hline
4 & \,\,\,0 & \,\,\,\,0 & \,\,\,0 \\\hline
\end{tabular}
\end{center}
\caption{Occurrence of the various values $x$ in each cycle.}
\label{cyclemat1}
\end{table}

\begin{enumerate}
\item we consider a subsequence $\{\hat\epsilon\} = \{\hat\epsilon_1, \hat\epsilon_2, \ldots, \hat\epsilon_L\} $ of length $L$ extracted from the sequence $\{\hat\mu\}$. 
\item we consider the partial sums $T_k = \sum_{j=1}^k \hat\epsilon_j$ of our sequence. We denote by $T'$ the total sum $T_L = \sum_{j=1}^L \hat\epsilon_n$ where we have added zeros before and after 
$T_L$, i.e. $T'=0, T_1, T_2, \ldots T_L, 0$. 
\item we count the number $J$ of zero crossings in $T'$, where a zero crossing is a value of zero in $T'$ that occurs after the starting point. $J$ is also the number of cycles in $T'$, where a cycle of $T'$ is a subsequence of $T'$ consisting of an occurrence of zero, followed by no-zero values and ending with another zero. Of course the ending zero in one cycle may the beginning zero in another cycle. For statistical reasons which will explain better below\footnote{For a $\chi^2$ test, the number of variables must be at least larger than $5$ and this is the reason of the constraint on $J$.}, we have to be sure than the number $J$ of zeros in the sequence $T'$ is larger than 800. Looking at the Figure \ref{cycle}, in the plot there are 
$j=3$ cycles given by the sequences  (0, 1,0), (0,-1,-2,-1,0) and (0,1,2,1,2,0). 
\item for each cycle, we compute how many times the sequence $T'$ passes by each of the integer values $x$, in the range $- 4 \leq x \leq -1$ and $1 \leq x \leq 4$. For the sequence shown in Figure \ref{cycle}, we have the situation summarised in Table \ref{cyclemat1}.   
\item for each of the eight states of $x$, we compute $\nu_k(x)$= the total number of cycles in which state $x$ occurs exactly $k$ times among all cycles, for $k=0, 1, \ldots, 5$ (for $k =5$, all frequencies $\geq 5$ are stored in $\nu_5(x)$). We have the sum rule $\sum_{k=0}^5 \nu_k(x) = J$. For the example of Figure \ref{cycle}, the situation is summarised in Table  \ref{cyclemat2}. 

 \begin{table}[t]
\begin{center}
\begin{tabular}{| c | c | c | c  | c | c | c| }\hline
\multicolumn{1}{| c |}{ } &
\multicolumn{6} {| c |} {Number of Cycles}\\ \cline{ 2-7}
State x & \,\,\,\,0 \, & \,\,\,\, 1 \,  & \,\,\,\, 2 & \,\,\,\, 3 & \,\,\,\, 4 & \,\,\,5 \\ \hline
-4 &  \,\,\,3 &  \,\,\,0 &  \,\,\,0 & \,\,\,0 & \,\,\,0  & \,\,\,0\\\hline
-3  & \,\,\,3  & \,\,\,0 & \,\,\,0  & \,\,\,0 & \,\,\,0 & \,\,\,0 \\\hline
-2 & \,\,\,2 & \,\,\,1 & \,\,\,0 & \,\,\,0 & \,\,\,0 & \,\,\,0 \\\hline
-1 & \,\,\,1 & \,\,\,0 &\,\,\, 1  &\,\,\,0 & \,\,\,0 & \,\,\,0\\\hline
1 & \,\,\,1 & \,\,\,1 & \,\,\,1 & \,\,\,0 & \,\,\,0 & \,\,\,0\\\hline
2 & \,\,\,2 & \,\,\,1 & \,\,\,0 & \,\,\,0 & \,\,\,0 & \,\,\,0\\\hline
3 & \,\,\,3 & \,\,\,0 & \,\,\,0 & \,\,\,0 & \,\,\,0 & \,\,\,0\\\hline
4 & \,\,\,3 & \,\,\,\,0 & \,\,\,0 & \,\,\,0 & \,\,\,0 & \,\,\,0\\\hline
\end{tabular}
\end{center}
\caption{Occurrence of the various values $x$ in each cycle.}
\label{cyclemat2}
\end{table}

\begin{figure}[b]
\centering\includegraphics[width=1.0\textwidth]{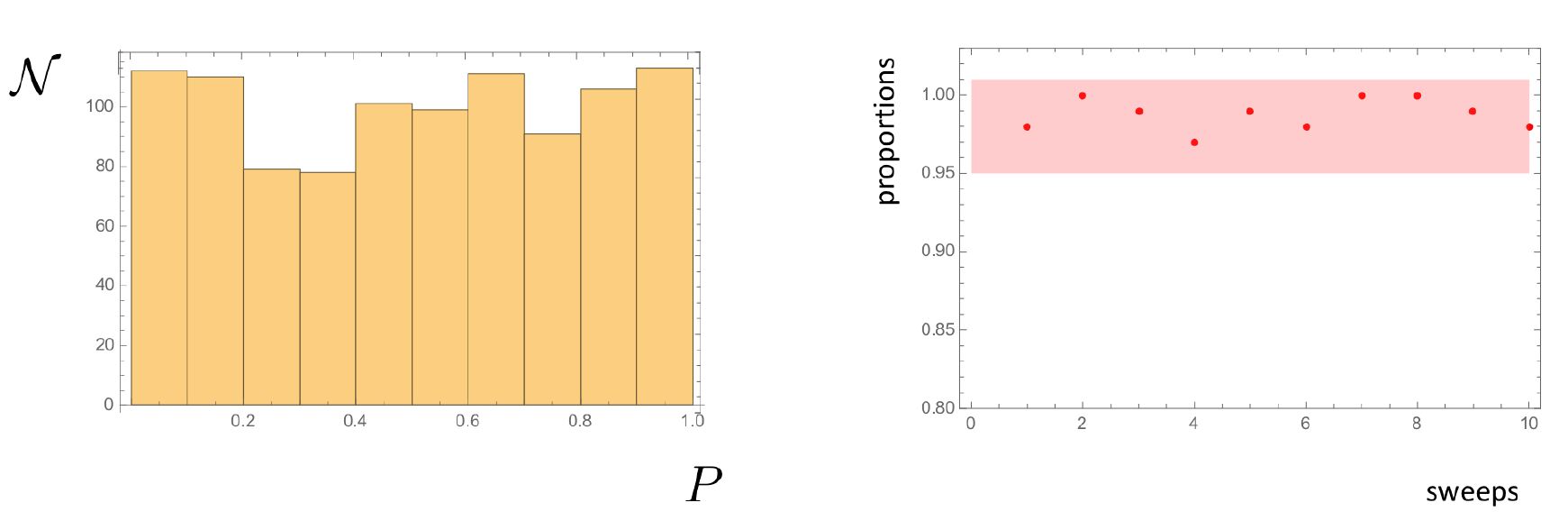}
\caption{(a) Histogram of the $P$-values for the Cumulative Sums ($x=3$)  for $1000$ sequences (b) Proportions of sequence that pass the statistical test $P$-value $\geq 0.01$. In red the lower and higher values of the confidence interval.}
\label{cumsumres}
\end{figure}

\item The theoretical prediction of these frequencies (here denoted as $\pi_k$) are known \cite{baron} and are given by 
\begin{eqnarray} 
&& \pi_0 \,=\, 1 - \frac{1}{2 |x|} \nonumber \\
&& \pi_k \,=\, \frac{1}{4 x^2} \left(1 - \frac{1}{2 |x|} \right)^{k-1} 
\hspace{3mm}
,
\hspace{3mm} k = 1,2, 3, 4 \\
&& \pi_5 \,=\, \frac{1}{2 |x|} \left(1 - \frac{1}{2 |x|}\right)^4 \nonumber 
\end{eqnarray}
\item Using these theoretical and the observed frequencies we can define the $\chi^2(x)$ variables (for each value $x$), given by 
\beq
\chi^2(x) \,=\, \sum_{k=0}^5 \frac{(\nu_k(x) - J \,\pi_k(x))^2}{J \,\pi_k(x)} \,\,\,.
\eeq
This variable is expected to have a $\chi^2$-distribution with $5$ degrees of freedom. This is valid as far as $J \,{\rm min}\,\pi_k(x) \geq 5$. Since ${\rm min}\,\pi_k(x) = 0.00735$ 
(corresponding to $\pi_5(4))$, this lead us to the constraint on $J$ and to the necessity take into consideration only those sequences with $J \geq 800$. 
\item The corresponding $P$-value of this test is given by 
\beq
{\rm P-value} \,=\, Q\left(\frac{5}{2},\frac{\chi^2(x)}{2}\right) \,\,\,.
\eeq
\end{enumerate}
Here we present in Figure \ref{cumsumres} 
the distribution of the $P$-value for $x=3$ since all other values present similar behaviour. We have analyzed  $1000$ non-overlapping subsequences of length 
$L=1.400.000$ extracted by the sequence $\{\hat\mu_n\}$ in the interval $(L_1,L_2) = (10^{14},10^{16})$. As evident from Figure \ref{cumsumres}, our sequence 
$\{{\mathcal S}_n\}$ passes successfully also this test.

\section{Time distribution tests}
\begin{figure}[b]
\centering\includegraphics[width=0.5\textwidth]{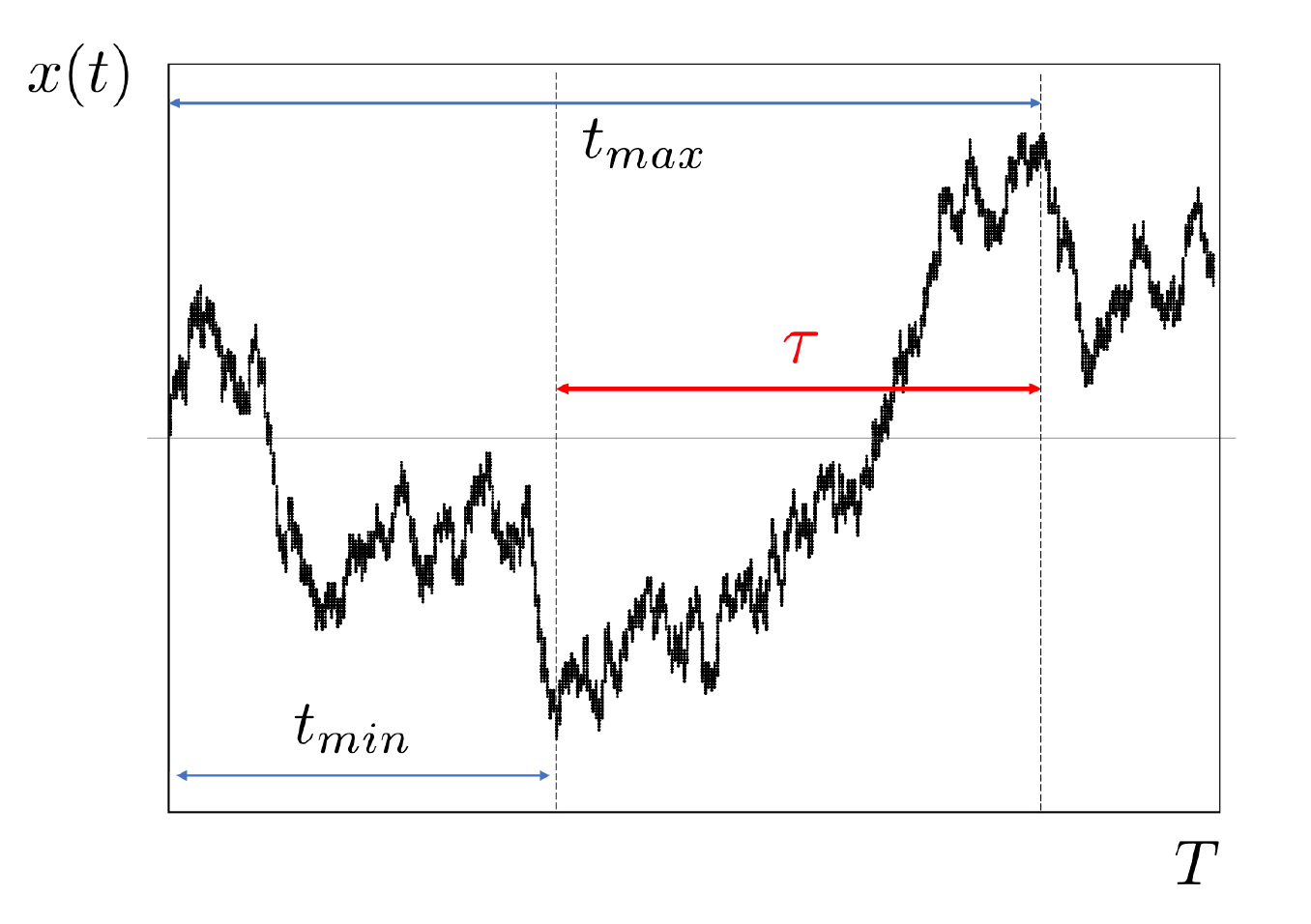}
\caption{Random Walk, whose length is $T$, and the times $t_{min}$, $t_{max}$ and $\tau$. 
}
\label{minmaxrandom}
\end{figure}

Consider a true one-dimensional random walk (see Figure \ref{minmaxrandom}) whose starting point is at the origin and which is long $T$. In the interval $(0, T)$ we can identify the time $t_{max}$ in which the random walk reaches its maximum and the time $t_{min}$ in which it reaches its minimum. In the following we denote by $\tau$ the time interval between the maximum and the minimum, $\tau = t_{max} - t_{min}$. 

\subsection{Probability distribution of ${\bf t_{max}}$ and ${\bf t_{min}}$} 
The probability distributions $P(t_{max} | T)$ and $P(t_{min} | T)$ are exactly known for the pure random walk case \cite{Feller,Yuval}: first of all, for symmetry reasons, they are equal and given by 
\beq
P(x | T) \,=\, \frac{T}{\pi} \,\frac{1}{\sqrt{x (1-x)}} \,\,\,,
\label{distrimm}
\eeq
where $x$ is the scaling variable $x=t_{max}/T = t_{min}/T$. The distribution is symmetric with respect $x \rightarrow (1-x)$ and its first moments are 
\beq
\begin{array}{cccc}
\langle x \rangle & \,=\, & 1/2 &\,=\, 0.5  \\
\langle x^2 \rangle & \,=\, & 3/8 &\,=\,0.375 \\
\langle x^3 \rangle & \,=\, & 5/16 &\,=\,0.3125  \\
\langle x^4 \rangle & \,=\, & 35/128  &\,=\,0.2734 \\
\langle x^5 \rangle & \,=\, & 63/256 & \,=\,0.2461\\
\langle x^6 \rangle & \,=\,& 231/1024 &\,=\,0.2255
\end{array}
\eeq

We have sampled many different intervals of the sequence $\{\mathcal{S}_n\}$ and compared the experimental distribution of $t_{min}$ versus the theoretical one given in 
eq.\,(\ref{distrimm}). 
Here we show the result relative to an interval of $10^8$ values around $x=10^{15}$: we have divided this interval into $N=2\times 10^4$ segments of length $T = 5\times 10^3$ and 
for each of these intervals we have determined $t_{min}$. The comparison between the theoretical and experimental distributions are shown in Figure \ref{minran}: the $\chi^2$ variable of this example is $\chi^2 = 22.35$ (for $20$ degrees of freedom) and the relative $P$-value =$0.60$, confirming the matching between the distribution coming from random walk and Mertens function. The experimental values of the first moments for the data shown in Figure \ref{minran} are 
\begin{eqnarray}
&& \langle x \rangle_{ex} \sim 0.4894
\hspace{3mm}
,
\hspace{3mm}
\langle x^2 \rangle_{ex} \sim 0.369
\hspace{3mm}
,
\hspace{3mm}
\langle x^3 \rangle_{ex}  \sim 0.307 \\
&&
\langle x^4 \rangle_{ex}  \sim 0.267 
\hspace{3mm}
,
\hspace{3mm}
\langle x^5 \rangle_{ex}  \sim 0.240
\hspace{3mm}
,
\hspace{3mm}
\langle x^6 \rangle_{ex}  \sim 0.219 \nonumber
\end{eqnarray}
which differ just of few percents from the theoretical ones. {We have also performed the Kolgomorv-Smirnov test to check the goodness of the expected function
(\ref{targettt}) versus the data: with a sample of $n$ points, the Kolmogorov-Smirnov variable is defined as 
\beq
D_n = {\rm sup}_x |U_n - U(x)|\,\,\,,
\label{KSSS}
\eeq
where $U(x)$ is the cumulative distribution function associated to (\ref{distrimm}) while $U_n$ is the ``experimental" cumulative distribution function with $n$ points. In our case, 
$n=20$ and we got $D_{20} = 0.0041343$ which gives a probability larger than $99.9\%$ that the experimental data follow indeed the expected theoretical function (\ref{distrimm}).

\begin{figure}[t]
\centering\includegraphics[width=0.6\textwidth]{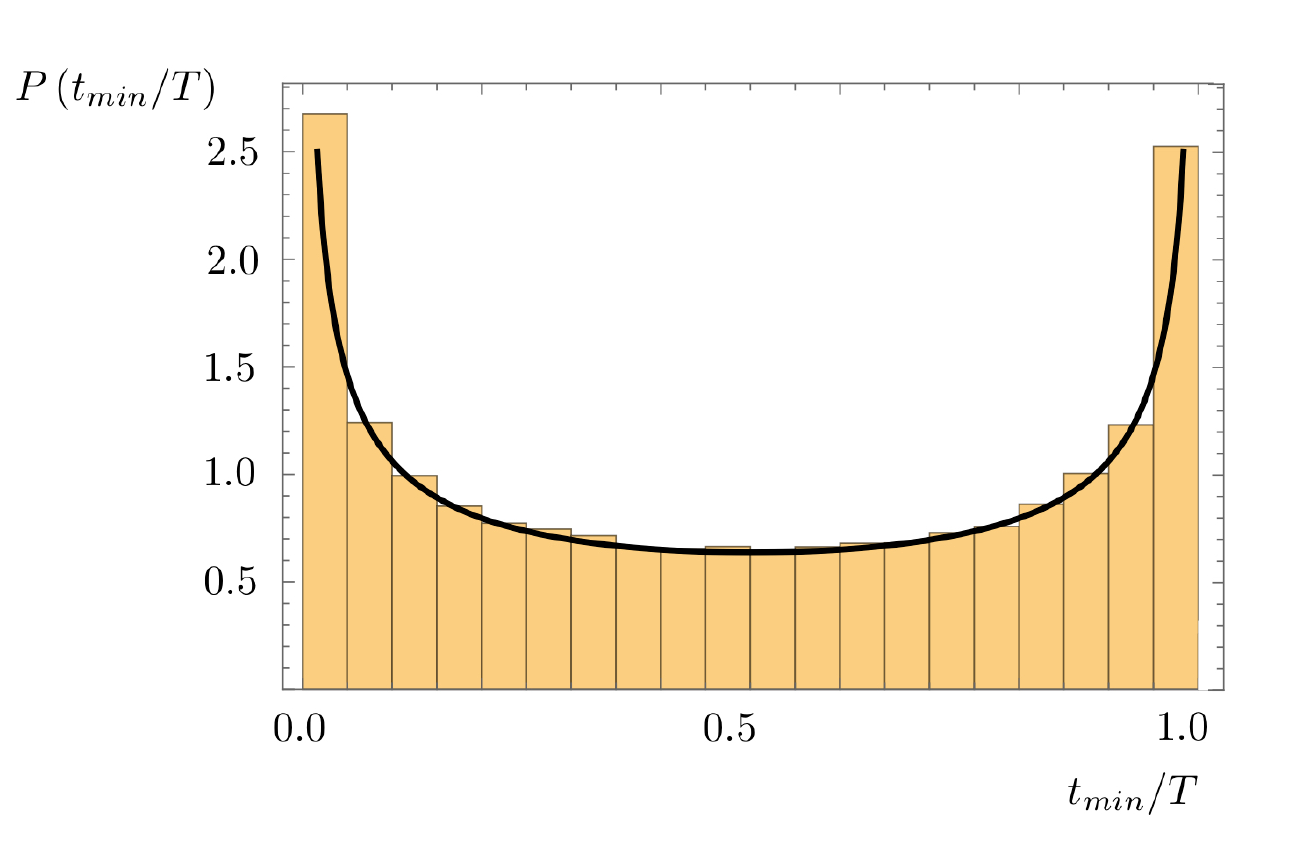}
\caption{The scaled distribution $ P(\tau/T)$ (black solid line) together with the experimental distribution of the $t_{min}/T$ variable. }
\label{minran}
\end{figure}
\begin{figure}[b]
\centering\includegraphics[width=0.6\textwidth]{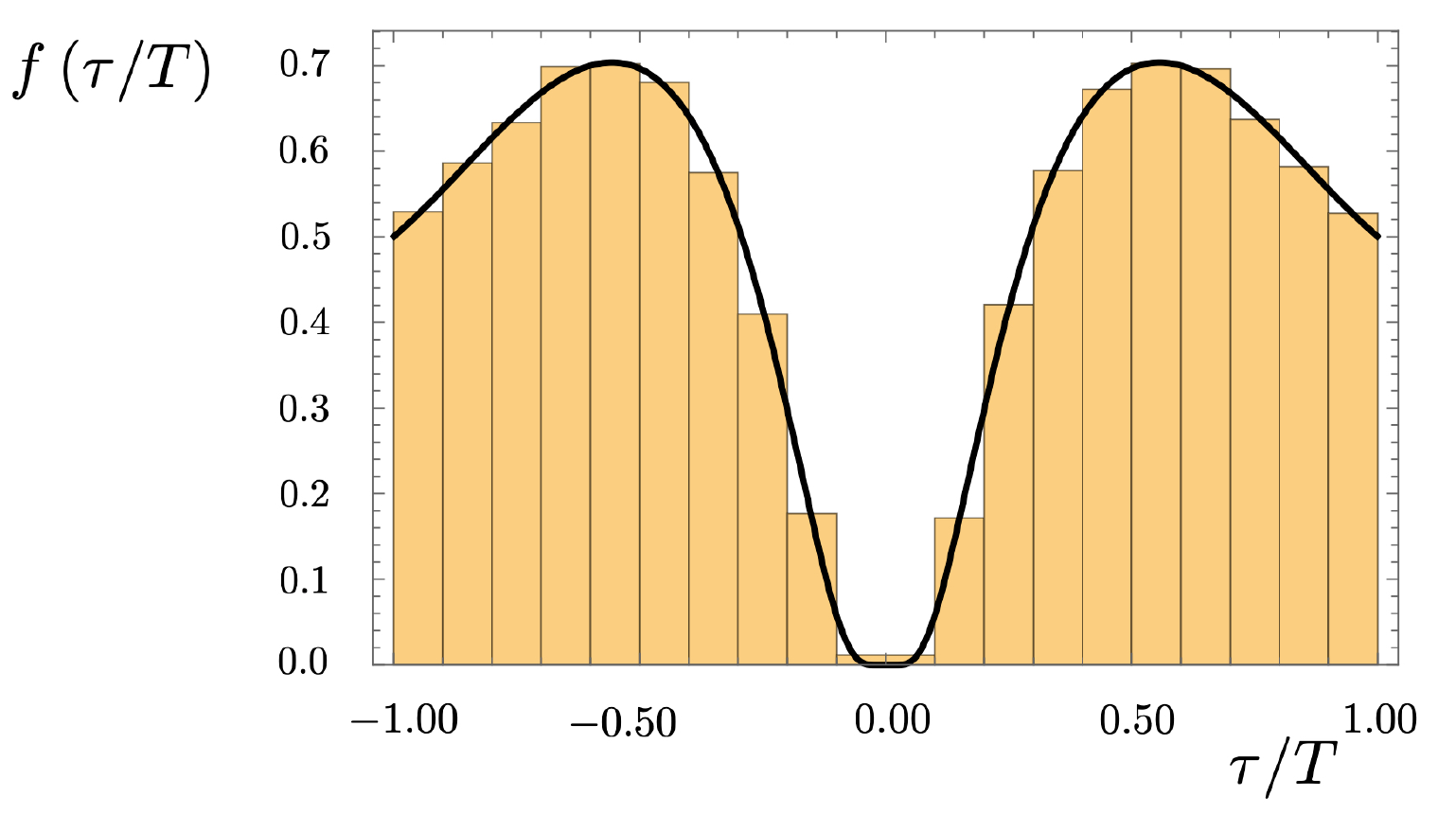}
\caption{The scaling function $f(\tau/T)$ given in eq.(\ref{majumdarform}) (black solid line) together with the experimental distribution of the $\tau/T$ variable. }
\label{taumajumdar}
\end{figure}
\subsection{Probability distribution of ${\bf \tau}$}
Recently Mori et al. \cite{MoriMajumdar} have determined the exact probability distribution of the time interval $\tau = t_{max} - t_{min}$. Denoting this function as $P(\tau | T)$, it has the scaling form 
\beq
P(\tau | T) \,=\,\frac{1}{T} \,f\left(\frac{\tau}{T}\right) \,\,\,,
\label{targettt}
\eeq
where the exact expression of the scaling function $f(x)$ is given by 
\beq
f(x) \,=\,\frac{2 (1- |x|)}{x^2} \, \sum_{m=-\infty}^{\infty} \frac{2 m +1}{\sinh\left(\left(2 m+1\right) \pi \sqrt{\frac{1-|x|}{|x|}}\right)} \,\,\,.
\label{majumdarform}
\eeq  
As shown in \cite{MoriMajumdar}, this scaling function satisfies the integral equation 
\beq
\int_0^1 \frac{f(y)}{1+ u y} \,=\,\int_0^\infty \frac{1}{\sinh z} \,\tanh^2\left(\frac{z}{2 \sqrt{1+u}}\right)\,\,\,,
\eeq
\begin{table}[b]
\begin{center}
\begin{tabular}{|c|c|c|c|}\hline
 \, $\langle |\tau|^n\rangle$  &  
  \, theoretical value & experimental value & relative error \\\hline
$\langle |\tau|^1\rangle$  & $ 0.5908$ & $0.5909$ & $0.13\times 10^{-3}$\\ 
$\langle |\tau|^2\rangle$  & $0.4009 $ & $0.4009$ & $ 0.03 \times 10^{-3}$\\
$\langle |\tau|^3\rangle$  & $  0.2972 $ & $0.2972 $ & $0.15 \times 10^{-3}$\\
$\langle |\tau|^4\rangle$  &  $0.2339$ & $0.2338$ & $ 0.38\times 10^{-3}$\\
$\langle |\tau|^5\rangle$&   $0.1918$ & $0.1917$ & $ 0.62 \times 10^{-3}$  \\
$\langle |\tau|^6\rangle$ &   $ 0.1621$ & $0.1620$ & $0.88 \times 10^{-3}$\\
$\langle |\tau|^7\rangle$   & $0.1401$ & $0.1400$& $ 1.15 \times 10^{-3}$\\ 
$\langle |\tau|^8\rangle$  & $0.1233$ & $0.1231$ & $ 1.43\times 10^{-3}$\\ 
$\langle |\tau|^9\rangle$ &  $0.1100$ & $0.1098$ & $ 1.74 \times 10^{-3} $\\
$\langle |\tau|^{10}\rangle$ &  $ 0.0992$ & $0.099$ & $ 2.08 \times 10^{-3}$\\
\hline
\end{tabular}
\end{center}
\caption{
 Moments of the probability distribution of the time between the maximum and minimum of Mertens function.} 
\label{momentsdistr}
\end{table}
\noindent
which is extremely useful for computing the various moments of the scaling variable $x = \tau/T$. Indeed, one has simply to expand both terms in powers of $u$ and compare the corresponding coefficients. In this way, the first exact expressions of the moments are 
\begin{eqnarray}
\langle | x | \rangle \,&=\,& \frac{4 \log(2) -1}{3} \,=\,0.5908.. \nonumber \\ 
\langle | x |^2 \rangle \,&=\,& \frac{7 \zeta(3) -2}{16} \,=\,0.4009..  \\ 
\langle | x |^3 \rangle \,&=\,& \frac{147 \zeta(3) -34}{480} \,=\,0.2972.. \nonumber \\ 
\langle | x |^4 \rangle \,&=\,& \frac{1701 \zeta(3) - 930 \zeta(5) -182}{3840} \,=\,0.2339 \nonumber 
\end{eqnarray}
where $\zeta(z)$ are the Riemann zeta function. The numerical value of other moments are reported in Table \ref{momentsdistr}. Also in this case 
we have sampled many different intervals of the sequence $\mathcal{S}_n$ and compared the experimental distribution of $\tau/T$ versus the theoretical one given in (\ref{majumdarform}). Here we show the result relative to an interval of $10^8$ values around $x=10^{15}$: we have divided this interval into $N=2\times 10^4$ segments of length $T = 5\times 10^3$ and for each of these intervals we have determined $\tau$. The comparison between the theoretical and experimental distributions are shown in Figure \ref{taumajumdar}: the $\chi^2$ variable of this example is $\chi^2 = 24.75$ (with $20$ degrees of freedom) and the relative $P$-value is $P$-value equals $0.73$, confirming the matching between the distribution coming from random walk and the Mertens function. Also in this case we have performed the Kolgomorv-Smirnov test to check the goodness of the expected function (\ref{targettt}) versus the data: with a sample of $n$ points, the Kolmogorov-Smirnov variable is defined as usual as 
\beq
D_n = sup_x |U_n - U(x)|\,\,\,,
\label{KSSS2}
\eeq
where, in this case, $U(x)$ is the cumulative distribution function associated to (\ref{targettt}) while $U_n$ is the ``experimental" cumulative distribution function with $n$ points. In our case, $n=20$ and we got $D_{20} = 0.00386475$ which gives a probability larger than $99.9\%$ that the experimental data follow indeed the expected theoretical function (\ref{targettt}). There is also an impressive agreement between the theoretical and the experimental values of the first ten moments as shown in Table \ref{momentsdistr}.

\def\C{C}
\def\Mhat{\hat{M}}
\def\U{\Mhat}
\def\s{\ell}

\section{Central limit theorem and a probabilistic conclusion on Mertens function}
In light of all the successful outputs of all the statistical tests that we have presented so far, it must appear natural that the Mertens function $\hat M(n)$ has a probabilistic 
distribution given by the gaussian normal law. 

As stressed several times, $\hat M(n)$ is completely deterministic with no probabilistic aspect per se. However, as explained in Section \ref{StatisticalSuite}, one can generate an ensemble and study its probability distribution by taking ``stroboscopic images" of this quantity. To this aim, we need to adopt the block variables $B_L(l)$ previously defined, see eq.\,(\ref{groupvariables}): these block variables are parts of the Mertens function and provide stroboscopic images of this function.  These quantities are the analog of the stroboscopic images which are needed to address the single Brownian motion problem \cite{Yuval,Mazo,Rudnick,Levy,diffusion}. 

The original interval $(1,n)$ relative to the start and ending points of the series $\hat M(n)$ can be broken into a large set of $N$ non-overlapping and well separated intervals of length $L$, of which we compute the corresponding block variables $B_L((l)$: it is this collection of block variables that forms the {\em set of events}, i.e. the ensemble $\CE$ relative to the sums of $L$ consecutive terms $\muhat(n)$
(see Section \ref{eenseble})
\beq
\label{CEdeff}
\CE = \{ B_L (\s) \},  ~~~~{\rm with~~}  I_L(\s)   \in {\cal G}_N
\eeq
Of course, given the infinite sequence of the M\"{o}bius coefficients, we could systematically enlarge the interval $n \rightarrow n'$, with $n' \gg n$, and take as new block variables those with length $n$ equal to the length of the original series $\hat M(n)$, and keep going with $n' \rightarrow n''$ and so on and so forth. Hence, with the ensemble defined as in eq.\,(\ref{CEdef}), i.e. made of non-overlapping and well separated intervals of length $L$, we compute the corresponding block variables of these intervals and sample them.   
Based on the exact result of vanishing mean of the sequence $\{\hat\mu_n\}$ and the absence of correlation between the restricted M\"{o}bius coeffients coming from the statistical analysis done in the previous sections 
\beq
\langle \hat\mu(n) \rangle \,=\,0 
\hspace{3mm}
,
\hspace{3mm}
\langle \hat\mu(n) \hat\mu(m) \rangle \,=\, \delta_{n,m}
\label{wellaverages}
\eeq
we expect that the probability distribution $P(z)$ of the random quantity
\beq
z \,=\, \frac{B_L}{\sqrt{L}}
\eeq
associated to the block variables is just given  by the normal distribution
\beq
P(z) \rightarrow \CN_{0,1} (z) \,\,\,.
\label{Unormaldist}
\eeq
This is clearly seen in Figure \ref{UNormal}, where we show a histogram relative to one of the samplings described above vs the normal distribution law $\CN_{0,1}(z)$.
The conclusion that the sampling of the restricted Mertens function satisfied the normal distribution is well supported by a large set of appropriate P-values of standard numerical analysis tests, such as Anderson-Darling, Cramer-von Mises, Pearson $\chi^2$, etc. We have also performed the Kolgomorov-Smirnov test relative to $n=30$ sampling point and we got $D_n= 0.0021$ for $n=30$ (see eq. (\ref{KSSS}) for the definition of $D_n$), which gives a probability larger than $99.9\%$ than the data follow indeed a normal law distribution.  
\begin{figure}[t]
\centering\includegraphics[width=.9\textwidth]{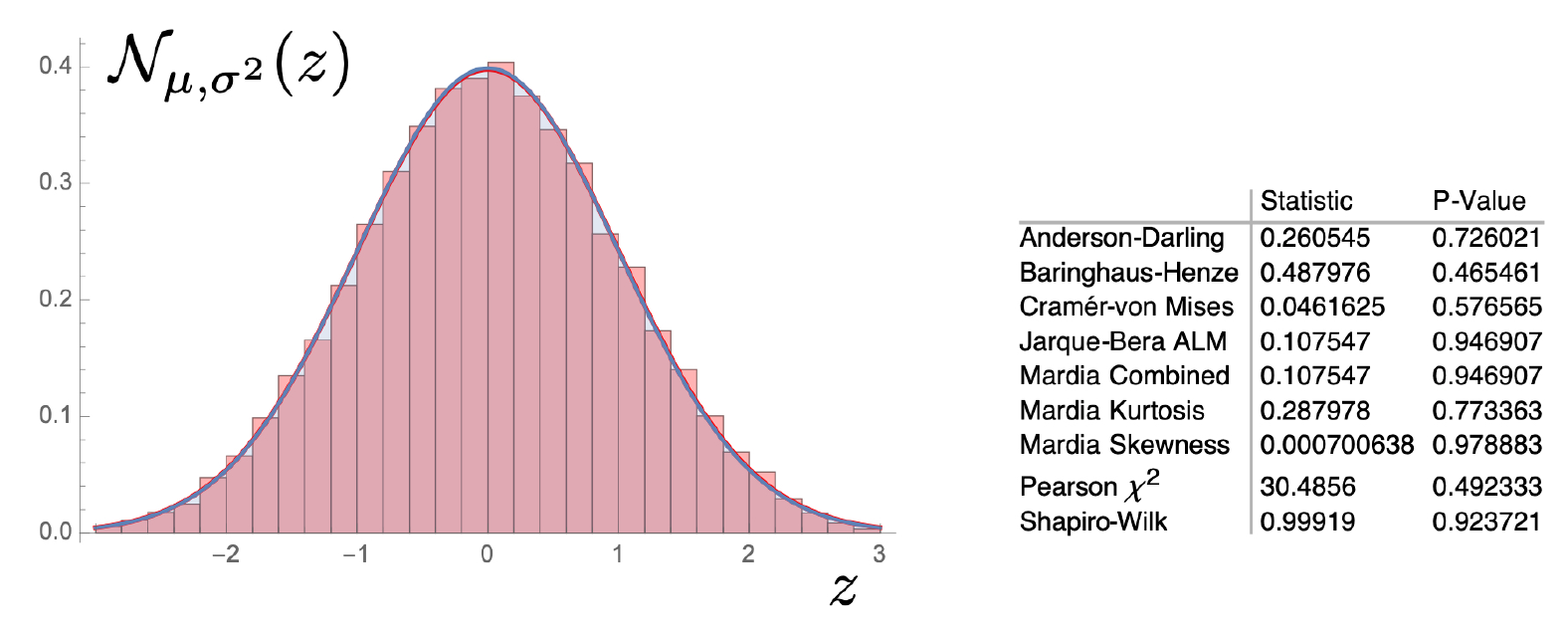}
\caption{Numerical evidence for the normal distribution proposed in eq.\eqref{Unormaldist}.     What is shown is a histogram of the properly normalized block variables $z=B_L/\sqrt{L}$. The ensemble $\CE$ corresponds to block variables of length $L= 10,000$, with $100.000$ intervals randomly separated in average $D = 1000$ generated starting from 
$\s = 10^{14}$. The blue curve is the normal distribution $\CN_{0,1} (z)$. 
On the right hand side, the table of P-values of confidence of the data to the normal distribution.}
\label{UNormal}
\end{figure}

It is well understood that the above gaussian distribution does {\it not} imply that $B_L = O(\sqrt{L})$.  In fact, the random walk behavior behind the gaussian distribution is known to imply
\beq
\label{limsup0}
\limsupL ~ \frac{|B_L|}{\sqrt{L}}  = \infty\,\,\,.
\eeq
For the remainder of this section, we describe two ways to deal with this fact as far as attempting to address the GRH based on \eqref{Unormaldist}.

First, based on the above normal distribution, the block variables $B_L$ of length $L$ always scale as 
$B_L = O(L^{1/2+\epsilon})$ 
for arbitrarily small $\epsilon >0$.  
In  probabilistic language, in the limit $L \rightarrow \infty$ this behavior occurs  with probability equal to 1. Indeed, using the gaussian  distribution (\ref{Unormaldist}), in the limit $L\rightarrow\infty$ we have  
 \barray
 \nonumber 
 \prob \probbracL  |B_L |    <   d  \, \sqrt{L} \probbracR  &=&
  \inv{\sqrt{2 \pi}} \int_{-d}^d  dx\,  e^{-x^2/2}  
 \\ 
 \label{PrBN} 
 &=& 1 -  \frac{e^{- d^2 /2}} {\sqrt{ 2 \pi } } \(  \frac{2}{d}  +  O\(  \inv{d^2 }  \)  \).
 \earray
  Chose $d  = \kappa L^\epsilon$ for any $\kappa >0$.   
 Then for any $\epsilon > 0 $,   
 \beq
 \label{ProbOne}
 \lim_{L \to \infty}  \prob \probbracL   B_L  = O(L^{1/2+ \epsilon}) \probbracR  = 1    \,\,\,.
 \eeq
 It is important to note that the above probability is not equal to $1$  as $L \to \infty$ if $\epsilon =0$.   
Since every element in the set ${\cal G}_N$  behaves the same as far as its growth is concerned, the almost surely true  probabilistic  statement  above, meaning with probability equal to $1$,   is promoted  to the surely true statement that
 \beq
 \label{Oep} 
 B_L = O(L^{1/2 + \epsilon}) ~~~{\rm  for~ any}   ~ \epsilon >0.     
 \eeq
 without assuming the RH.  

A second and more detailed  approach is also based on the above results.   Based on the above normal law,   one can place a different kind of bound not involving the above $\epsilon$,  by invoking  the law of iterated logs. If, as argued above, {$B_L (\ell)$}  is a random walk,  then the law states 
\beq
\label{iterated}
\limsupL ~ \frac{ |B_L|}{ \sqrt{2 L \log \log L} }   =  1    ~~~~~a. s. 
\eeq
where $a.\,  s.$ means almost surely.   
As stated above, the bound \eqref{iterated} is sufficient for showing the validity of the GRH since the $\log$'s do not spoil convergence of the integral \eqref{Mellin}. 
The analysis done for the block variables $B_L$ which sample on  all possible scales the restricted Mertens function $\hat M(n)$ extends of course to $\hat M(n)$ itself, so that we can conclude 
\beq
\label{iteratedmmm}
\limsupn ~ \frac{ |\hat M(n)|}{ \sqrt{2 n \log \log n} }   =  1    ~~~~~a. s. 
\eeq

\newpage

\addcontentsline{toc}{section}{Part D}
\begin{center}
{\Large {\bf PART D}}
\end{center}

\section{Conclusions}

There are different ways to approach the Riemann Hypothesis, as discussed for instance in \cite{Borwein}. 
 In this paper, we have pursued an approach which consists of identifying the first singularity,  {approaching the critical line  ${\mathbb Re} \,s =1/2$ from the right,}  of the inverse Mellin transform of the restricted Mertens function $\hat M(n)$ for  the square-free numbers $\sqf_n$.  Such a function $\hat M(n)$ is a sum of $\pm 1$ terms relative to the M\"{o}bius coefficients $\hat \mu(n)$ of the square-free numbers. Our efforts in this paper have been to show the stochastic behaviour of the restricted Mertens function $\hat M(n)$, despite its deterministic nature. To this aim, we have developed a global approach to the function $\hat M(n)$ based on a series of probabilistic results concerning the prime number distribution along the series of square-free numbers, the average number of prime divisors, the Erd\H{o}s-Kac theorem for square-free numbers, etc. All these results support  the existence of a  normal  law  satisfied by $\hat M(n)$. 
 

It is better to stress once again that there is  nothing which prevents a deterministic sequence from being a realization of a random process,  and the restricted M\"{o}bius coefficients and the associated Mertens function seem to be a very good example of this fact. This has been strikingly confirmed by our further numerical analysis aimed to check the local stochasticity of the restricted M\"{o}bius coefficients. Indeed, driven by the natural curiosity to see how good the  random nature of the Mertens function $\hat M(n)$ and its M\"{o}bius coefficients $\hat\mu(n)$ is satisfied by finite samples of these sequence, we have performed a massive statistical analysis in a range up to $10^{16}$ natural numbers, applying numerous randomness tests of increasing precision and complexity: {\em all} these tests were successfully satisfied with a level of reliability of $99\%$, as quantified by the chi-square distributions. 

For a theoretical physicist, all these results leave essentially no room to doubt the validity of the Riemann Hypothesis: 
the probabilistic approach pursued here has led to the conclusion that 
 \beq
\limsupn  ~\frac{| \hat M(n)|}{\sqrt{2 n \log\log n}} \,= \, 1 
\hspace{3mm} 
, 
\hspace{3mm}
{\rm a. s. } 
\label{sqrtloglog555}
\eeq
This is a probabilistic result that applies to all random sequences made of $\pm 1$, therefore also to the 
Mertens function $\hat M(n)$ once its random nature has been established.  We expect that our bound  (\ref{sqrtloglog555}) is the best one can obtain  for the Mertens function 
{\em without assuming the RH}: as a matter of fact, in this paper we have reversed the argument, i.e. we have shown that we can use the bound (\ref{sqrtloglog555}) to arrive to the validity of the RH. The successful outputs of all the statistical tests performed in this paper show that the sequence of the M\"{o}bius coefficients $\hat \mu(n)$, restricted to the square-free numbers, is, for all practical purposes, a random sequence and it can be used as a perfect random number generator in any computer simulation, an interesting by-product of our analysis. 

However, putting aside any probabilistic considerations, one may be curious about the asymptotic behavior of the actual Mertens function $M(x)$. As easily imagined, this behavior  is presently unknown but an estimate of the asymptotic behaviour of this function can be achieved {\em assuming the validity of the Riemann Hypothesis}. This is, for instance, the content of a series of bounds obtained by various authors. Let's define 
\begin{eqnarray}
&&  {\mathcal M}_S(x)  \equiv \sqrt{x} \, \exp\left[(\log x)^{1/2} \, (\log\log x)^{5/2+\epsilon}\right] \,\nonumber\\
&& {\mathcal M}_R (x) \equiv (2 x  \log \log x)^{1/2} \\
&& {\mathcal M}_G(x) \equiv  \sqrt{x} \,(\log \log \log x)^{5/4}  \nonumber 
\end{eqnarray}
Soundararajan, for instance, assuming the RH, arrived to the bound \cite{SounMert,Balazard}
\beq
M(x) \ll {\mathcal M}_S(x)  \,\,\,.
\eeq
A more stringent bound comes from an unpublished conjecture by Gonek \cite{Gonek1}, who stated that 
\beq
\label{GonekConj}
M(x) = O ({\mathcal M}_G(x)) \,\,\,.
\eeq
The  above bound remains an unpublished conjecture, however further evidence for it was given by Ng \cite{Ng}, and it is likely to be correct.       
Notice that, if it is indeed correct, it does not contradict the bound (\ref{sqrtloglog555}) which was based on the normal distribution we described in this paper. As a matter of fact, 
the probability estimate ${\mathcal M}_R(x)$ of the (restricted) Mertens function  based on the law of iterated logs described in this paper satisfies the inequalities 
\beq
 {\mathcal M}_G(x) <  {\mathcal M}_R(x) <  {\mathcal M}_S(x) \,,
\eeq
for large $x$, 
and it is therefore in agreement with both bounds. 

In this paper we have also argued that the validity of the RH can imply the validity of the GRH, and we have argued this on the basis of two hypotheses: (i) that the restricted Mertens function goes as $|\hat M(x) |\sim x^{1/2 + \epsilon}$ and (ii) that the restricted M\"{o}bius coefficients $\hat \mu(n)$ behave as random independent variables. Both hypotheses 
seem indeed satisfied. It would be of course interesting to find an alternative proof of the inequality (\ref{Mqqr}) based on a global analysis of the functions $\hat M_{r}(x)$, similarly to the analysis done in Part B for the function $\hat M(x)$. 

In conclusion, using a probabilistic approach and the robustness of probabilistic theorems that followed, the validity of the Riemann and the Generalised Riemann Hypothesis comes quite naturally from the central limit theorem and the dominant role of the gaussian distribution law in the probabilistic realms. Even taking a more conservative point of view, 
it remains however true that, in view of the large battery of analytic results and statistical tests here presented, while a violation of the RH may be still possible, it is an event 
which appears however to be {wildly} improbable. 


\vspace{10mm}
\noindent
{\bf {\Large Acknowledgements}}

\vspace{3mm}
\noindent
We would like to thank Andrea Gambassi, Steve Gonek, Satya Majumdar, Silvia Pappalardi, Kannan Soundararajan and Don Zagier for useful discussions on different aspects of the problem. We are particularly grateful to German Sierra for his suggestions and a thorough and deep reading of the manuscript. 
This work has been carried on for several years, during which we have enjoyed the hospitality of several institutes. GM would like to thank in particular: the Simons Center for Geometry and Physics in Stony Brook, the International Institute of Physics in Natal and Institute Henri Poincare in Paris for the nice hospitality and partial support during part of this work. AL thank SISSA for the kind hospitality and partial support during the early stage of the work. 

\newpage

\appendix

\section{The mean of the Mertens function}\label{AA}

Let's start from an identity which involves the M\"{o}bius coefficients $\mu(n)$, namely \cite{Apostol}
\beq
\sum_{d | m} \mu(d) \,=\,\delta_{m,1} \,\,\,, 
\label{summatory}
\eeq
where the sum over $d$ is made on the numbers $d$ which divide $m$. This identity is important for proving that the mean of the Mertens function vanishes. 
\beq
\lim_{x\rightarrow \infty} \frac{1}{x} \sum_{m=1}^x \mu(m) \,=\, 0 \,\,\,. 
\label{averageMertens11}
\eeq
To prove this result, let's start by showing that  for any $x \geq 1$ it holds 
\beq
\left| \sum_{m=1}^x \frac{\mu(m)}{m} \right| \,\leq\,  1\,\,\,.
\label{inequality111}
\eeq
Consider the sum $S(N) \,=\,\sum_{m \leq N} \delta_{m,1}$, which is of course equal to $1$.  Thanks to the identity (\ref{summatory}), we can write $S(N)$ as 
\begin{eqnarray}
S(N) &\,=\, & \sum_{m \leq N} \sum_{d | m} \mu(d)  \,=\, \sum_{d \leq N} \mu(d) \left[\frac{N}{d}\right] \, =\, \\
& \, = \, & N \,\sum_{d \leq N} \frac{\mu(d)}{d} - \sum_{d \leq N} \mu(d) \,\left\{\frac{N}{d}\right\} \nonumber\,\,\,,
\end{eqnarray}
where $[ a ]$ and $\{a\} $ are respectively the integer and the fractional parts of the number $a$. Using that $\{N/d\} =0$ if $d=N$ and the bound $|\mu(d) \{N/d\}| < 1$ for all 
$1\leq d < N$, for the last term we have $ |\sum_{d \leq N} \mu(d) \,\left\{N/d\right\} | < (N-1)$ and therefore 
\beq
\left| N \sum_{d \leq N} \frac{\mu(d)}{d} \right| \leq (N-1) + | S(N) | \,=\,(N-1) + 1 \,=\, N\,\,\,.
\eeq
So, taking $x = N$, we have indeed the bound (\ref{inequality111}). Consider now the function 
\beq
A(x) \,=\, \sum_{m \leq x} \frac{\mu(m)}{m} \,\,\,,
\eeq
and write the Mertens function $M(x)$ as 
\begin{equation} 
M(x)\, =\, \sum_{m\leq x} \mu(m) \,=\, \sum_{m\leq x} \frac{\mu(m)}{m} \, n \,=\, x\, A(x) - \int_1^x A(t) \, dt \,\,\,,
\end{equation}
namely
\beq
\frac{M(x)}{x} \,=\,A(x) - \frac{1}{x} \int_1^x A(t) \, dt\,\,\,. 
\eeq
Taking $x\rightarrow \infty$, we have that $A(x) \rightarrow 0$ since it is just the value of the M\"{o}bius function $\tilde\mu(s)$ at $s=1$, while for the second term, given an arbitrary $\epsilon >0$, there exists a constant $c$ (which depends on $\epsilon$) such that $|A(x)| < \epsilon$ if $x > c$. Hence, we have 
\beq
\left|\frac{1}{x} \int_1^x A(t) dt \right| \leq \left|\frac{1}{x} \int_1^c A(t) \, dt\right| + \left|\frac{1}{x} \int_c^x A(t) dt \right| \,\leq \, \frac{c-1}{x} + \frac{\epsilon (x-c)}{x} \,\,\,, 
\eeq
 and thererfore, in the limit $x \rightarrow \infty$ 
 \beq
 \lim_{x\rightarrow \infty} {\rm sup}\, \left|\frac{1}{x} \int_1^x A(t) \, dt \right| \,\leq \, \epsilon \,\,\,.
 \eeq
 Since $\epsilon$ can be arbitrarily small, we have then the prove of (\ref{averageMertens11}).

\end{document}